\documentclass[oneside,english]{amsart}
\usepackage[T1]{fontenc}
\usepackage[latin1]{inputenc}
\usepackage{amssymb}

\makeatletter
 \theoremstyle{plain}
\newtheorem{thm}{Theorem}[section]
  \theoremstyle{plain}
  \newtheorem{lem}[thm]{Lemma}
  \theoremstyle{remark}
  \newtheorem*{rem*}{Remark}
  \theoremstyle{definition}
  \newtheorem{defn}[thm]{Definition}
  \theoremstyle{plain}
  \newtheorem{prop}[thm]{Proposition}
  \theoremstyle{plain}
  \newtheorem*{thm*}{Theorem}
  \theoremstyle{remark}
  \newtheorem{rem}[thm]{Remark}
  \theoremstyle{plain}
  \newtheorem{cor}[thm]{Corollary}

\usepackage{babel}
\makeatother
\begin{document}

\title{Plancherel-Rotach Asymptotics for $q$-Series }

\curraddr{School of Mathematics\\
Guangxi Normal University\\
Guilin City, Guangxi 541004\\
P. R. China.}



\begin{abstract}
In this work we study the Plancherel-Rotach type asymptotics for selected
$q$-series and $q$-orthogonal polynomials with complex scalings.
The $q$-series we cover are Euler's $q$-exponential, Ramanujan function,
Jackson's $q$-Bessel function of second kind, Ismail-Masson orthogonal
polynomials, Stieltjes-Wigert polynomials and $q$-Laguerre polynomials.
For a fixed $q$ with $0<q<1$, in each case the main term of the
asymptotic formulas may contain Ramanujan function or theta function
depending on the value of scaling parameter. Furthemore, when the
scaling parameter is in certain strip of the complex plane, its number
theoretical property completely determines the order of the error
term. In each cases, we also investigate the asymptotic behavior of
the mentioned $q$-series when $q$ approaching $1$ in a restricted
manner. These asymptotic formulas may provide insights to new random
matrix models.
\end{abstract}

\author{Ruiming Zhang}

\email{ruimingzhang@yahoo.com}

\maketitle

\section{Introduction\label{sec:Introduction}}

The Plancherel-Rotach asymptotics for classical orthogonal polynomials
are essential to obtain universality results in random matrix theory
\cite{Mehta}. The associated random matrix models for $q$-orthogonal
polynomials are still unknown today. It might be interesting to calculate
the Plancherel-Rotach type asymptotics for $q$-orthogonal polynomials
to gain some insights to the related random matrix models.

The Plancherel\emph{-}Rotach type asymptotics for $q$-orthogonal
polynomials started in \cite{Ismail7}, where we studied Plancherel\emph{-}Rotach
type asymptotics for three families of $q$-orthogonal polynomials
with real logarithmic scalings. They are Ismail-Masson polynomials
$\left\{ h_{n}(x|q)\right\} _{n=0}^{\infty}$ , Stieltjes-Wigert polynomials
$\left\{ S_{n}(x;q)\right\} _{n=0}^{\infty}$ and $q$-Laguerre polynomials
$\left\{ L_{n}^{(\alpha)}(x;q)\right\} _{n=0}^{\infty}$. These $q$-orthogonal
polynomials are associated with indeterminant moment problems \cite{Akhiezer}.
The asymptotics reveal a remarkable pattern which is quite different
to the pattern associated with classical Plancherel\emph{-}Rotach
asymptotics \cite{Szego,Ismail2,Deift1,Deift2}. The main term of
asymptotics may contain Ramanujan function $A_{q}(z)$ or theta function
according to the value of the scaling parameter. 

In this work we will investigate Plancherel\emph{-}Rotach type asymptotics
with complex scalings for some $q$-series, first with fixed $q\in(0,1)$,
then let $q$ approaching $1$ in some restricted manners. The $q$-series
treated in this work are Euler's $q$-exponential $E_{q}(z)$, Ramanujan
function $A_{q}(z)$, Jackson's $q$-Bessel functions of second kind
$J_{\nu}^{(2)}(z;q)$, Ismail-Masson orthogonal polynomials $\left\{ h_{n}(x|q)\right\} _{n=0}^{\infty}$,
Stieltjes-Wigert polynomials $\left\{ S_{n}(x;q)\right\} _{n=0}^{\infty}$
and $q$-Laguerre polynomials $\left\{ L_{n}^{(\alpha)}(x;q)\right\} _{n=0}^{\infty}$.
The results reveal an even more complicated pattern. The error terms
manifest a link between special functions and number theory. The method
we use in this study is a discrete version of the classical Laplace
method in asymptotics \cite{Ismail7}. This version of discrete Laplace
method may also be applied to study the asymptotics for other $q$-series
such as general confluent basic hypergeometric functions. We prefer
to publish these results at another time. 

Through out this work, We shall assume that $0<q<1$ unless otherwise
stated. We also assume that $s=\sigma+it$ is a complex number with
$\sigma=\Re(s)$ and $t=\Im(s)$. All the $\log$ and power functions
are taken as their principle branches. The asymptotic formulas here
are mostly for conceptual purposes, so we won't pursue the best orders
nor the best implied constants.

\subsection{$q$-series}

For any complex number $a$, we define \cite{Andrews4,Gasper,Ismail2,Koekoek}\begin{equation}
(a;q)_{\infty}:=\prod_{k=0}^{\infty}(1-aq^{k})\label{eq:1}\end{equation}
 and the $q$-shifted factorial of $a,a_{1},\dotsc a_{m}$ are defined
by

\begin{equation}
(a;q)_{n}:=\frac{(a;q)_{\infty}}{(aq^{n};q)_{\infty}},\quad(a_{1},\dotsc,a_{m};q)_{n}:=\prod_{k=1}^{m}(a_{k};q)_{n}\label{eq:2}\end{equation}
 for all integers $n\in\mathbb{Z}$ amd $m\in\mathbb{N}$. 

\begin{lem}
\label{lem:1}Given any complex number $z$, we assume that\begin{equation}
0<\frac{\left|z\right|q^{n}}{1-q}<\frac{1}{2}\label{eq:3}\end{equation}
 for some positive integer $n$. Then, \begin{equation}
\frac{(z;q)_{\infty}}{(z;q)_{n}}=(zq^{n};q)_{\infty}:=1+r_{1}(z;n),\label{eq:4}\end{equation}
 with\begin{equation}
\left|r_{1}(z;n)\right|\le\frac{2\left|z\right|q^{n}}{1-q},\label{eq:5}\end{equation}
 and \begin{equation}
\frac{(z;q)_{n}}{(z;q)_{\infty}}=\frac{1}{(zq^{n};q)_{\infty}}:=1+r_{2}(z;n),\label{eq:6}\end{equation}
 with \begin{equation}
\left|r_{2}(z;n)\right|\le\frac{2\left|z\right|q^{n}}{(1-q)}.\label{eq:7}\end{equation}

\end{lem}
\begin{proof}
From the $q$-binomial theorem \cite{Andrews4,Gasper,Ismail2,Koekoek}\begin{equation}
\frac{(Az;q)_{\infty}}{(z;q)_{\infty}}=\sum_{k=0}^{\infty}\frac{(A;q)_{k}}{(q;q)_{k}}z^{k},\quad A,z\in\mathbb{C},\quad|z|<1,\label{eq:8}\end{equation}
and the inequality \begin{equation}
(q;q)_{k}\ge(1-q)^{k}\label{eq:9}\end{equation}
 for $k=0,1,...$ we obtain \begin{eqnarray}
r_{2}(z;n) & = & \sum_{k=0}^{\infty}\frac{\left(zq^{n}\right)^{k+1}}{(q;q)_{k+1}},\label{eq:10}\end{eqnarray}
and \begin{align}
\left|r_{2}(z;n)\right| & \le\sum_{k=0}^{\infty}\frac{\left(\left|z\right|q^{n}\right)^{k+1}}{(q;q)_{k+1}}\le\frac{\left|z\right|q^{n}}{(1-q)}\sum_{k=0}^{\infty}\left(\frac{\left|z\right|q^{n}}{1-q}\right)^{k}\le\frac{2\left|z\right|q^{n}}{(1-q)}.\label{eq:11}\end{align}
 Apply a limiting case of \eqref{eq:8}, \begin{equation}
(z;q)_{\infty}=\sum_{k=0}^{\infty}\frac{q^{k(k-1)/2}}{(q;q)_{k}}(-z)^{k}\quad z\in\mathbb{C},\label{eq:12}\end{equation}
 and the inequalities,\begin{equation}
\frac{1-q^{k}}{1-q}\ge kq^{k-1},\quad\frac{(q;q)_{k}}{(1-q)^{k}}\ge k!q^{k(k-1)/2}\label{eq:13}\end{equation}
 for $k=0,1,\dotsc$ we get\begin{align}
r_{1}(z;n) & =\sum_{k=1}^{\infty}\frac{q^{k(k-1)/2}(-zq^{n})^{k}}{(q;q)_{k}},\label{eq:14}\end{align}
 and \begin{align}
\left|r_{1}(z;n)\right| & \le\sum_{k=0}^{\infty}\frac{(\left|z\right|q^{n})^{k+1}}{(1-q)^{k+1}}\frac{(1-q)^{k+1}q^{k(k+1)/2}}{(q;q)_{k+1}}\label{eq:15}\\
 & \le\sum_{k=0}^{\infty}\frac{(\left|z\right|q^{n})^{k+1}}{(1-q)^{k+1}}\frac{1}{(k+1)!}\le\frac{\left|z\right|q^{n}\sqrt{e}}{1-q}<\frac{2\left|z\right|q^{n}}{1-q}.\nonumber \end{align}

\end{proof}
\begin{rem*}
In this work, we either have a fixed $q$ with $0<q<1$, $q=\exp(-d\pi n^{-a})$
or $q=\exp\left(-\frac{d}{\log n}\right)$, with $0<a<1$ and $d>0$.
In these cases, the condition \eqref{eq:3} is clearly satisfied for
$n$ large. 
\end{rem*}
The Dedekind $\eta(\tau)$ is defined as \cite{Rademarcher}\begin{equation}
\eta(\tau):=e^{\pi i\tau/12}\prod_{k=1}^{\infty}(1-e^{2\pi ik\tau}),\label{eq:16}\end{equation}
 or\begin{equation}
\eta(\tau)=q^{1/12}(q^{2};q^{2})_{\infty},\label{eq:17}\end{equation}
 for \begin{equation}
q=e^{\pi i\tau},\quad\Im(\tau)>0.\label{eq:18}\end{equation}
 It has the transformation formula\begin{equation}
\eta\left(-\frac{1}{\tau}\right)=\sqrt{\frac{\tau}{i}}\eta(\tau).\label{eq:19}\end{equation}
 The four Jacobi theta functions are defined as \begin{align}
\theta_{1}(v|\tau) & :=-i\sum_{k=-\infty}^{\infty}(-1)^{k}q^{(k+1/2)^{2}}e^{(2k+1)\pi iv},\label{eq:20}\\
\theta_{2}(v|\tau) & :=\sum_{k=-\infty}^{\infty}q^{(k+1/2)^{2}}e^{(2k+1)\pi iv},\label{eq:21}\\
\theta_{3}(v|\tau) & :=\sum_{k=-\infty}^{\infty}q^{k^{2}}e^{2k\pi iv},\label{eq:22}\\
\theta_{4}(v|\tau) & :=\sum_{k=-\infty}^{\infty}(-1)^{k}q^{k^{2}}e^{2k\pi iv},\label{eq:23}\end{align}
 where\begin{equation}
q=e^{\pi i\tau},\quad\Im(\tau)>0.\label{eq:24}\end{equation}
 For our convenience, we also use the following notations\begin{equation}
\theta_{\lambda}(z;q):=\theta_{\lambda}(v|\tau),\quad z=e^{2\pi iv},\quad q=e^{\pi i\tau}\label{eq:25}\end{equation}
for $\lambda=1,2,3,4$. The Jacobi's triple product identities are\begin{align}
\theta_{1}(v|\tau) & =2q^{1/4}\sin\pi v(q^{2};q^{2})_{\infty}(q^{2}e^{2\pi iv};q^{2})_{\infty}(q^{2}e^{-2\pi iv};q^{2})_{\infty},\label{eq:26}\\
\theta_{2}(v|\tau) & =2q^{1/4}\cos\pi v(q^{2};q^{2})_{\infty}(-q^{2}e^{2\pi iv};q^{2})_{\infty}(-q^{2}e^{-2\pi iv};q^{2})_{\infty},\label{eq:27}\\
\theta_{3}(v|\tau) & =(q^{2};q^{2})_{\infty}(-qe^{2\pi iv};q^{2})_{\infty}(-qe^{-2\pi iv};q^{2})_{\infty},\label{eq:28}\\
\theta_{4}(v|\tau) & =(q^{2};q^{2})_{\infty}(qe^{2\pi iv};q^{2})_{\infty}(qe^{-2\pi iv};q^{2})_{\infty}.\label{eq:29}\end{align}
 The Jacobi $\theta$ functions satisfy transformations\begin{align}
\theta_{1}\left(\frac{v}{\tau}\mid-\frac{1}{\tau}\right) & =-i\sqrt{\frac{\tau}{i}}e^{\pi iv^{2}/\tau}\theta_{1}\left(v\mid\tau\right),\label{eq:30}\\
\theta_{2}\left(\frac{v}{\tau}\mid-\frac{1}{\tau}\right) & =\sqrt{\frac{\tau}{i}}e^{\pi iv^{2}/\tau}\theta_{4}\left(v\mid\tau\right),\label{eq:31}\\
\theta_{3}\left(\frac{v}{\tau}\mid-\frac{1}{\tau}\right) & =\sqrt{\frac{\tau}{i}}e^{\pi iv^{2}/\tau}\theta_{3}\left(v\mid\tau\right),\label{eq:32}\\
\theta_{4}\left(\frac{v}{\tau}\mid-\frac{1}{\tau}\right) & =\sqrt{\frac{\tau}{i}}e^{\pi iv^{2}/\tau}\theta_{2}\left(v\mid\tau\right).\label{eq:33}\end{align}

\begin{lem}
\label{lem:2}Let \begin{equation}
0<a<1,\quad n\in\mathbb{N},\label{eq:34}\end{equation}
 and\begin{equation}
q=e^{-2\pi\gamma^{-1}n^{-a}},\quad\gamma>0,\label{eq:35}\end{equation}
 then,\begin{equation}
(q;q)_{\infty}=\sqrt{\gamma n^{a}}\exp\left(\pi n^{-a}/(12\gamma)-\pi\gamma n^{a}/12\right)\left\{ 1+\mathcal{O}\left(e^{-2\pi\gamma n^{a}}\right)\right\} ,\label{eq:36}\end{equation}
 and\begin{equation}
\frac{1}{(q;q)_{\infty}}=\frac{\exp\left(\pi\gamma n^{a}/12-\pi n^{-a}/(12\gamma)\right)}{\sqrt{\gamma n^{a}}}\left\{ 1+\mathcal{O}\left(e^{-2\pi\gamma n^{a}}\right)\right\} \label{eq:37}\end{equation}
 as $n\to\infty$. 
\end{lem}
\begin{proof}
From formulas \eqref{eq:15} and \eqref{eq:18} we get \begin{align}
 & (q;q)_{\infty}=\exp\left(\pi n^{-a}/(12\gamma)\right)\eta\left(n^{-a}\gamma^{-1}i\right)\label{eq:38}\\
 & =\sqrt{\gamma n^{a}}\exp\left(\pi n^{-a}/(12\gamma)\right)\eta(\gamma n^{a}i)\nonumber \\
 & =\sqrt{\gamma n^{a}}\exp\left(\pi n^{-a}/(12\gamma)-\pi\gamma n^{a}/12\right)\prod_{k=1}^{\infty}(1-e^{-2\pi\gamma kn^{a}})\nonumber \\
 & =\sqrt{\gamma n^{a}}\exp\left(\pi n^{-a}/(12\gamma)-\pi\gamma n^{a}/12\right)\left\{ 1+\mathcal{O}\left(e^{-2\pi\gamma n^{a}}\right)\right\} \nonumber \end{align}
 and\begin{equation}
\frac{1}{(q;q)_{\infty}}=\frac{\exp\left(\pi\gamma n^{a}/12-\pi n^{-a}/(12\gamma)\right)}{\sqrt{\gamma n^{a}}}\left\{ 1+\mathcal{O}\left(e^{-2\pi\gamma n^{a}}\right)\right\} \label{eq:39}\end{equation}
 as $n\to\infty$. 
\end{proof}
Similarly, we have:

\begin{lem}
\label{lem:3}Assume that \begin{equation}
q=\exp\left(-\frac{2\pi}{\gamma\log n}\right),\quad n\ge2,\quad\gamma>0,\label{eq:40}\end{equation}
 then,\begin{equation}
(q;q)_{\infty}=n^{-\pi\gamma/12}\sqrt{\gamma\log n}\exp\left(\frac{\pi}{12\gamma\log n}\right)\left\{ 1+\mathcal{O}\left(n^{-2\pi\gamma}\right)\right\} ,\label{eq:41}\end{equation}
 and\begin{equation}
\frac{1}{(q;q)_{\infty}}=\frac{n^{\pi\gamma/12}}{\sqrt{\gamma\log n}}\exp\left(-\frac{\pi}{12\gamma\log n}\right)\left\{ 1+\mathcal{O}\left(n^{-2\pi\gamma}\right)\right\} \label{eq:42}\end{equation}
 as $n\to\infty$. 
\end{lem}
We will also make use of the trivial inequality \begin{equation}
|e^{z}-1|\le|z|e^{|z|}\label{eq:43}\end{equation}
 for any $z\in\mathbb{C}$.

\subsection{Generalized Irrational Measure}

For an irrational number $\theta$, Chebyshev's theorem implies that
for any real number $\beta$ , there exist infinitely many pairs of
integers $n$ and $m$ with $n>0$ such that \cite{Hua} \begin{equation}
n\theta=m+\beta+\gamma_{n}\quad{\rm with}\quad|\gamma_{n}|\le\frac{3}{n}.\label{eq:44}\end{equation}
Clearly, Chebyshev's theorem says that the arithmetic progression
$\left\{ n\theta\right\} _{n\in\mathbb{Z}}$ is ergodic in $\mathbb{R}$. 

\begin{defn}
\label{def:generalized irrationality measure}Given real numbers $\theta_{j}$
and $\beta_{j}$ for $j=1,...,N$, the generalized \emph{}irrationality
\emph{}measure $\omega(\theta_{1},\dots,\theta_{N}|\beta_{1},\dots,\beta_{N})$
of $\theta_{1},\dots,\theta_{N}$ associated with $\beta_{1},\dots,\beta_{N}$
is defined as the least upper bound of the set of real numbers $r$
such that there exist infinitely many positive integers $n$ and integers
$m_{1},\dots,m_{N}$ such that\begin{equation}
|n\theta_{j}-\beta_{j}-m_{j}|<\frac{1}{n^{r-1}}.\label{eq:45}\end{equation}
 for $j=1,\dots,N$. 
\end{defn}
\begin{prop}
Let $\theta$ be an irrational number, then for any real number $\beta$,
its generalized irrational measure associated with $\beta$ is $\omega(\theta|\beta)\ge2$
. 
\end{prop}
\begin{proof}
The assertion is a direct consequence of the Chebyshev's theorem. 
\end{proof}
Recall that the irrationality \emph{}measure (or Liouville\emph{-}Roth
\emph{}constant) $\mu(\theta)$ of a real number $\theta$ is defined
as the least upper bound of the set of real numbers $r$ such that
\cite{Wikipdedia} \begin{equation}
0<|n\theta-m|<\frac{1}{n^{r-1}}\label{eq:46}\end{equation}
is satisfied by an infinite number of integer pairs $(n,m)$ with
$n>0$. It is clear that we have \begin{equation}
\omega(\theta|0)=\mu(\theta).\label{eq:47}\end{equation}
A real algebraic number $\theta$ of degree $\ell$ if it is a root
of an irreducible polynomial of degree $\ell$ with integer coefficients.
Liouville's theorem in number theory says that for a real algebraic
number $\theta$ of degree $\ell$, there exists a positive constant
$K(\theta)$ such that for any integer $m$ and $n>0$ we have \begin{equation}
|n\theta-m|>\frac{K(\theta)}{n^{\ell-1}}.\label{eq:48}\end{equation}
A Liouville number is a real number $\theta$ such that for any positive
integer $\ell$ there exist infinitely many integers $n$ and $m$
with $n>1$ such that \cite{Wikipdedia}\begin{equation}
0<|n\theta-m|<\frac{1}{n^{\ell-1}},\label{eq:49}\end{equation}
It is well known that even though the set of all Liouville numbers
is of Lebesgue measure zero, topologically, almost all real numbers
are Liouville numbers.

\begin{prop}
The generalized irrational measure $\omega(\theta|\beta)$ has the
following properties:
\begin{enumerate}
\item For any real algebraic number $\theta$ of degree $\ell$, one has
$\omega(\theta|0)\le\ell$.
\item For a Liouville number $\theta$, one has $\omega(\theta|0)=\infty$.
\end{enumerate}
\end{prop}
\begin{proof}
The first assertion follows from the definition \ref{def:generalized irrationality measure}
and \eqref{eq:47} while the last assertion follows directly from
the definitions of the generalized irrational measure and the Liouville
numbers.
\end{proof}
It is clear that the quadratic irrationals $\theta$ such as $\sqrt{2}$
have the generalized irrational measure $\omega(\theta|0)=2$ for
any real number $\beta$. 

For any real number $\theta$, we consider the set \begin{equation}
\mathbb{S}(\theta)=\left\{ \left\{ n\theta\right\} :n\in\mathbb{N}\right\} .\label{eq:50}\end{equation}
Obviously, the set $\mathbb{S}(\theta)$ is a finite subset of $[0,1)$
for $\theta\in\mathbb{Q}$. When $\theta\notin\mathbb{Q}$, the set
$\mathbb{S}(\theta)$ is dense in $(0,1)$, which is a consequence
of Chebyshev's theorem.

\begin{lem}
We have the following:
\begin{enumerate}
\item Assume that $\theta>0$ and $\theta\in\mathbb{Q}$. For any real number
$\lambda$ with $\left\{ \lambda\right\} \in\mathbb{S}(\theta)$,
there exist infinitely many pair of integers $n$ and $m$ with $n>0$
such that\begin{equation}
n\theta=m+\lambda,\quad m\in\mathbb{N}.\label{eq:51}\end{equation}

\item Let $\theta$ be an irrational number, then for any fixed real number
$\beta$ there is a positive number $\rho\ge1$ such that \begin{equation}
n\theta=m+\beta+\gamma_{n},\quad m\in\mathbb{Z},\quad|\gamma_{n}|\le\frac{1}{n^{\rho}}\label{eq:52}\end{equation}
 holds for infinitely many positive integers $n$. Furthermore, $m\in\mathbb{N}$
if $\theta>0$.
\item Given two rational numbers $\tau,\theta$ with $\tau<0$, there are
$\lambda$, $\lambda_{1}$ with $\left\{ \lambda\right\} \in\mathbb{S}(-\tau),\,\left\{ \lambda_{1}\right\} \in\mathbb{S}(\theta)$
and infinitely many integers $n$, $m$ with $n>0$ such that\begin{equation}
-n\tau=m+\lambda,\quad m\in\mathbb{N},\quad n\theta=m_{1}+\lambda_{1},\quad m_{1}\in\mathbb{Z}.\label{eq:53}\end{equation}

\item Assume $\tau<0$. If only one of $\tau$ and $\theta$ is rational,
say, $\tau$ is rational and $\theta$ is irrational, then for any
fixed real number $\beta$, there is a positive number $\rho\ge1$
and some rational number $\lambda$ with $\left\{ \lambda\right\} \in S(-\tau)$
such that \begin{equation}
n\theta=m_{1}+\beta+b_{n},\quad m_{1}\in\mathbb{Z},\quad|b_{n}|<\frac{1}{n^{\rho}}\label{eq:54}\end{equation}
 and\begin{equation}
-n\tau=m+\lambda,\quad m\in\mathbb{N}\label{eq:55}\end{equation}
 hold for infinitely many positive integers $n$.
\item If both $-\tau>0$ and $\theta$ are irrational and there exist real
numbers $\beta_{1}$, $\beta_{2}$ such that \begin{equation}
\omega(-\tau,\theta|\beta_{1},\beta_{2})>1,\label{eq:56}\end{equation}
then for any $\rho$ with\begin{equation}
0<\rho<\omega(-\tau,\theta|\beta_{1},\beta_{2})-1,\label{eq:57}\end{equation}
 there exist infinitely many integers $n$ and $m$ with $n>0$ such
that\begin{equation}
n\theta=m_{1}+\beta_{2}+b_{n},\quad m_{1}\in\mathbb{Z},\quad|b_{n}|<\frac{1}{n^{\rho}},\label{eq:58}\end{equation}
and\begin{equation}
-n\tau=m+\beta_{1}+a_{n},\quad m\in\mathbb{N},\quad|a_{n}|<\frac{1}{n^{\rho}}.\label{eq:59}\end{equation}

\end{enumerate}
\end{lem}
\begin{proof}
The first two cases are trivial, they are direct consequences of definition
of $\omega(\theta|\beta)$ and the Chebyshev's theorem. 

To prove the third assertion, we only need to consider the case that
$\lambda,\lambda_{1}\in(0,1)$. Let $\tau,\theta\in\mathbb{Q}$, if
$\theta$ or $\tau$ is an integer, then \eqref{eq:53} reduces to
\eqref{eq:51}. Assume that $\tau<0$ and $\theta$ are not integers,
let \begin{equation}
-\tau=\frac{q_{1}}{p_{1}},\quad p_{1}\ge1,\quad(p_{1},q_{1})=1,\label{eq:60}\end{equation}
 and\begin{equation}
\theta=\frac{q_{2}}{p_{2}},\quad p_{2}\ge1,\quad(p_{2},q_{2})=1.\label{eq:61}\end{equation}
 If \begin{equation}
d=(p_{1},p_{2})=1,\label{eq:62}\end{equation}
 and for any \begin{equation}
\lambda=\frac{v_{1}}{p_{1}}\in\mathbb{S}(-\tau),\quad\lambda_{1}=\frac{v_{2}}{p_{2}}\in\mathbb{S}(\theta),\label{eq:63}\end{equation}
 let $u_{1}$, $u_{2}$ be any two integers satisfying\begin{equation}
v_{1}\equiv q_{1}u_{1}(\mod p_{1}),\quad v_{2}\equiv q_{2}u_{2}(\mod p_{2}),\label{eq:64}\end{equation}
 then, by Chinese remainder theorem, there exist infinitely many positive
integers $n,n_{1},n_{2}$ such that\begin{equation}
n=p_{1}n_{1}+u_{1},\quad n=p_{2}n_{2}+u_{2}.\label{eq:65}\end{equation}
 Hence, \begin{equation}
-n\tau=n_{1}q_{1}+\frac{q_{1}u_{1}}{p_{1}}=m+\lambda\label{eq:66}\end{equation}
and\begin{equation}
n\theta=n_{2}q_{2}+\frac{q_{2}u_{2}}{p_{2}}=m_{1}+\lambda_{1}.\label{eq:67}\end{equation}
 In the case that\begin{equation}
d=(p_{1},p_{2})>1,\label{eq:68}\end{equation}
 we let\begin{equation}
p'_{1}=\frac{p_{1}}{d},\quad p'_{2}=\frac{p_{2}}{d}.\label{eq:69}\end{equation}
 From the previous discussion we see that for \begin{equation}
\lambda\in\left\{ 0,\frac{1}{p'_{1}},\dots,\frac{p'_{1}-1}{p'_{1}}\right\} \subset\mathbb{S}(-\tau)\label{eq:70}\end{equation}
 and \begin{equation}
\lambda_{1}\in\left\{ 0,\frac{1}{p'_{2}},\dots,\frac{p'_{2}-1}{p'_{2}}\right\} \subset\mathbb{S}(\theta),\label{eq:71}\end{equation}
 there exists infinitely many positive integers $n'$, $m'$ such
that\begin{equation}
-n'\tau=m+\lambda,\quad n'\theta=m{}_{1}+\lambda_{1},\quad m{}_{1},m\in\mathbb{Z},\label{eq:72}\end{equation}
where the new $n'=dn$ for each $n$ in \eqref{eq:64}. Given any
two real numbers $\tau<0$ and $\theta$, if just one of them is irrational,
say, $\tau$ is rational and $\theta$ is irrational, then by \eqref{eq:51},
there are infinitely many positive integers $n$ satisfying equations
of \eqref{eq:53}. Since $\mathbb{S}(-\tau)$ is finite set, then
there exist some $\lambda\in\mathbb{S}(-\tau)$ and infinitely many
positive integers $n$ satisfying all the equations of \eqref{eq:53}
and \eqref{eq:54}.

In the case when both $\tau<0$ and $\theta$ are irrational, \eqref{eq:57}
and \eqref{eq:58} follow from the definition of $\omega(-\tau,\theta|\beta_{1},\beta_{2}$).
\end{proof}

\section{$q$-Exponential Function $E_{q}(z)$}

The Euler's $q$-Exponential is defined by \cite{Andrews4,Gasper,Ismail2,Koekoek}\begin{equation}
E_{q}(z):=(-z;q)_{\infty}=\sum_{k=0}^{\infty}\frac{q^{k(k-1)/2}}{(q;q)_{k}}z^{k},\quad z\in\mathbb{C}.\label{eq:73}\end{equation}
For any complex number $z$, we notice that\begin{equation}
E_{q}((1-q)z)=\sum_{k=0}^{\infty}\frac{(1-q)^{k}}{(q;q)_{k}}q^{k(k-1)/2}z^{k}.\label{eq:74}\end{equation}
 By applying \eqref{eq:13} we get \begin{equation}
\frac{(1-q)^{k}q^{k(k-1)/2}}{(q;q)_{k}}\le\frac{1}{k!}\label{eq:75}\end{equation}
for any nonnegative integer $k$. Hence by Lebesgue dominated theorem
we have\begin{equation}
\lim_{q\to1}E_{q}((1-q)z)=\sum_{k=0}^{\infty}\frac{z^{k}}{k!}=e^{z},\label{eq:76}\end{equation}
and this is the reason why $E_{q}(z)$ is called a $q$-exponential.
Indeed, it is one of several $q$-analogues of $e^{z}$ in $q$-series.
From \eqref{eq:73} and \eqref{eq:74}, it is also clear that \begin{equation}
|E_{q}((1-q)z)|\le e^{|z|},\quad z\in\mathbb{C}.\label{eq:77}\end{equation}
From \eqref{eq:72} we also have \begin{equation}
|(E_{q}(z)|\le\sum_{k=0}^{\infty}\frac{q^{k/2}}{(q;q)_{k}}\left\{ q^{k^{2}/2}\left(\frac{|z|}{q}\right)^{k}\right\} \label{eq:78}\end{equation}
 and for any nonzero complex number $z$, the terms 

\begin{equation}
q^{k^{2}/2}\left(\frac{|z|}{q}\right)^{k},\quad k=0,1,...\label{eq:79}\end{equation}
 are bounded by \begin{equation}
\exp\left\{ -\frac{\log^{2}|z/q|}{2\log q}\right\} .\label{eq:80}\end{equation}
 Then we have\begin{equation}
|(E_{q}(z)|\le\exp\left\{ -\frac{\log^{2}|z/q|}{2\log q}\right\} \sum_{k=0}^{\infty}\frac{q^{k/2}}{(q;q)_{k}}\label{eq:81}\end{equation}
 or\begin{equation}
|(E_{q}(z)|\le\frac{\exp\left\{ -\frac{\log^{2}|z/q|}{2\log q}\right\} }{(\sqrt{q};q)_{\infty}}.\label{eq:82}\end{equation}
From Jacobi triple product formulas we obtain\begin{align}
\left|\theta_{3}(z;q)\right| & =\left|(q^{2},-qz,-q/z;q^{2})_{\infty}\right|\label{eq:83}\\
 & \le\frac{(q^{2};q^{2})_{\infty}}{(q;q^{2})_{\infty}^{2}}\exp\left\{ -\frac{1}{4\log q}\left[\log^{2}(|z|/q)+\log^{2}(|z|q)\right]\right\} .\nonumber \end{align}
 Since\begin{equation}
\log^{2}(|z|/q)+\log^{2}(|z|q)=2\log^{2}|z|+2\log^{2}q,\label{eq:84}\end{equation}
 then,\begin{equation}
|\theta_{3}(z;q)|\le\frac{(q^{2};q^{2})_{\infty}}{(q;q^{2})_{\infty}^{2}\sqrt{q}}\exp\left\{ -\frac{\log^{2}|z|}{2\log q}\right\} .\label{eq:85}\end{equation}
Similarly, we have\begin{eqnarray}
\left|\theta_{4}(z;q)\right| & \le & \frac{(q^{2};q^{2})_{\infty}}{(q;q^{2})_{\infty}^{2}\sqrt{q}}\exp\left\{ -\frac{\log^{2}|z|}{2\log q}\right\} ,\label{eq:86}\\
\left|\theta_{1}(z;q)\right| & \le & \frac{2\sqrt[4]{q}(q^{2};q^{2})_{\infty}\cosh(\frac{1}{2}\log|z|)}{(q;q^{2})_{\infty}^{2}}\exp\left\{ -\frac{\log^{2}|z|}{2\log q}\right\} ,\label{eq:87}\\
\left|\theta_{2}(z;q)\right| & \le & \frac{2\sqrt[4]{q}(q^{2};q^{2})_{\infty}\cosh(\frac{1}{2}\log|z|)}{(q;q^{2})_{\infty}^{2}}\exp\left\{ -\frac{\log^{2}|z|}{2\log q}\right\} .\label{eq:88}\end{eqnarray}

\subsection{Asymptotics for $E_{q}(z)$ \label{sec:q-exponential asymptotics}}

Assume that \begin{equation}
s=\tau+\frac{2\pi i\theta}{\log q},\quad\tau,\theta\in\mathbb{R}.\label{eq:89}\end{equation}
 Then for any nonzero complex number $z$, we have \begin{eqnarray}
E_{q}(-q^{ns+1/2}z) & = & \sum_{k=0}^{\infty}\frac{q^{k^{2}/2}}{(q;q)_{k}}\left(-ze^{2\pi in\theta}q^{n\tau}\right)^{k}.\label{eq:90}\end{eqnarray}
 For the Euler $q$-Exponential $E_{q}(z)$ we have the following
results:

\begin{thm*}
\label{thm:q-exponential}For any nonzero complex number $z$, let\begin{equation}
j_{n}=\left\lfloor \frac{q^{4}\log n}{-\log q}\right\rfloor ,\quad k_{n}=\left\lfloor \frac{-n\tau}{2}\right\rfloor ,\label{eq:91}\end{equation}
we have the following results for $E_{q}(z)$:
\begin{enumerate}
\item If $\tau>0$, then we have \begin{equation}
E_{q}(-q^{ns+1/2}z)=1+r_{qe}(n|1),\label{eq:92}\end{equation}
and\begin{equation}
|r_{qe}(n|1)|\le\frac{|z|q^{n\tau+1/2}}{1-q}\exp\left\{ \frac{|z|q^{n\tau+1/2}}{1-q}\right\} .\label{eq:93}\end{equation}

\item Assume that $\tau=0$. If for some fixed real numbers $\beta$ and
$\rho\ge1$ there are infinitely many positive integers $n$ such
that \begin{equation}
n\theta=m+\beta+b_{n},\quad|b_{n}|<\frac{1}{n^{\rho}},\quad m\in\mathbb{Z},\label{eq:94}\end{equation}
 then,\begin{equation}
E_{q}(-q^{ns+1/2}z)=E_{q}(-zq^{1/2}e^{2\pi i\beta};q)_{\infty}+r_{qe}(n|3),\label{eq:95}\end{equation}
and\begin{equation}
|r_{qe}(n|3)|\le24\exp\left(\frac{|z|\sqrt{q}}{1-q}\right)\left\{ \frac{\log n}{n^{\rho}}+\frac{1}{j_{n}!}\left(\frac{|z|\sqrt{q}}{1-q}\right)^{j_{n}}\right\} \label{eq:96}\end{equation}
for $n$ sufficiently large. 
\item Assume that $\tau<0$. If for some fixed real numbers $\lambda$ and
$\lambda_{1}$ there are infinitely many positive integers $n$ such
that \begin{equation}
-n\tau=m+\lambda,\quad m\in\mathbb{N},\quad n\theta=m_{1}+\lambda_{1},\quad m_{1}\in\mathbb{Z}.\label{eq:97}\end{equation}
Then\begin{equation}
\frac{E_{q}(-q^{ns+1/2}z)(q;q)_{\infty}}{(-ze^{2\pi in\theta})^{m}q^{m(n\tau+m/2)}}=\theta_{4}\left(z^{-1}q^{\lambda}e^{-2\pi i\lambda_{1}};\sqrt{q}\right)+r_{qe}(n|4),\label{eq:98}\end{equation}
and \begin{align}
|r_{qe}(n|4)| & \le4\theta_{3}\left(|z|^{-1}q^{\lambda};\sqrt{q}\right)\left\{ \frac{q^{k_{n}}}{1-q}+\frac{q^{k_{n}^{2}/2+\lambda k_{n}}}{\left|z\right|^{k_{n}}}\right\} \label{eq:99}\end{align}
 for $n$ sufficiently large.
\item Assume that $\tau<0$. If for some fixed real numbers $\beta$, $\lambda$,
$\rho\ge1$ there exist infinitely many positive integers $n$ such
that\begin{equation}
-n\tau=m+\lambda,\quad m\in\mathbb{N},\label{eq:100}\end{equation}
 and\begin{equation}
n\theta=m_{1}+\beta+b_{n},\quad|b_{n}|<\frac{1}{n^{\rho}},\quad m_{1}\in\mathbb{Z},\label{eq:101}\end{equation}
 then,\begin{equation}
\frac{E_{q}(-q^{ns+1/2}z)(q;q)_{\infty}}{(-ze^{2\pi in\theta})^{m}q^{m(n\tau+m/2)}}=\theta_{4}\left(z^{-1}q^{\lambda}e^{-2\pi i\beta};\sqrt{q}\right)+r_{qe}(n|5),\label{eq:102}\end{equation}
 and\begin{align}
|r_{qe}(n|5)| & \le48\theta_{3}\left(|z|^{-1}q^{\lambda};\sqrt{q}\right)\left\{ \frac{q^{j_{n}^{2}/2+\lambda j_{n}}}{\left|z\right|^{j_{n}}}+\left|z\right|^{j_{n}}q^{j_{n}^{2}/2-\lambda j_{n}}+\frac{q^{k_{n}}}{1-q}+\frac{\log n}{n^{\rho}}\right\} \label{eq:103}\end{align}
 for $n$ sufficiently large.
\item Assume that $\tau<0$. If for some fixed real numbers $\beta$, $\lambda$,
$\rho\ge1$ there are infinitely many positive integers $n$ such
that\begin{equation}
-n\tau=m+\beta+a_{n},\quad|a_{n}|<\frac{1}{n^{\rho}},\quad m\in\mathbb{N},\quad n\theta=m_{1}+\lambda,\quad m_{1}\in\mathbb{Z},\label{eq:104}\end{equation}
then,\begin{equation}
\frac{E_{q}(-q^{ns+1/2}z)(q;q)_{\infty}}{(-ze^{2\pi in\theta})^{m}q^{m(n\tau+m/2)}}=\theta_{4}\left(z^{-1}q^{\beta}e^{-2\pi i\lambda};\sqrt{q}\right)+r_{qe}(n|6),\label{eq:105}\end{equation}
 and\begin{align}
|r_{qe}(n|6)| & \le6\theta_{3}(|z|^{-1}q^{\beta};\sqrt{q})\left\{ \frac{\log n}{n^{\rho}}+\frac{q^{k_{n}}}{1-q}+\left|z\right|^{j_{n}}q^{j_{n}^{2}/2-(\beta+1)j_{n}}+\frac{q^{j_{n}^{2}/2+(\beta-1)j_{n}}}{\left|z\right|^{j_{n}}}\right\} \label{eq:106}\end{align}
for $n$ sufficiently large.
\item Assume that $\tau<0$. If for some fixed real numbers $\beta_{1}$,
$\beta_{2}$ and $\rho>0$ there exist infinitely many positive integers
$n$ such that\begin{equation}
-n\tau=m+\beta_{1}+a_{n},\quad|a_{n}|<\frac{1}{n^{\rho}},\quad m\in\mathbb{N},\label{eq:107}\end{equation}
 and\begin{equation}
n\theta=m_{1}+\beta_{2}+b_{n},\quad|b_{n}|<\frac{1}{n^{\rho}},\quad m_{1}\in\mathbb{Z}.\label{eq:108}\end{equation}
Then,\begin{equation}
\frac{E_{q}(-q^{ns+1/2}z)(q;q)_{\infty}}{(-ze^{2\pi in\theta})^{m}q^{m(n\tau+m/2)}}=\theta_{4}\left(z^{-1}q^{\beta_{1}}e^{-2\pi i\beta_{2}};\sqrt{q}\right)+r_{qe}(n|7),\label{eq:109}\end{equation}
 and\begin{align}
|r_{qe}(n|7)| & \le54\theta_{3}(|z|^{-1}q^{\beta_{1}};\sqrt{q})\left\{ \frac{\log n}{n^{\rho}}+\frac{q^{k_{n}}}{1-q}+q^{j_{n}^{2}/2-(\beta_{1}+1)j_{n}}|z|^{j_{n}}+\frac{q^{j_{n}^{2}/2+(\beta_{1}-1)j_{n}}}{|z|^{j_{n}}}\right\} \label{eq:110}\end{align}
 for $n$ sufficiently large.
\end{enumerate}
\end{thm*}
\begin{rem}
If only one of $\tau$ and $\theta$ is an irrational number, say
$\tau$, then the real number $\beta$ could be fixed as $0$. If
$\tau$ is an algebraic of degree $\ell$, then the order of error
term can be no better than $\mathcal{O}(n^{1-\ell})$, while if the
number is a Liouville number, then, the order is better than $\mathcal{O}(n^{-r})$
for any real number $r$. In case of an arbitrary irrational number
$\tau$ and an arbitrary real number $\beta$, the error term is $\mathcal{O}(n^{-1+\epsilon})$
for any small $\epsilon>0$ by the Chebyshev's theorem. This remark
applies to all the asymptotics in this work.
\end{rem}
In the following corollaries we assume that \begin{equation}
z:=e^{2\pi u},\quad u\in\mathbb{R}.\label{eq:111}\end{equation}

\begin{cor}
\label{cor:q-exponential-power}Assume that \begin{equation}
q=\exp(-2n^{-a}\pi),\quad0<a<\frac{1}{2},\quad n\in\mathbb{N},\label{eq:112}\end{equation}
we have the following results for $E_{q}(z)$:
\begin{enumerate}
\item Assume that $\tau>0$, we have\begin{equation}
E_{q}(-\exp2\pi(u+n^{1-a}\tau+n^{-a}/2))=1+\mathcal{O}\left(n^{a}e^{-2\pi\tau n^{1-a}}\right)\label{eq:113}\end{equation}
for $n$ sufficiently large.
\item Assume that $\tau<0$, if for some fixed real number $\lambda$ there
are infinitely many positive integer $n$ such that \begin{equation}
-n\tau=m+\lambda,\quad m\in\mathbb{N},\label{eq:114}\end{equation}
then,\begin{align}
 & E_{q}(-\exp2\pi(u+\tau n^{1-a}+n^{-a}/2))\label{eq:115}\\
 & =\frac{2\exp\left(\pi n^{-a}(n^{a}u-n\tau)^{2}\right)\left\{ \cos\pi(n^{a}u+\lambda)+\mathcal{O}\left(e^{-2\pi n^{a}}\right)\right\} }{(-1)^{\lambda+n\tau}\exp\left(\pi n^{a}/6+\pi n^{-a}/12\right)}\nonumber \end{align}
 for $n$ sufficiently large.
\end{enumerate}
\end{cor}
Let \begin{equation}
q=\exp(-\frac{2\pi}{\gamma\log n}),\quad n\ge2,\quad\gamma>0,\label{eq:116}\end{equation}
 then we have the following:

\begin{cor}
\label{cor:q-exponential-log}we have the following results for $E_{q}(z)$:
\begin{enumerate}
\item Assume that $\tau>0$, we have\begin{equation}
E_{q}\left(-\exp2\pi\left(u-\frac{n\tau+1/2}{\gamma\log n}\right)\right)=1+\mathcal{O}\left(\log n\exp\left(-\frac{2n\pi\tau}{\gamma\log n}\right)\right)\label{eq:117}\end{equation}
for $n$ sufficiently large.
\item Assume that $\tau<0$ and for some fixed real number $\lambda$ there
are infinitely many positive integer $n$ such that \begin{equation}
-n\tau=m+\lambda,\quad m\in\mathbb{N},\label{eq:118}\end{equation}
then,\begin{align}
 & E_{q}\left(-\exp2\pi\left(u-\frac{n\tau+1/2}{\gamma\log n}\right)\right)\label{eq:119}\\
 & =(-1)^{\tau n+\lambda}\exp\left(\frac{\pi(\gamma u\log n-n\tau)^{2}}{\gamma\log n}-\frac{\pi}{12\gamma\log n}\right)\nonumber \\
 & \times2n^{-\gamma\pi/6}\left\{ \cos\pi(\gamma u\log n+\lambda)+\mathcal{O}(n^{-2\pi\gamma})\right\} \nonumber \end{align}
for $n$ sufficiently large.
\item Assume that $\tau<0$, if for some fixed real numbers $\beta$, $\rho\ge1$
and $\lambda$ there exist infinitely many positive integers $n$
such that\begin{equation}
-n\tau=m+\beta+a_{n},\quad|a_{n}|<\frac{1}{n^{\rho}},\quad m\in\mathbb{N},\label{eq:120}\end{equation}
then,\begin{align}
 & E_{q}\left(-\exp2\pi\left(u-\frac{n\tau+1/2}{\gamma\log n}\right)\right)\label{eq:121}\\
 & =\exp\left(\frac{\pi(\gamma u\log n-n\tau)^{2}}{\gamma\log n}-\frac{\pi}{12\gamma\log n}\right)\nonumber \\
 & \times\frac{2\left\{ \cos\pi(\gamma u\log n+\beta)+\mathcal{O}(n^{-8\rho/9}\log n)\right\} }{(-1)^{m}n^{2\rho/27}}\nonumber \end{align}
for $n$ sufficiently large, where \begin{equation}
\gamma=\frac{4\rho}{9\pi}.\label{eq:122}\end{equation}

\end{enumerate}
\end{cor}

\subsection{Proofs for Theorem \ref{thm:q-exponential} }

Assume that $\tau<0$ and \begin{equation}
-n\tau=m+c_{n},\quad m\in\mathbb{N},\quad n\theta=m_{1}+d_{n},\quad m_{1}\in\mathbb{Z},\label{123}\end{equation}
 then,\begin{align}
E(-q^{ns+1/2}z) & =\sum_{k=0}^{\infty}\frac{q^{k^{2}/2}}{(q;q)_{k}}\left(-ze^{2\pi in\theta}q^{n\tau}\right)^{k}\label{eq:124}\\
 & =\sum_{k=0}^{m}\frac{q^{k^{2}/2}}{(q;q)_{k}}\left(-ze^{2\pi in\theta}q^{n\tau}\right)^{k}+\sum_{k=m+1}^{\infty}\frac{q^{k^{2}/2}}{(q;q)_{k}}\left(-ze^{2\pi in\theta}q^{n\tau}\right)^{k}\nonumber \\
 & =s_{1}+s_{2}.\nonumber \end{align}
 Reverse summation order in $s_{1}$,\begin{equation}
\frac{(q;q)_{\infty}s_{1}}{(-ze^{2\pi in\theta})^{m}q^{m(n\tau+m/2)}}=\sum_{k=0}^{m}q^{k^{2}/2}(-q^{c_{n}}z^{-1}e^{-2\pi id_{n}})^{k}e(k,n)\label{eq:125}\end{equation}
 with\begin{equation}
e(k,n)=\frac{(q;q)_{\infty}}{(q;q)_{m-k}},\label{eq:126}\end{equation}
 then,\begin{equation}
|e(k,n)|\le1\label{eq:127}\end{equation}
 for $0\le k\le m$. An application of Lemma \ref{lem:1} gives \begin{equation}
|e(k,n)-1|=\left|r_{1}(q;m-k)\right|\le\frac{2q^{2+k_{n}}}{1-q}.\label{eq:128}\end{equation}
for $0\le k\le k_{n}-1$ and $n$ sufficiently large.

We shift the summation index from $k$ to $k+m$ in $s_{2}$,\begin{equation}
\frac{(q;q)_{\infty}s_{2}}{(-ze^{2\pi in\theta})^{m}q^{m(n\tau+m/2)}}=\sum_{k=1}^{\infty}q^{k^{2}/2}(-zq^{-c_{n}}e^{2\pi id_{n}})^{k}f(k,n),\label{eq:129}\end{equation}
with\begin{equation}
f(k,n)=\frac{(q;q)_{\infty}}{(q;q)_{m+k}},\label{eq:130}\end{equation}
 then,\begin{equation}
|f(k,n)|\le1,\label{eq:131}\end{equation}
 and\begin{equation}
|f(k,n)-1|=\left|r_{1}(q;m+k)\right|\le\frac{2q^{k_{n}+2}}{1-q}\label{eq:132}\end{equation}
 for $k\in\mathbb{N}$ and $n$ sufficiently large.

\subsubsection{Proof for case 1}

If we write\begin{equation}
E_{q}(-q^{ns+1/2}z)=1+r_{qe}(n|1),\label{eq:133}\end{equation}
then,\begin{equation}
r_{qe}(n|1)=\sum_{k=1}^{\infty}\frac{q^{k^{2}/2}}{(q;q)_{k}}\left(-ze^{2\pi in\theta}q^{n\tau}\right)^{k}.\label{eq:134}\end{equation}
 Hence,\begin{align}
|r_{qe}(n|1)| & \le\sum_{k=1}^{\infty}\frac{q^{k^{2}/2}}{(q;q)_{k}}(|z|q^{n\tau})^{k}\le\sum_{k=1}^{\infty}\frac{(|z|q^{n\tau+1/2})^{k}}{k!(1-q)^{k}}\frac{(1-q)^{k}k!q^{k(k-1)/2}}{(q;q)_{k}}\label{eq:135}\\
 & \le\frac{|z|q^{n\tau+1/2}}{1-q}\sum_{k=0}^{\infty}\frac{(|z|q^{n\tau+1/2})^{k}}{k!(1-q)^{k}}\le\frac{|z|q^{n\tau+1/2}}{1-q}\exp\left\{ \frac{|z|q^{n\tau+1/2}}{1-q}\right\} \nonumber \end{align}
 by an application of \eqref{eq:13}.

\subsubsection{Proof for case 2}

It is clear that \begin{equation}
1<j_{n}<\frac{n^{\rho}}{2\pi}\label{eq:136}\end{equation}
 for $n$ sufficiently large, then,

\begin{align}
 & E(-q^{ns+1/2}z)=\sum_{k=0}^{\infty}\frac{q^{k^{2}/2}}{(q;q)_{k}}\left(-ze^{2\pi i\beta}\right)^{k}e^{2\pi ikb_{n}}\label{eq:137}\\
 & =\sum_{k=0}^{\infty}\frac{q^{k^{2}/2}}{(q;q)_{k}}\left(-ze^{2\pi i\beta}\right)^{k}-\sum_{k=j_{n}}^{\infty}\frac{q^{k^{2}/2}}{(q;q)_{k}}\left(-ze^{2\pi i\beta}\right)^{k}\nonumber \\
 & +\sum_{k=0}^{j_{n}-1}\frac{q^{k^{2}/2}}{(q;q)_{k}}\left(-ze^{2\pi i\beta}\right)^{k}\left\{ e^{2\pi ikb_{n}}-1\right\} +\sum_{k=j_{n}}^{\infty}\frac{q^{k^{2}/2}}{(q;q)_{k}}\left(-ze^{2\pi i\beta}\right)^{k}e^{2\pi ikb_{n}}\nonumber \\
 & =E_{q}(-zq^{1/2}e^{2\pi i\beta})+s_{1}+s_{2}+s_{3}.\nonumber \end{align}
Then, \begin{align}
|s_{1}+s_{3}| & \le2\sum_{k=j_{n}}^{\infty}\frac{q^{k^{2}/2}}{(q;q)_{k}}|z|^{k}\le2\sum_{k=j_{n}}^{\infty}\frac{1}{k!}\left(\frac{|z|\sqrt{q}}{1-q}\right)^{k}\label{eq:138}\\
 & \le\frac{\left(\frac{|z|\sqrt{q}}{1-q}\right)^{j_{n}}}{j_{n}!}\sum_{k=0}^{\infty}\frac{\left(\frac{|z|\sqrt{q}}{1-q}\right)^{k}}{k!}\le\frac{2\exp\left(\frac{|z|\sqrt{q}}{1-q}\right)}{j_{n}!}\left(\frac{|z|\sqrt{q}}{1-q}\right)^{j_{n}}.\nonumber \end{align}
 For $n$ sufficiently large, \begin{align}
|s_{2}| & \le2\pi j_{n}|b_{n}|e^{2\pi j_{n}|b_{n}|}\sum_{k=0}^{\infty}\frac{q^{k^{2}/2}|z|^{k}}{(q;q)_{k}}\le24\exp\left(\frac{|z|\sqrt{q}}{1-q}\right)\frac{\log n}{n^{\rho}}.\label{eq:139}\end{align}
 Hence we have proved\begin{equation}
E_{q}(-q^{ns+1/2}z)=E_{q}(-zq^{1/2}e^{2\pi i\beta})+r_{qe}(n|3),\label{eq:140}\end{equation}
 where\begin{equation}
r_{qe}(n|3)=s_{1}+s_{2}+s_{3},\label{eq:141}\end{equation}
 and \begin{equation}
|r_{qe}(n|3)|\le24\exp\left(\frac{|z|\sqrt{q}}{1-q}\right)\left\{ \frac{\log n}{n^{\rho}}+\frac{1}{j_{n}!}\left(\frac{|z|\sqrt{q}}{1-q}\right)^{j_{n}}\right\} \label{eq:142}\end{equation}
 for $n$ sufficiently large.

\subsubsection{Proof for case 4}

Notice that

\begin{align}
 & \frac{(q;q)_{\infty}s_{1}}{(-ze^{2\pi in\theta})^{m}q^{m(n\tau+m/2)}}=\sum_{k=0}^{m}q^{k^{2}/2}(-q^{\lambda}z^{-1}e^{-2\pi i\lambda_{1}})^{k}e(k,n)\label{eq:143}\\
 & =\sum_{k=0}^{\infty}q^{k^{2}/2}(-q^{\lambda}z^{-1}e^{-2\pi i\lambda_{1}})^{k}-\sum_{k=k_{n}}^{\infty}q^{k^{2}/2}(-q^{\lambda}z^{-1}e^{-2\pi i\lambda_{1}})^{k}\nonumber \\
 & +\sum_{k=0}^{k_{n}-1}q^{k^{2}/2}(-q^{\lambda}z^{-1}e^{-2\pi i\lambda_{1}})^{k}\left\{ e(k,n)-1\right\} +\sum_{k=k_{n}}^{m}q^{k^{2}/2}(-q^{\lambda}z^{-1}e^{-2\pi i\lambda_{1}})^{k}e(k,n)\nonumber \\
 & =\sum_{k=0}^{\infty}q^{k^{2}/2}(-q^{\lambda}z^{-1}e^{-2\pi i\lambda_{1}})^{k}+s_{11}+s_{12}+s_{13}.\nonumber \end{align}
 Thus,

\begin{align}
|s_{11}+s_{13}| & \le2\sum_{k=k_{n}}^{\infty}q^{k^{2}/2}\left(\frac{q^{\lambda}}{|z|}\right)^{k}\le2\frac{q^{k_{n}^{2}/2+\lambda k_{n}}}{\left|z\right|^{k_{n}}}\sum_{k=0}^{\infty}q^{k^{2}/2}\left(\frac{q^{\lambda+k_{n}}}{|z|}\right)^{k}\label{eq:144}\\
 & \le2\theta_{3}\left(|z|^{-1}q^{\lambda};\sqrt{q}\right)\frac{q^{k_{n}^{2}/2+\lambda k_{n}}}{\left|z\right|^{k_{n}}},\nonumber \end{align}
 and for $n$ sufficiently large,\begin{align}
|s_{12}| & \le\frac{2q^{k_{n}}}{1-q}\sum_{k=0}^{\infty}q^{k^{2}/2}\left(\frac{q^{\lambda}}{|z|}\right)^{k}\le\frac{2\theta_{3}\left(|z|^{-1}q^{\lambda};\sqrt{q}\right)}{1-q}q^{k_{n}}.\label{eq:145}\end{align}
 Let \begin{equation}
r_{1}(n)=s_{11}+s_{12}+s_{13},\label{eq:146}\end{equation}
then,\begin{equation}
\frac{(q;q)_{\infty}s_{1}}{(-ze^{2\pi in\theta})^{m}q^{m(n\tau+m/2)}}=\sum_{k=0}^{\infty}q^{k^{2}/2}(-q^{\lambda}z^{-1}e^{-2\pi i\lambda_{1}})^{k}+r_{1}(n),\label{eq:147}\end{equation}
and\begin{align}
|r_{1}(n)| & \le2\theta_{3}\left(|z|^{-1}q^{\lambda};\sqrt{q}\right)\left\{ \frac{q^{k_{n}}}{1-q}+\frac{q^{k_{n}^{2}/2+\lambda k_{n}}}{\left|z\right|^{k_{n}}}\right\} \label{eq:148}\end{align}
for $n$ sufficiently large.

In the second sum,\begin{align}
 & \frac{(q;q)_{\infty}s_{2}}{(-ze^{2\pi in\theta})^{m}q^{m(n\tau+m/2)}}=\sum_{k=1}^{\infty}q^{k^{2}/2}(-zq^{-\lambda}e^{2\pi i\lambda_{1}})^{k}f(k,n)\label{eq:149}\\
 & =\sum_{k=1}^{\infty}q^{k^{2}/2}(-zq^{-\lambda}e^{2\pi i\lambda_{1}})^{k}+\sum_{k=1}^{\infty}q^{k^{2}/2}(-zq^{-\lambda}e^{2\pi i\lambda_{1}})^{k}\left\{ f(k,n)-1\right\} \nonumber \\
 & =\sum_{k=-\infty}^{-1}q^{k^{2}/2}(-q^{\lambda}z^{-1}e^{-2\pi i\lambda_{1}})^{k}+r_{2}(n),\nonumber \end{align}
 and\begin{align}
|r_{2}(n)| & \le\frac{2q^{k_{n}}}{1-q}\sum_{k=1}^{\infty}q^{k^{2}/2}(q^{-\lambda}|z|)^{k}\le\frac{2\theta_{3}\left(|z|^{-1}q^{\lambda};\sqrt{q}\right)}{1-q}q^{k_{n}}\label{eq:150}\end{align}
 for $n$ sufficiently large. Hence,\begin{equation}
\frac{E_{q}(-q^{ns+1/2}z)(q;q)_{\infty}}{(-ze^{2\pi in\theta})^{m}q^{m(n\tau+m/2)}}=\theta_{4}\left(z^{-1}q^{\lambda}e^{-2\pi i\lambda_{1}};\sqrt{q}\right)+r_{qe}(n|4),\label{eq:151}\end{equation}
 with\begin{align}
|r_{qe}(n|4)| & \le4\theta_{3}\left(|z|^{-1}q^{\lambda};\sqrt{q}\right)\left\{ \frac{q^{k_{n}}}{1-q}+\frac{q^{k_{n}^{2}/2+\lambda k_{n}}}{\left|z\right|^{k_{n}}}\right\} \label{eq:152}\end{align}
 for $n$ sufficiently large.

\subsubsection{Proof for case 5}

In this case we have\begin{align}
 & \frac{(q;q)_{\infty}s_{1}}{(-ze^{2\pi in\theta})^{m}q^{m(n\tau+m/2)}}=\sum_{k=0}^{m}q^{k^{2}/2}(-q^{\lambda}z^{-1}e^{-2\pi i\beta})^{k}e^{-2\pi ikb_{n}}e(k,n)\label{eq:153}\\
 & =\sum_{k=0}^{\infty}q^{k^{2}/2}(-q^{\lambda}z^{-1}e^{-2\pi i\beta})^{k}-\sum_{k=j_{n}}^{\infty}q^{k^{2}/2}(-q^{\lambda}z^{-1}e^{-2\pi i\beta})^{k}\nonumber \\
 & +\sum_{k=0}^{j_{n}-1}q^{k^{2}/2}(-q^{\lambda}z^{-1}e^{-2\pi i\beta})^{k}\left\{ e(k,n)-1\right\} +\sum_{k=0}^{j_{n}-1}q^{k^{2}/2}(-q^{\lambda}z^{-1}e^{-2\pi i\beta})^{k}e(k,n)\left\{ e^{-2\pi ikb_{n}}-1\right\} \nonumber \\
 & +\sum_{k=j_{n}}^{m}q^{k^{2}/2}(-q^{\lambda}z^{-1}e^{-2\pi i\beta})^{k}e^{-2\pi ikb_{n}}e(k,n)\nonumber \\
 & =\sum_{k=0}^{\infty}q^{k^{2}/2}(-q^{\lambda}z^{-1}e^{-2\pi i\beta})^{k}+s_{11}+s_{12}+s_{13}+s_{14}.\nonumber \end{align}
Then,\begin{align}
|s_{11}+s_{14}| & \le2\sum_{k=j_{n}}^{\infty}q^{k^{2}/2}\left(\frac{q^{\lambda}}{|z|}\right)^{k}\le2\frac{q^{j_{n}^{2}/2+\lambda j_{n}}}{\left|z\right|^{j_{n}}}\sum_{k=0}^{\infty}q^{k^{2}/2}\left(\frac{q^{\lambda+j_{n}}}{|z|}\right)^{k}\label{eq:154}\\
 & \le2\theta_{3}\left(|z|^{-1}q^{\lambda};\sqrt{q}\right)\frac{q^{j_{n}^{2}/2+\lambda j_{n}}}{\left|z\right|^{j_{n}}}.\nonumber \end{align}
We also have\begin{align}
|s_{13}| & \le2\pi j_{n}|b_{n}|e^{2\pi j_{n}|b_{n}|}\sum_{k=0}^{\infty}q^{k^{2}/2}\left(\frac{q^{\lambda}}{|z|}\right)^{k}\le24\theta_{3}\left(|z|^{-1}q^{\lambda};\sqrt{q}\right)\frac{\log n}{n^{\rho}},\label{eq:155}\end{align}
and\begin{equation}
|s_{12}|\le\frac{2q^{k_{n}}}{1-q}\theta_{3}\left(|z|^{-1}q^{\lambda};\sqrt{q}\right)\label{eq:156}\end{equation}
 for $n$ sufficiently large. 

Let \begin{equation}
r_{1}(n)=s_{11}+s_{12}+s_{13}+s_{14},\label{eq:157}\end{equation}
then\begin{equation}
\frac{(q;q)_{\infty}s_{1}}{(-ze^{2\pi in\theta})^{m}q^{m(n\tau+m/2)}}=\sum_{k=0}^{\infty}q^{k^{2}/2}(-z^{-1}q^{\lambda}e^{-2\pi i\beta})^{k}+r_{1}(n),\label{eq:158}\end{equation}
and\begin{equation}
|r_{1}(n)|\le24\theta_{3}\left(|z|^{-1}q^{\lambda};\sqrt{q}\right)\left\{ \frac{q^{j_{n}^{2}/2+\lambda j_{n}}}{\left|z\right|^{j_{n}}}+\frac{q^{k_{n}}}{1-q}+\frac{\log n}{n^{\rho}}\right\} \label{eq:159}\end{equation}
for $n$ sufficiently large.

Similarly,\begin{align}
 & \frac{(q;q)_{\infty}s_{2}}{(-ze^{2\pi in\theta})^{m}q^{m(n\tau+m/2)}}=\sum_{k=1}^{\infty}q^{k^{2}/2}(-zq^{-\lambda}e^{2\pi i\beta})^{k}f(k,n)e^{2\pi ikb_{n}}\label{eq:160}\\
 & =\sum_{k=1}^{\infty}q^{k^{2}/2}(-zq^{-\lambda}e^{2\pi i\beta})^{k}-\sum_{k=j_{n}}^{\infty}q^{k^{2}/2}(-zq^{-\lambda}e^{2\pi i\beta})^{k}\nonumber \\
 & +\sum_{k=1}^{j_{n}-1}q^{k^{2}/2}(-zq^{-\lambda}e^{2\pi i\beta})^{k}\left\{ f(k,n)-1\right\} +\sum_{k=1}^{j_{n}-1}q^{k^{2}/2}(-zq^{-\lambda}e^{2\pi i\beta})^{k}f(k,n)\left\{ e^{2\pi ikb_{n}}-1\right\} \nonumber \\
 & +\sum_{k=j_{n}}^{\infty}q^{k^{2}/2}(-zq^{-\lambda}e^{2\pi i\beta})^{k}f(k,n)e^{2\pi ikb_{n}}\nonumber \\
 & =\sum_{k=-\infty}^{-1}q^{k^{2}/2}(-z^{-1}q^{\lambda}e^{-2\pi i\beta})^{k}+s_{21}+s_{22}+s_{23}+s_{24}.\nonumber \end{align}
 Then,\begin{align}
|s_{21}+s_{24}| & \le2\sum_{k=j_{n}}^{\infty}q^{k^{2}/2}(|z|q^{-\lambda})^{k}\le2\left|z\right|^{j_{n}}q^{j_{n}^{2}/2-\lambda j_{n}}\sum_{k=0}^{\infty}q^{k^{2}/2}(|z|q^{-\lambda+j_{n}})^{k}\label{eq:161}\\
 & \le2\theta_{3}\left(|z|^{-1}q^{\lambda};\sqrt{q}\right)\left|z\right|^{j_{n}}q^{j_{n}^{2}/2-\lambda j_{n}}.\nonumber \end{align}
 For sufficiently large $n$ we have\begin{align}
|s_{22}| & \le\frac{2q^{k_{n}}}{1-q}\sum_{k=1}^{\infty}q^{k^{2}/2}(|z|q^{-\lambda})^{k}\le\frac{2q^{k_{n}}}{1-q}\theta_{3}\left(|z|^{-1}q^{\lambda};\sqrt{q}\right),\label{eq:162}\end{align}
 and\begin{align}
|s_{23}| & \le2\pi j_{n}|b_{n}|e^{2\pi j_{n}|b_{n}|}\sum_{k=1}^{\infty}q^{k^{2}/2}(|z|q^{-\lambda})^{k}\le24\theta_{3}\left(|z|^{-1}q^{\lambda};\sqrt{q}\right)\frac{\log n}{n^{\rho}}\label{eq:163}\end{align}
for $n$ sufficiently large. 

Let \begin{equation}
r_{2}(n)=s_{21}+s_{22}+s_{23}+s_{24},\label{eq:164}\end{equation}
then,\begin{equation}
\frac{(q;q)_{\infty}s_{2}}{(-ze^{2\pi in\theta})^{m}q^{m(n\tau+m/2)}}=\sum_{k=-\infty}^{-1}q^{k^{2}/2}(-z^{-1}q^{\lambda}e^{-2\pi i\beta})^{k}+r_{2}(n),\label{eq:165}\end{equation}
with\begin{equation}
|r_{2}(n)|\le24\theta_{3}\left(|z|^{-1}q^{\lambda};\sqrt{q}\right)\left\{ \frac{\log n}{n^{\rho}}+\frac{q^{k_{n}}}{1-q}+\left|z\right|^{j_{n}}q^{j_{n}^{2}/2-\lambda j_{n}}\right\} \label{eq:166}\end{equation}
 for sufficiently large $n$.

Therefore,\begin{equation}
\frac{E_{q}(-q^{ns+1/2}z)(q;q)_{\infty}}{(-ze^{2\pi in\theta})^{m}q^{m(n\tau+m/2)}}=\theta_{4}\left(z^{-1}q^{\lambda}e^{-2\pi i\beta};\sqrt{q}\right)+r_{qe}(n|5)\label{eq:167}\end{equation}
 with\begin{align}
|r_{qe}(n|5)| & \le48\theta_{3}\left(|z|^{-1}q^{\lambda};\sqrt{q}\right)\left\{ \frac{q^{j_{n}^{2}/2+\lambda j_{n}}}{\left|z\right|^{j_{n}}}+\left|z\right|^{j_{n}}q^{j_{n}^{2}/2-\lambda j_{n}}+\frac{q^{k_{n}}}{1-q}+\frac{\log n}{n^{\rho}}\right\} \label{eq:168}\end{align}
 for $n$ sufficiently large.

\subsubsection{Proof for case 6}

Notice that\begin{align}
 & \frac{(q;q)_{\infty}s_{1}}{(-ze^{2\pi in\theta})^{m}q^{m(n\tau+m/2)}}=\sum_{k=0}^{m}q^{k^{2}/2}(-q^{\beta}z^{-1}e^{-2\pi i\lambda})^{k}q^{ka_{n}}e(k,n)\label{eq:169}\\
 & =\sum_{k=0}^{\infty}q^{k^{2}/2}(-q^{\beta}z^{-1}e^{-2\pi i\lambda})^{k}-\sum_{k=j_{n}}^{\infty}q^{k^{2}/2}(-q^{\beta}z^{-1}e^{-2\pi i\lambda})^{k}\nonumber \\
 & +\sum_{k=0}^{j_{n}-1}q^{k^{2}/2}(-q^{\beta}z^{-1}e^{-2\pi i\lambda})^{k}\left\{ e(k,n)-1\right\} +\sum_{k=0}^{j_{n}-1}q^{k^{2}/2}(-q^{\beta}z^{-1}e^{-2\pi i\lambda})^{k}e(k,n)\left\{ q^{ka_{n}}-1\right\} \nonumber \\
 & +\sum_{k=j_{n}}^{m}q^{k^{2}/2}(-q^{\beta}z^{-1}e^{-2\pi i\lambda})^{k}q^{ka_{n}}e(k,n)=\sum_{k=0}^{\infty}q^{k^{2}/2}(-q^{\beta}z^{-1}e^{-2\pi i\lambda})^{k}+s_{11}+s_{12}+s_{13}+s_{14}.\nonumber \end{align}
Since\begin{equation}
|q^{ka_{n}}-1|\le j_{n}|a_{n}|\log q^{-1}e^{j_{n}|a_{n}|\log q^{-1}}\le\frac{3\log n}{n^{\rho}},\label{eq:170}\end{equation}
and\begin{equation}
|q^{ka_{n}}|\le e^{j_{n}|a_{n}|\log q^{-1}}<3\label{eq:171}\end{equation}
 for $0\le k\le j_{n}-1$ and\begin{equation}
|q^{ka_{n}}|\le q^{-2k}\label{eq:172}\end{equation}
 for any $k\ge0$, then,

\begin{align}
|s_{11}+s_{14}| & \le2\sum_{k=j_{n}}^{\infty}q^{k^{2}/2}(q^{\beta-1}|z|^{-1})^{k}\le\frac{2q^{j_{n}^{2}/2+(\beta-1)j_{n}}}{\left|z\right|^{j_{n}}}\sum_{k=0}^{\infty}q^{k^{2}/2}(q^{\beta-1+j_{n}}|z|^{-1})^{k}\label{eq:173}\\
 & \le2\theta_{3}\left(|z|^{-1}q^{\beta};\sqrt{q}\right)\frac{q^{j_{n}^{2}/2+(\beta-1)j_{n}}}{\left|z\right|^{j_{n}}},\nonumber \end{align}
and \begin{equation}
|s_{12}|\le\frac{2q^{k_{n}}}{1-q}\sum_{k=0}^{\infty}q^{k^{2}/2}(q^{\beta}|z|^{-1})^{k}\le2\theta_{3}(|z|^{-1}q^{\beta};\sqrt{q})\frac{q^{k_{n}}}{1-q},\label{eq:174}\end{equation}
\begin{equation}
|s_{13}|\le\frac{3\log n}{n^{\rho}}\sum_{k=0}^{\infty}q^{k^{2}/2}(q^{\beta}|z|^{-1})^{k}\le3\theta_{3}(|z|^{-1}q^{\beta};\sqrt{q})\frac{\log n}{n^{\rho}}\label{eq:175}\end{equation}
for $n$ sufficiently large. Hence,

\begin{equation}
\frac{(q;q)_{\infty}s_{1}}{(-ze^{2\pi in\theta})^{m}q^{m(n\tau+m/2)}}=\sum_{k=0}^{\infty}q^{k^{2}/2}(-q^{\beta}z^{-1}e^{-2\pi i\lambda})^{k}+r_{1}(n),\label{eq:176}\end{equation}
 where\begin{equation}
r_{1}(n)=s_{11}+s_{12}+s_{13}+s_{14},\label{eq:177}\end{equation}
and\begin{equation}
|r_{1}(n)|\le3\theta_{3}(|z|^{-1}q^{\beta};\sqrt{q})\left\{ \frac{q^{k_{n}}}{1-q}+\frac{\log n}{n^{\rho}}+\frac{q^{j_{n}^{2}/2+(\beta-1)j_{n}}}{\left|z\right|^{j_{n}}}\right\} \label{eq:178}\end{equation}
 for $n$ sufficiently large. Similarly,

\begin{align}
 & \frac{(q;q)_{\infty}s_{2}}{(-ze^{2\pi in\theta})^{m}q^{m(n\tau+m/2)}}=\sum_{k=1}^{\infty}q^{k^{2}/2}(-zq^{-\beta}e^{2\pi i\lambda})^{k}f(k,n)q^{-ka_{n}}\label{eq:179}\\
 & =\sum_{k=-1}^{-\infty}q^{k^{2}/2}(-q^{\beta}z^{-1}e^{-2\pi i\lambda})^{k}-\sum_{k=j_{n}}^{\infty}q^{k^{2}/2}(-zq^{-\beta}e^{2\pi i\lambda})^{k}\nonumber \\
 & +\sum_{k=1}^{j_{n}-1}q^{k^{2}/2}(-zq^{-\beta}e^{2\pi i\lambda})^{k}\left\{ f(k,n)-1\right\} +\sum_{k=1}^{j_{n}-1}q^{k^{2}/2}(-zq^{-\beta}e^{2\pi i\lambda})^{k}f(k,n)\left\{ q^{-ka_{n}}-1\right\} \nonumber \\
 & +\sum_{k=j_{n}}^{\infty}q^{k^{2}/2}(-zq^{-\beta}e^{2\pi i\lambda})^{k}f(k,n)q^{-ka_{n}}\nonumber \\
 & =\sum_{k=-\infty}^{-1}q^{k^{2}/2}(-q^{\beta}z^{-1}e^{-2\pi i\lambda})^{k}+s_{21}+s_{22}+s_{23}+s_{24}.\nonumber \end{align}
Then,\begin{align}
|s_{21}+s_{24}| & \le2\sum_{k=j_{n}}^{\infty}q^{k^{2}/2}(|z|q^{-1-\beta})^{k}\le2\left|z\right|^{j_{n}}q^{j_{n}^{2}/2-(\beta+1)j_{n}}\sum_{k=0}^{\infty}q^{k^{2}/2}(|z|q^{-1-\beta+j_{n}})^{k}\label{eq:180}\\
 & \le2\theta_{3}(|z|^{-1}q^{\beta};\sqrt{q})\left|z\right|^{j_{n}}q^{j_{n}^{2}/2-(\beta+1)j_{n}},\nonumber \end{align}
\begin{equation}
|s_{22}|\le\frac{2q^{k_{n}}}{1-q}\sum_{k=1}^{\infty}q^{k^{2}/2}(|z|q^{-\beta})^{k}\le\frac{2q^{k_{n}}}{1-q}\theta_{3}(|z|^{-1}q^{\beta};\sqrt{q}),\label{eq:181}\end{equation}
and\begin{equation}
|s_{23}|\le\frac{3\log n}{n^{\rho}}\sum_{k=0}^{\infty}q^{k^{2}/2}(|z|q^{-\beta})^{k}\le\frac{3\log n}{n^{\rho}}\theta_{3}(|z|^{-1}q^{\beta};\sqrt{q})\label{eq:182}\end{equation}
for $n$ sufficiently large. Hence\begin{equation}
\frac{(q;q)_{\infty}s_{2}}{(-ze^{2\pi in\theta})^{m}q^{m(n\tau+m/2)}}=\sum_{k=-\infty}^{-1}q^{k^{2}/2}(-q^{\beta}z^{-1}e^{-2\pi i\lambda})^{k}+r_{2}(n),\label{eq:183}\end{equation}
 where\begin{equation}
r_{2}(n)=s_{21}+s_{22}+s_{23}+s_{24},\label{eq:184}\end{equation}
 and\begin{equation}
|r_{2}(n)|\le3\theta_{3}(|z|^{-1}q^{\beta};\sqrt{q})\left\{ \frac{\log n}{n^{\rho}}+\frac{q^{k_{n}}}{1-q}+\left|z\right|^{j_{n}}q^{j_{n}^{2}/2-(\beta+1)j_{n}}\right\} \label{eq:185}\end{equation}
 for $n$ sufficiently large. Therefore,

\begin{equation}
\frac{E_{q}(-q^{ns+1/2}z)(q;q)_{\infty}}{(-ze^{2\pi in\theta})^{m}q^{m(n\tau+m/2)}}=\theta_{4}\left(z^{-1}q^{\beta}e^{-2\pi i\lambda};\sqrt{q}\right)+r_{qe}(n|6),\label{eq:186}\end{equation}
 with\begin{align}
|r_{qe}(n|6)| & \le6\theta_{3}(|z|^{-1}q^{\beta};\sqrt{q})\left\{ \frac{\log n}{n^{\rho}}+\frac{q^{k_{n}}}{1-q}+\left|z\right|^{j_{n}}q^{j_{n}^{2}/2-(\beta+1)j_{n}}+\frac{q^{j_{n}^{2}/2+(\beta-1)j_{n}}}{\left|z\right|^{j_{n}}}\right\} \label{eq:187}\end{align}
 for $n$ sufficiently large.

\subsubsection{Proof for case 7}

Observe that \begin{align}
 & \frac{(q;q)_{\infty}s_{1}}{(-ze^{2\pi in\theta})^{m}q^{m(n\tau+m/2)}}=\sum_{k=0}^{m}q^{k^{2}/2}(-q^{\beta_{1}}z^{-1}e^{-2\pi i\beta_{2}})^{k}q^{ka_{n}}e^{-2\pi ikb_{n}}e(k,n)\label{eq:188}\\
 & =\sum_{k=0}^{\infty}q^{k^{2}/2}(-q^{\beta_{1}}z^{-1}e^{-2\pi i\beta_{2}})^{k}-\sum_{k=j_{n}}^{\infty}q^{k^{2}/2}(-q^{\beta_{1}}z^{-1}e^{-2\pi i\beta_{2}})^{k}\nonumber \\
 & +\sum_{k=0}^{j_{n}-1}q^{k^{2}/2}(-q^{\beta_{1}}z^{-1}e^{-2\pi i\beta_{2}})^{k}\left\{ q^{ka_{n}}-1\right\} +\sum_{k=0}^{j_{n}-1}q^{k^{2}/2}(-q^{\beta_{1}}z^{-1}e^{-2\pi i\beta_{2}})^{k}q^{ka_{n}}\left\{ e(k,n)-1\right\} \nonumber \\
 & +\sum_{k=0}^{j_{n}-1}q^{k^{2}/2}(-q^{\beta_{1}}z^{-1}e^{-2\pi i\beta_{2}})^{k}q^{ka_{n}}\left\{ e^{-2\pi ikb_{n}}-1\right\} e(k,n)\nonumber \\
 & +\sum_{k=j_{n}}^{m}q^{k^{2}/2}(-q^{\beta_{1}}z^{-1}e^{-2\pi i\beta_{2}})^{k}q^{ka_{n}}e^{-2\pi ikb_{n}}e(k,n)\nonumber \\
 & =\sum_{k=0}^{\infty}q^{k^{2}/2}(-q^{\beta_{1}}z^{-1}e^{-2\pi i\beta_{2}})^{k}+s_{11}+s_{12}+s_{13}+s_{14}+s_{15}.\nonumber \end{align}
For $n$ sufficiently large we have

\begin{eqnarray}
|s_{11}+s_{15}| & \le & 2\sum_{k=j_{n}}^{\infty}q^{k^{2}/2}(q^{\beta_{1}-1}|z|^{-1})^{k}\le\frac{2q^{j_{n}^{2}/2+(\beta_{1}-1)j_{n}}}{|z|^{j_{n}}}\sum_{k=0}^{\infty}q^{k^{2}/2}(q^{\beta_{1}-1+j_{n}}|z|^{-1})^{k}\label{eq:189}\\
 & \le & \frac{2q^{j_{n}^{2}/2+(\beta_{1}-1)j_{n}}}{|z|^{j_{n}}}\theta_{3}(|z|^{-1}q^{\beta_{1}};\sqrt{q}),\nonumber \end{eqnarray}
 \begin{equation}
|s_{12}|\le\frac{3\log n}{n^{\rho}}\sum_{k=0}^{\infty}q^{k^{2}/2}(q^{\beta_{1}}|z|^{-1})^{k}\le\frac{3\log n}{n^{\rho}}\theta_{3}(|z|^{-1}q^{\beta_{1}};\sqrt{q}),\label{eq:190}\end{equation}
\begin{equation}
|s_{13}|\le\frac{2q^{k_{n}}}{1-q}\sum_{k=0}^{\infty}q^{k^{2}/2}(q^{\beta_{1}}|z|^{-1})^{k}\le\frac{2q^{k_{n}}}{1-q}\theta_{3}(|z|^{-1}q^{\beta_{1}};\sqrt{q}),\label{eq:191}\end{equation}
 and\begin{equation}
|s_{14}|\le\frac{24\log n}{n^{\rho}}\theta_{3}(|z|^{-1}q^{\beta_{1}};\sqrt{q}).\label{eq:192}\end{equation}
Hence,\begin{equation}
\frac{(q;q)_{\infty}s_{1}}{(-ze^{2\pi in\theta})^{m}q^{m(n\tau+m/2)}}=\sum_{k=0}^{\infty}q^{k^{2}/2}(-q^{\beta_{1}}z^{-1}e^{-2\pi i\beta_{2}})^{k}+r_{1}(n),\label{eq:193}\end{equation}
where,\begin{equation}
r_{1}(n)=s_{11}+s_{12}+s_{13}+s_{14}+s_{15},\label{eq:194}\end{equation}
and\begin{equation}
|r_{1}(n)|\le27\theta_{3}(|z|^{-1}q^{\beta_{1}};\sqrt{q})\left\{ \frac{q^{j_{n}^{2}/2+(\beta_{1}-1)j_{n}}}{|z|^{j_{n}}}+\frac{q^{k_{n}}}{1-q}+\frac{\log n}{n^{\rho}}\right\} \label{eq:195}\end{equation}
for $n$ sufficiently large.

Similarly,

\begin{align}
 & \frac{(q;q)_{\infty}s_{2}}{(-ze^{2\pi in\theta})^{m}q^{m(n\tau+m/2)}}=\sum_{k=1}^{\infty}q^{k^{2}/2}(-zq^{-\beta_{1}}e^{2\pi i\beta_{2}})^{k}f(k,n)e^{2\pi ikb_{n}}q^{-ka_{n}}\label{eq:196}\\
 & =\sum_{k=-\infty}^{-1}q^{k^{2}/2}(-q^{\beta_{1}}z^{-1}e^{-2\pi i\beta_{2}})^{k}-\sum_{k=j_{n}}^{\infty}q^{k^{2}/2}(-zq^{-\beta_{1}}e^{2\pi i\beta_{2}})^{k}\nonumber \\
 & +\sum_{k=1}^{j_{n}-1}q^{k^{2}/2}(-zq^{-\beta_{1}}e^{2\pi i\beta_{2}})^{k}\left\{ f(k,n)-1\right\} +\sum_{k=1}^{j_{n}-1}q^{k^{2}/2}(-zq^{-\beta_{1}}e^{2\pi i\beta_{2}})^{k}f(k,n)\left\{ q^{-ka_{n}}-1\right\} \nonumber \\
 & +\sum_{k=1}^{j_{n}-1}q^{k^{2}/2}(-zq^{-\beta_{1}}e^{2\pi i\beta_{2}})^{k}f(k,n)q^{-ka_{n}}\left\{ e^{2\pi ikb_{n}}-1\right\} \nonumber \\
 & +\sum_{k=j_{n}}^{\infty}q^{k^{2}/2}(-zq^{-\beta_{1}}e^{2\pi i\beta_{2}})^{k}f(k,n)e^{2\pi ikb_{n}}q^{-ka_{n}}\nonumber \\
 & =\sum_{k=-\infty}^{-1}q^{k^{2}/2}(-q^{\beta_{1}}z^{-1}e^{-2\pi i\beta_{2}})^{k}+s_{21}+s_{22}+s_{23}+s_{24}+s_{25}.\nonumber \end{align}

For $n$ sufficiently large we have

\begin{eqnarray}
|s_{21}+s_{25}| & \le & 2\sum_{k=j_{n}}^{\infty}q^{k^{2}/2}(|z|q^{-\beta_{1}-1})^{k}\le2q^{j_{n}^{2}/2-(\beta_{1}+1)j_{n}}|z|^{j_{n}}\sum_{k=0}^{\infty}q^{k^{2}/2}(|z|q^{j_{n}-\beta_{1}-1})^{k}\label{eq:197}\\
 & \le & 2q^{j_{n}^{2}/2-(\beta_{1}+1)j_{n}}|z|^{j_{n}}\theta_{3}(|z|^{-1}q^{\beta_{1}};\sqrt{q}),\nonumber \end{eqnarray}
\begin{equation}
|s_{22}|\le\frac{2q^{k_{n}}}{1-q}\sum_{k=1}^{\infty}q^{k^{2}/2}(|z|q^{-\beta_{1}})^{k}\le\frac{2q^{k_{n}}}{1-q}\theta_{3}(|z|^{-1}q^{\beta_{1}};\sqrt{q}),\label{eq:198}\end{equation}
\begin{equation}
|s_{23}|\le\frac{3\log n}{n^{\rho}}\sum_{k=1}^{\infty}q^{k^{2}/2}(|z|q^{-\beta_{1}})^{k}\le\frac{3\log n}{n^{\rho}}\theta_{3}(|z|^{-1}q^{\beta_{1}};\sqrt{q})\label{eq:199}\end{equation}
and\begin{eqnarray}
|s_{24}| & \le & \frac{24\log n}{n^{\rho}}\theta_{3}(|z|^{-1}q^{\beta_{1}};\sqrt{q}).\label{eq:200}\end{eqnarray}
Hence,\begin{equation}
\frac{(q;q)_{\infty}s_{2}}{(-ze^{2\pi in\theta})^{m}q^{m(n\tau+m/2)}}=\sum_{k=-\infty}^{-1}q^{k^{2}/2}(-q^{\beta_{1}}z^{-1}e^{-2\pi i\beta_{2}})^{k}+r_{2}(n),\label{eq:201}\end{equation}
 where,

\begin{equation}
r_{2}(n)=s_{21}+s_{22}+s_{23}+s_{24}+s_{25},\label{eq:202}\end{equation}
and\begin{equation}
|r_{2}(n)|\le27\theta_{3}(|z|^{-1}q^{\beta_{1}};\sqrt{q})\left\{ \frac{\log n}{n^{\rho}}+\frac{q^{k_{n}}}{1-q}+q^{j_{n}^{2}/2-(\beta_{1}+1)j_{n}}|z|^{j_{n}}\right\} .\label{eq:203}\end{equation}
Therefore,\begin{equation}
\frac{E_{q}(-q^{ns+1/2}z)(q;q)_{\infty}}{(-ze^{2\pi in\theta})^{m}q^{m(n\tau+m/2)}}=\theta_{4}\left(z^{-1}q^{\beta_{1}}e^{-2\pi i\beta_{2}};\sqrt{q}\right)+r_{qe}(n|7),\label{eq:204}\end{equation}
where,\begin{equation}
r_{qe}(n|7)=r_{1}(n)+r_{2}(n),\label{eq:205}\end{equation}
 with\begin{equation}
|r_{qe}(n|7)|\le54\theta_{3}(|z|^{-1}q^{\beta_{1}};\sqrt{q})\left\{ \frac{\log n}{n^{\rho}}+\frac{q^{k_{n}}}{1-q}+q^{j_{n}^{2}/2-(\beta_{1}+1)j_{n}}|z|^{j_{n}}+\frac{q^{j_{n}^{2}/2+(\beta_{1}-1)j_{n}}}{|z|^{j_{n}}}\right\} \label{eq:206}\end{equation}
for $n$ sufficiently large.

\subsection{Proofs For Corollary \ref{cor:q-exponential-power} }

\subsubsection{Proof for case 1 }

We notice that \begin{equation}
\frac{q^{n\tau+1/2}}{1-q}=\mathcal{O}\left(n^{a}e^{-2\pi\tau n^{1-a}}\right)\label{eq:207}\end{equation}
 for $n$ sufficiently large. Then from \eqref{eq:92} and \eqref{eq:93}
we have,

\begin{equation}
E_{q}(-\exp2\pi(u-n^{1-a}\tau-n^{-a}/2))=1+\mathcal{O}\left(n^{a}e^{-2\pi\tau n^{1-a}}\right)\label{eq:208}\end{equation}
as $n\to\infty$.

\subsubsection{Proof for case 2}

From \eqref{eq:91}, \pageref{eq:97} and \eqref{eq:98} we have\begin{equation}
\left(\frac{q^{k_{n}}}{1-q}+\frac{q^{k_{n}^{2}/2+\lambda k_{n}}}{\left|z\right|^{k_{n}}}\right)e^{\pi n^{a}/4}=\mathcal{O}\left(e^{-2\pi n^{a}}\right)\label{eq:209}\end{equation}
as $n\to\infty$. From \eqref{eq:31} and \eqref{eq:32} we have\begin{align}
\theta_{3}\left(z^{-1}q^{\lambda};\sqrt{q}\right) & =\theta_{3}(n^{-a}\lambda i+ui|n^{-a}i)\label{eq:210}\\
 & =\sqrt{n^{a}}e^{\pi n^{-a}(n^{a}u+\lambda)^{2}}\theta_{3}(n^{a}u+\lambda|n^{a}i)\nonumber \\
 & =\sqrt{n^{a}}e^{\pi n^{-a}(n^{a}u+\lambda)^{2}}\left\{ 1+\mathcal{O}\left(e^{-\pi n^{a}}\right)\right\} ,\nonumber \end{align}
and\begin{align}
\theta_{4}\left(z^{-1}q^{\lambda};\sqrt{q}\right) & =\theta_{4}(n^{-a}\lambda i+ui|n^{-a}i)\label{eq:211}\\
 & =\sqrt{n^{a}}e^{\pi n^{-a}(n^{a}u+\lambda)^{2}}\theta_{2}(n^{a}u+\lambda|n^{a}i)\nonumber \\
 & =2\sqrt{n^{a}}\exp\left(\pi n^{-a}(n^{a}u+\lambda)^{2}-\frac{\pi n^{a}}{4}\right)\nonumber \\
 & \cos\pi(n^{a}u+\lambda)\left\{ 1+\mathcal{O}\left(e^{-2\pi n^{a}}\right)\right\} ,\nonumber \end{align}
 and\begin{equation}
\frac{1}{(q;q)_{\infty}}=\frac{\exp\left(\pi n^{a}/12-\pi n^{-a}/12\right)}{\sqrt{n^{a}}}\left\{ 1+\mathcal{O}\left(e^{-2\pi n^{a}}\right)\right\} \label{eq:212}\end{equation}
 for $n$ sufficiently large.

Then,\begin{align}
 & E_{q}(-\exp2\pi(u-n^{1-a}\tau-n^{-a}/2))\label{eq:213}\\
 & =\frac{2\exp\left(\pi n^{-a}(n^{a}u-n\tau)^{2}\right)\left\{ \cos\pi(n^{a}u+\lambda)+\mathcal{O}\left(e^{-2\pi n^{a}}\right)\right\} }{(-1)^{\lambda+n\tau}\exp\left(\pi n^{a}/6+\pi n^{-a}/12\right)}\nonumber \end{align}
 for $n$ sufficiently large.

\subsection{Proof for Corollary \ref{cor:q-exponential-log}}

\subsubsection{Proof for case 1.}

Since\begin{equation}
q^{n\tau+1/2}=\mathcal{O}\left(\log n\exp\left(-\frac{2n\pi\tau}{\gamma\log n}\right)\right)\label{eq:214}\end{equation}
 as $n\to\infty$, then

\begin{equation}
E_{q}\left(-\exp2\pi\left(u-\frac{n\tau+1/2}{\gamma\log n}\right)\right)=1+\mathcal{O}\left(\log n\exp\left(-\frac{2n\pi\tau}{\gamma\log n}\right)\right)\label{eq:215}\end{equation}
 for $n$ sufficiently large.

\subsubsection{Proof for case 2.}

In this case we have\begin{align}
\theta_{3}(z^{-1}q^{\lambda};\sqrt{q}) & =\theta_{3}\left(ui+\frac{\lambda i}{\gamma\log n}\mid\frac{i}{\gamma\log n}\right)\label{eq:216}\\
 & =\sqrt{\gamma\log n}\exp\frac{\left(\gamma u\log n+\lambda\right)^{2}}{\gamma\log n}\theta_{3}\left(\gamma u\log n+\lambda\mid i\gamma\log n\right)\nonumber \\
 & =\sqrt{\gamma\log n}\exp\frac{\left(\gamma u\log n+\lambda\right)^{2}}{\gamma\log n}\left\{ 1+\mathcal{O}(n^{-\pi\gamma})\right\} ,\nonumber \end{align}
 and\begin{align}
\theta_{4}(z^{-1}q^{\lambda};\sqrt{q}) & =\theta_{4}\left(ui+\frac{\lambda i}{\gamma\log n}\mid\frac{i}{\gamma\log n}\right)\label{eq:217}\\
 & =\sqrt{\gamma\log n}\exp\frac{\left(\gamma u\log n+\lambda\right)^{2}}{\gamma\log n}\theta_{2}\left(\gamma u\log n+\lambda\mid i\gamma\log n\right)\nonumber \\
 & =2n^{-\pi\gamma/4}\sqrt{\gamma\log n}\exp\frac{\left(\gamma u\log n+\lambda\right)^{2}}{\gamma\log n}\nonumber \\
 & \times\cos\pi(\gamma u\log n+\lambda)\left\{ 1+\mathcal{O}(n^{-2\pi\gamma})\right\} ,\nonumber \end{align}
and\begin{equation}
\frac{1}{(q;q)_{\infty}}=\frac{n^{\gamma\pi/12}\exp\left(-\frac{\pi}{12\gamma\log n}\right)}{\sqrt{\gamma\log n}}\left\{ 1+\mathcal{O}(n^{-2\pi\gamma})\right\} \label{eq:218}\end{equation}
as $n\to\infty$. Thus,\begin{align}
E_{q}\left(-\exp2\pi\left(u-\frac{n\tau+1/2}{\gamma\log n}\right)\right) & =\exp\left(\frac{\pi(\gamma u\log n-n\tau)^{2}}{\gamma\log n}-\frac{\pi}{12\gamma\log n}\right)\label{eq:219}\\
 & \times\frac{2\left\{ \cos\pi(\gamma u\log n+\lambda)+\mathcal{O}(n^{-2\pi\gamma})\right\} }{(-1)^{\tau n+\lambda}n^{\gamma\pi/6}}\nonumber \end{align}
 as $n\to\infty$.

\subsubsection{Proof for case 3.}

A similar calculation to get\begin{align}
E_{q}\left(-\exp2\pi\left(u-\frac{n\tau+1/2}{\gamma\log n}\right)\right) & =\exp\left(\frac{\pi(\gamma u\log n-n\tau)^{2}}{\gamma\log n}-\frac{\pi}{12\gamma\log n}\right)\label{eq:220}\\
 & \times\frac{2\left\{ \cos\pi(\gamma u\log n+\beta)+\mathcal{O}(n^{-8\rho/9}\log n)\right\} }{(-1)^{m}n^{2\rho/27}}\nonumber \end{align}
 as $n\to\infty$, where \begin{equation}
\gamma=\frac{4\rho}{9\pi}.\label{eq:221}\end{equation}

\section{Ramanujan Functions}

Ramanujan function $A_{q}(z)$ \eqref{eq:1.3}, which is also called
$q$-Airy function in the literature, appears repeatedly in Ramanujan's
work starting from the Rogers-Ramanujan identities, where $A_{q}(-1)$
and $A_{q}(-q)$ are expressed as infinite products, \cite{Andrews1},
to properties of and conjectures about its zeros, \cite{Andrews3,Andrews4,Hayman,Ismail5}.
It is called $q$-Airy function because it appears repeatedly in the
Plancherel-Rotach type asymptotics of $q$-orthogonal polynomials,
just like the classical Airy function in the classical Plancherel-Rotach
asymptotics of classical orthogonal polynomials \cite{Szego,Ismail2,Deift1,Deift2}. 

Apply  inequalities in \eqref{eq:13}, we have\begin{equation}
\left|A'_{q}(z)\right|\le\frac{q}{1-q}A_{q}(-|z|),\quad z\in\mathbb{C},\label{eq:222}\end{equation}
and \begin{equation}
\left|A_{q}((1-q)z)\right|\le e^{q|z|},\quad z\in\mathbb{C},\label{eq:223}\end{equation}
 or\begin{equation}
\left|A_{q}(z)\right|\le e^{q|z|/(1-q)},\quad z\in\mathbb{C}.\label{eq:224}\end{equation}
 Similar to $E_{q}(z)$ case, we could prove that\begin{equation}
\lim_{q\to1}A_{q}((1-q)z)=e^{-z},\quad z\in\mathbb{C}.\label{eq:225}\end{equation}
Hence, $A_{q}(z)$ is really one of many $q$-analogues of the exponential
function $e^{-z}$.

For any nonzero complex number $z$, then,\begin{align}
\left|A_{q}(z)\right| & \le\sum_{k=0}^{\infty}\frac{q^{k^{2}/2}}{(q;q)_{k}}\left(q^{k/2}\left|z\right|\right)^{k},\label{eq:226}\end{align}
clearly, the terms $q^{k^{2}/2}|z|^{k}$ are bounded by \begin{equation}
\exp\left\{ -\frac{\log^{2}\left|z\right|}{2\log q}\right\} ,\label{eq:227}\end{equation}
for $k=0,1,...$, thus,

\begin{align}
\left|A_{q}(z)\right| & \le\exp\left\{ -\frac{\log^{2}\left|z\right|}{2\log q}\right\} \sum_{k=0}^{\infty}\frac{q^{k^{2}/2}}{(q;q)_{k}}\label{eq:228}\\
 & \le(-\sqrt{q};q)_{\infty}\exp\left\{ -\frac{\log^{2}\left|z\right|}{2\log q}\right\} ,\nonumber \end{align}
 or we have\begin{align}
\left|A_{q}(z)\right| & \le(-\sqrt{q};q)_{\infty}\exp\left\{ -\frac{\log^{2}\left|z\right|}{2\log q}\right\} \label{eq:229}\end{align}
 for any nonzero complex number $z$.

\subsection{Asymptotic Formulas for $A_{q}(z)$\label{sec:ramanujan-asymptotics}}

Assume that \begin{equation}
s=2\tau+\frac{2\pi i\theta}{\log q},\quad\tau,\theta\in\mathbb{R},\label{eq:230}\end{equation}
 then for any complex number $z$ we have\begin{equation}
A_{q}(q^{ns}z)=\sum_{k=0}^{\infty}\frac{q^{k^{2}}}{(q;q)_{k}}(-q^{2n\tau}ze^{2n\pi i\theta})^{k}.\label{eq:231}\end{equation}

\begin{thm*}
\label{thm:ramanujan}For any nonzero complex number $z$, let\begin{equation}
j_{n}=\left\lfloor \frac{q^{4}\log n}{-\log q}\right\rfloor ,\quad k_{n}=\left\lfloor \frac{-n\tau}{2}\right\rfloor ,\label{eq:232}\end{equation}
we have the following for Ramanujan function $A_{q}(z)$:
\begin{enumerate}
\item If $\tau>0$, then, \begin{equation}
A_{q}(q^{ns}z)=1+r_{rf}(n|1),\label{eq:233}\end{equation}
and\begin{equation}
|r_{rf}(n|1)|\le\frac{|z|q^{2n\tau+1}}{1-q}\exp\left(\frac{|z|q^{2n\tau+1}}{1-q}\right).\label{eq:234}\end{equation}

\item Assume that $\tau=0$, if for some fixed real numbers $\beta$ and
$\rho\ge1$ there are infinitely many positive integers $n$ such
that \begin{equation}
n\theta=m+\beta+b_{n},\quad|b_{n}|<\frac{1}{n^{\rho}},\quad m\in\mathbb{Z},\label{eq:235}\end{equation}
 then,\begin{equation}
A_{q}(q^{ns}z)=A_{q}\left(ze^{2\pi i\beta}\right)+r_{rf}(n|3),\label{eq:236}\end{equation}
and\begin{equation}
|r_{rf}(n|3)|\le24\exp\left(\frac{|z|q}{1-q}\right)\left\{ \frac{\log n}{n^{\rho}}+\frac{q^{j_{n}}|z|^{j_{n}}}{j_{n}!(1-q)^{j_{n}}}\right\} \label{eq:237}\end{equation}
for $n$ sufficiently large. 
\item Assume that $\tau<0$. If for some fixed rational numbers $\lambda$
and $\lambda_{1}$ there are infinitely many positive integers $n$
such that \begin{equation}
-n\tau=m+\lambda,\quad m\in\mathbb{N},\quad n\theta=m_{1}+\lambda_{1},\quad m_{1}\in\mathbb{Z},\label{eq:238}\end{equation}
 then,\begin{align}
\frac{A_{q}(q^{ns}z)(q;q)_{\infty}}{(-ze^{2\pi in\theta})^{m}q^{m(2n\tau+m)}} & =\theta_{4}\left(z^{-1}q^{2\lambda}e^{-2\pi i\lambda_{1}};q\right)+r_{rf}(n|4),\label{eq:239}\end{align}
and \begin{align}
|r_{rf}(n|4)| & \le4\theta_{3}\left(|z|^{-1}q^{2\lambda};q\right)\left\{ \frac{q^{k_{n}}}{1-q}+\frac{q^{k_{n}^{2}+2\lambda k_{n}}}{\left|z\right|^{k_{n}}}\right\} \label{eq:240}\end{align}
for $n$ sufficiently large.
\item Assume that $\tau<0$. If for some fixed real numbers $\beta$, $\lambda$
and $\rho\ge1$ there exists infinitely many positive integers $n$
such that\begin{equation}
-n\tau=m+\lambda,\quad m\in\mathbb{N},\quad n\theta=m_{1}+\beta+b_{n},\quad|b_{n}|<\frac{1}{n^{\rho}},\quad m_{1}\in\mathbb{Z},\label{eq:241}\end{equation}
 then, \begin{align}
\frac{A_{q}(q^{ns}z)(q;q)_{\infty}}{(-ze^{2\pi in\theta})^{m}q^{m(2n\tau+m)}} & =\theta_{4}\left(z^{-1}q^{2\lambda}e^{-2\pi i\beta};q\right)+r_{rf}(n|5),\label{eq:242}\end{align}
 and\begin{align}
|r_{rf}(n|5)| & \le48\theta_{3}\left(|z|^{-1}q^{2\lambda};q\right)\left\{ \frac{q^{j_{n}^{2}+2\lambda j_{n}}}{\left|z\right|^{j_{n}}}+\left|z\right|^{j_{n}}q^{j_{n}^{2}-2\lambda j_{n}}+\frac{q^{k_{n}}}{1-q}+\frac{\log n}{n^{\rho}}\right\} \label{eq:243}\end{align}
 for $n$ sufficiently large. 
\item Assume that $\tau<0$. If for some fixed real numbers $\lambda$,
$\beta$ and $\rho\ge1$ there exist infinitely many positive integers
$n$ such that\begin{equation}
-n\tau=m+\beta+a_{n},\quad|a_{n}|<\frac{1}{n^{\rho}},\quad m\in\mathbb{N},\quad n\theta=m_{1}+\lambda\quad m_{1}\in\mathbb{Z},\label{eq:244}\end{equation}
 then,\begin{align}
\frac{A_{q}(q^{ns}z)(q;q)_{\infty}}{(-ze^{2\pi in\theta})^{m}q^{m(2n\tau+m)}} & =\theta_{4}\left(z^{-1}q^{2\beta}e^{-2\pi i\lambda};q\right)+r_{rf}(n|6),\label{eq:245}\end{align}
 and\begin{align}
|r_{rf}(n|6)| & \le12\theta_{3}\left(|z|^{-1}q^{2\beta};q\right)\left\{ \frac{q^{j_{n}^{2}+2(\beta-1)j_{n}}}{\left|z\right|^{j_{n}}}+\left|z\right|^{j_{n}}q^{j_{n}^{2}-2(\beta+1)j_{n}}+\frac{q^{k_{n}}}{1-q}+\frac{\log n}{n^{\rho}}\right\} \label{eq:246}\end{align}
 for $n$ sufficiently large. 
\item Assume that $\tau<0$. If for some fixed real numbers $\beta_{1},\beta_{2}$
and $\rho>0$ there exists infinitely many positive integers $n$
such that\begin{equation}
-n\tau=m+\beta_{1}+a_{n},\quad m\in\mathbb{N},\quad|a_{n}|<\frac{1}{n^{\rho}},\label{eq:247}\end{equation}
 and\begin{equation}
n\theta=m_{1}+\beta_{2}+b_{n},\quad m_{1}\in\mathbb{Z},\quad|b_{n}|<\frac{1}{n^{\rho}},\label{eq:248}\end{equation}
then,\begin{align}
\frac{A_{q}(q^{ns}z)(q;q)_{\infty}}{(-ze^{2\pi in\theta})^{m}q^{m(2n\tau+m)}} & =\theta_{4}\left(z^{-1}q^{2\beta_{1}}e^{-2\pi i\beta_{2}};q\right)+r_{rf}(n|7),\label{eq:249}\end{align}
 and\begin{align}
|r_{rf}(n|7)| & \le54\theta_{3}\left(|z|^{-1}q^{2\beta_{1}};q\right)\left\{ \frac{q^{j_{n}^{2}+2(\beta_{1}-1)j_{n}}}{\left|z\right|^{j_{n}}}+\left|z\right|^{j_{n}}q^{j_{n}^{2}-2(\beta_{1}+1)j_{n}}+\frac{q^{k_{n}}}{1-q}+\frac{\log n}{n^{\rho}}\right\} \label{eq:250}\end{align}
 for $n$ sufficiently large. 
\end{enumerate}
\end{thm*}
In the following corollaries, we let\begin{equation}
z:=e^{2\pi u},\quad u\in\mathbb{R}.\label{eq:251}\end{equation}

\begin{cor}
\label{cor:ramanujan-power}Assume that \begin{equation}
0<a<\frac{1}{2},\quad n\in\mathbb{N},\quad q=\exp(-\pi n^{-a}),\label{eq:252}\end{equation}
we have the following results for $A_{q}(z)$:
\begin{enumerate}
\item Assume that $\tau>0$, we have\begin{equation}
A_{q}(\exp2\pi(u-n^{1-a}\tau))=1+\mathcal{O}\left(n^{a}\exp(-2\pi\tau n^{1-a})\right)\label{eq:253}\end{equation}
 for $n$ sufficiently large.
\item Assume that $\tau<0$, if for some fixed real number $\lambda$, there
are infinitely many $m$ and $n>0$ such that \begin{equation}
-n\tau=m+\lambda,\quad m\in\mathbb{N}.\label{eq:254}\end{equation}
Then,\begin{align}
 & A_{q}\left(\exp\left(2\pi(u-n^{1-a}\tau)\right)\right)\label{eq:255}\\
 & =\frac{\sqrt{2}\exp\left\{ \pi n^{-a}(n^{a}u-\tau n)^{2}\right\} \left\{ \cos\pi(n^{a}u+\lambda)+\mathcal{O}(e^{-2\pi n^{a}})\right\} }{(-1)^{\tau n+\lambda}\exp\left\{ \frac{\pi n^{-a}}{24}+\frac{\pi n^{a}}{12}\right\} }\nonumber \end{align}
 as $n$ sufficiently large.
\end{enumerate}
\end{cor}
Similarly, we have the following:

\begin{cor}
\label{cor:ramanujan-log}Assume that\begin{equation}
q=\exp(-\frac{\pi}{\gamma\log n}),\quad n\ge2,\quad\gamma>0,\label{eq:256}\end{equation}
we have the following results for $A_{q}(z)$:
\begin{enumerate}
\item Assume that $\tau>0$, we have\begin{equation}
A_{q}\left(\exp2\pi\left(u-\frac{n\tau}{\gamma\log n}\right)\right)=1+\mathcal{O}\left(\log n\exp\left(-\frac{2\pi n\tau}{\gamma\log n}\right)\right)\label{eq:257}\end{equation}
 for $n$ sufficiently large.
\item Assume that $\tau<0$ and for some fixed real number $\lambda$, there
are infinitely many $m$ and $n>0$ such that \begin{equation}
-n\tau=m+\lambda,\quad m\in\mathbb{N}.\label{eq:258}\end{equation}
Then,\begin{align}
 & A_{q}\left(\exp2\pi\left(u-\frac{n\tau}{\gamma\log n}\right)\right)\label{eq:259}\\
 & =\sqrt{2}\frac{\exp\left(\frac{\pi(u\gamma\log n-\tau n)^{2}}{\gamma\log n}\right)\left\{ \cos\pi(u\gamma\log n+\lambda)+\mathcal{O}(n^{-2\pi\gamma})\right\} }{(-1)^{n\tau+\lambda}n^{\pi\gamma/12}\exp\left(\frac{\pi}{24\gamma\log n}\right)}\nonumber \end{align}
for $n$ sufficiently large.
\item Assume that $\tau<0$, if for some fixed real numbers $\beta$, $\rho>0$
and $\lambda$ there exist infinitely many positive integers $n$
such that\begin{equation}
-n\tau=m+\beta+a_{n},\quad|a_{n}|<\frac{1}{n^{\rho}},\quad m\in\mathbb{N},\label{eq:260}\end{equation}
then for each of such $n$, we have\begin{align}
 & A_{q}\left(\exp2\pi\left(u-\frac{n\tau}{\gamma\log n}\right)\right)\label{eq:261}\\
 & =\sqrt{2}\frac{\exp\left(\frac{\pi(u\gamma\log n-\tau n)^{2}}{\gamma\log n}\right)\left\{ \cos\pi(u\gamma\log n+\beta)+\mathcal{O}(n^{-8\rho/9}\log n)\right\} }{(-1)^{m}n^{\rho/27}\exp\left(\frac{\pi}{24\gamma\log n}\right)}\nonumber \end{align}
for $n$ sufficiently large, where\begin{equation}
\gamma=\frac{4\rho}{9\pi}.\label{eq:262}\end{equation}

\end{enumerate}
\end{cor}

\subsubsection{Proof for Theorem \ref{thm:q-bessel}}

Assume that $\tau<0$ and \begin{equation}
-n\tau=m+c_{n},\quad m\in\mathbb{N},\quad n\theta=m_{1}+d_{n},\quad m_{1}\in\mathbb{Z},\label{eq:263}\end{equation}
 then,\begin{align}
A_{q}(q^{ns}z) & =\sum_{k=0}^{\infty}\frac{q^{k^{2}}}{(q;q)_{k}}\left(-ze^{2\pi in\theta}q^{2n\tau}\right)^{k}\label{eq:264}\\
 & =\sum_{k=0}^{m}\frac{q^{k^{2}}}{(q;q)_{k}}\left(-ze^{2\pi in\theta}q^{2n\tau}\right)^{k}+\sum_{k=m+1}^{\infty}\frac{q^{k^{2}}}{(q;q)_{k}}\left(-ze^{2\pi in\theta}q^{2n\tau}\right)^{k}\nonumber \\
 & =s_{1}+s_{2}.\nonumber \end{align}
Reverse summation order in $s_{1}$,\begin{align}
 & \frac{(q;q)_{\infty}s_{1}}{(-ze^{2\pi in\theta})^{m}q^{m(2n\tau+m)}}=\sum_{k=0}^{m}q^{k^{2}}(-q^{2c_{n}}z^{-1}e^{-2\pi id_{n}})^{k}e(k,n),\label{eq:265}\end{align}
 where\begin{equation}
e(k,n)=\frac{(q;q)_{\infty}}{(q;q)_{m-k}},\label{eq:266}\end{equation}
 then\begin{equation}
|e(k,n)|\le1\label{eq:267}\end{equation}
 for $0\le k\le m$. From Lemma \ref{lem:1} we have \begin{equation}
|e(k,n)-1|=\left|r_{1}\left(q;m-k\right)\right|\le\frac{2q^{2+k_{n}}}{1-q}.\label{eq:268}\end{equation}
 for $0\le k\le k_{n}-1$ and $n$ sufficiently large.

We shift the summation index from $k$ to $k+m$ in $s_{2}$,\begin{align}
\frac{(q;q)_{\infty}s_{2}}{(-ze^{2\pi in\theta})^{m}q^{m(2n\tau+m)}} & =\sum_{k=1}^{\infty}q^{k^{2}}(-zq^{-2c_{n}}e^{2\pi id_{n}})^{k}f(k,n),\label{eq:269}\end{align}
where\begin{equation}
f(k,n)=\frac{(q;q)_{\infty}}{(q;q)_{m+k}},\label{eq:270}\end{equation}
 thus,\begin{equation}
|f(k,n)|\le1,\label{eq:271}\end{equation}
 and\begin{equation}
|f(k,n)-1|=\left|r_{1}\left(q;m+k\right)\right|\le\frac{2q^{2+k_{n}}}{1-q}\label{eq:272}\end{equation}
 for $k\in\mathbb{N}$ and $n$ sufficiently large. The rest are very
similar to the corresponding proofs for $E_{q}(z)$.

\subsection{Proofs for Corollary \ref{cor:ramanujan-power} }

Notice that\begin{equation}
\frac{1}{(q;q)_{\infty}}=\frac{\exp(\pi n^{a}/6-\pi n^{-a}/24)}{\sqrt{2n^{a}}}\left\{ 1+\mathcal{O}(e^{-4\pi n^{a}})\right\} ,\label{eq:273}\end{equation}
and\begin{align}
\theta_{4}\left(z^{-1}q^{2\lambda};q\right) & =\theta_{4}(ui+i\lambda n^{-a}|n^{-a}i)\label{eq:274}\\
 & =\sqrt{n^{a}}\exp\left\{ \frac{\pi(n^{a}u+\lambda)^{2}}{n^{a}}\right\} \theta_{2}\left(n^{a}u+\lambda\mid n^{a}i\right)\nonumber \\
 & =2\sqrt{n^{a}}\exp\left\{ \frac{\pi(n^{a}u+\lambda)^{2}}{n^{a}}-\frac{\pi n^{a}}{4}\right\} \nonumber \\
 & \times\cos\pi(n^{a}u+\lambda)\left\{ 1+\mathcal{O}(e^{-2\pi n^{a}})\right\} ,\nonumber \end{align}
and\begin{align}
\theta_{3}\left(|z|^{-1}q^{2\lambda};q\right) & =\theta_{3}(ui+\lambda in^{-a}|n^{-a}i)\label{eq:275}\\
 & =\sqrt{n^{a}}\exp\left\{ \frac{\pi(n^{a}u+\lambda)^{2}}{n^{a}}\right\} \theta_{3}\left(n^{a}u+\lambda\mid n^{a}i\right)\nonumber \\
 & =\sqrt{n^{a}}\exp\left\{ \frac{\pi(n^{a}u+\lambda)^{2}}{n^{a}}\right\} \left\{ 1+\mathcal{O}(e^{-\pi n^{a}})\right\} ,\nonumber \end{align}
The rest of the proof for this corollary are very similar to the $E_{q}(z)$
case.

\subsection{Proof for Corollary \ref{cor:ramanujan-log} }

Observe that\begin{equation}
\frac{1}{(q;q)_{\infty}}=\frac{n^{\pi\gamma/6}\exp\left(-\frac{\pi}{24\gamma\log n}\right)}{\sqrt{2\gamma\log n}}\left\{ 1+\mathcal{O}(n^{-4\pi\gamma})\right\} ,\label{eq:276}\end{equation}
 and\begin{align}
\theta_{3}(z^{-1}q^{2\lambda};q) & =\theta_{3}\left(ui+i\frac{\lambda}{\gamma\log n}\mid\frac{i}{\gamma\log n}\right)\label{eq:277}\\
 & =\sqrt{\gamma\log n}\exp\left(\frac{\pi(u\gamma\log n+\lambda)^{2}}{\gamma\log n}\right)\theta_{3}(u\gamma\log n+\lambda\mid i\gamma\log n)\nonumber \\
 & =\sqrt{\gamma\log n}\exp\left(\frac{\pi(u\gamma\log n+\lambda)^{2}}{\gamma\log n}\right)\left\{ 1+\mathcal{O}(n^{-\pi\gamma})\right\} ,\nonumber \end{align}
and\begin{align}
\theta_{4}(z^{-1}q^{2\lambda};q) & =\theta_{4}\left(ui+i\frac{\lambda}{\gamma\log n}\mid\frac{i}{\gamma\log n}\right)\label{eq:278}\\
 & =\sqrt{\gamma\log n}\exp\left(\frac{\pi(u\gamma\log n+\lambda)^{2}}{\gamma\log n}\right)\theta_{2}(u\gamma\log n+\lambda\mid i\gamma\log n)\nonumber \\
 & =2n^{-\pi\gamma/4}\sqrt{\gamma\log n}\exp\left(\frac{\pi(u\gamma\log n+\lambda)^{2}}{\gamma\log n}\right)\nonumber \\
 & \times\cos\pi(u\gamma\log n+\lambda)\left\{ 1+\mathcal{O}(n^{-2\pi\gamma})\right\} \nonumber \end{align}
as $n\to\infty$. The rest of proof is similar to the proof for $E_{q}(z)$
case.

\section{$q$-Bessel Function of Second Kind}

Jackson's $q$-Bessel function of second kind $J_{\nu}^{(2)}(z;q)$
is defined as \cite{Ramanujan,Ismail2,Andrews1,Gasper,Koekoek} \begin{equation}
J_{\nu}^{(2)}(z;q):=\frac{(q^{\nu+1};q)_{\infty}}{(q;q)_{\infty}}\sum_{k=0}^{\infty}\frac{q^{k^{2}+k\nu}(-1)^{k}}{(q,q^{\nu+1};q)_{k}}\left(\frac{z}{2}\right)^{2k+\nu}.\label{eq:279}\end{equation}
 In this chapter, we always assume that $\nu>-1$. Clearly, \begin{equation}
\left|J_{\nu}^{(2)}(z;q)\right|\le\frac{(|z|/2)^{\nu}}{(q;q)_{\infty}}A_{q}\left(-\frac{|z|^{2}q^{\nu}}{4}\right),\label{eq:280}\end{equation}
 and\begin{equation}
\left|J_{\nu}^{(2)}(z;q)\right|\le\frac{(|z|/2)^{\nu}}{(q;q)_{\infty}}\exp\left\{ \frac{q^{1+\nu}|z|^{2}}{4(1-q)}\right\} .\label{eq:281}\end{equation}
for any complex number $z$. 

For any nonzero complex number $z$, \eqref{eq:223} implies\begin{equation}
\left|J_{\nu}^{(2)}(z;q)\right|\le\frac{(-\sqrt{q};q)_{\infty}}{(q;q)_{\infty}}\left(\frac{|z|}{2}\right)^{\nu}\exp\left\{ -\frac{\log^{2}\left(\left|z\right|^{2}q^{\nu}/4\right)}{2\log q}\right\} .\label{eq:282}\end{equation}

\subsection{Asymptotics for $J_{\nu}^{(2)}(z;q)$ \label{sec:q-bessel asymptotics}}

From \eqref{eq:267} we have\begin{equation}
\frac{J_{\nu}^{(2)}(2q^{ns-\nu/2}z;q)(q;q)_{\infty}}{(q^{\nu+1};q)_{\infty}\left(q^{ns-\nu/2}z\right)^{\nu}}=\sum_{k=0}^{\infty}\frac{q^{k^{2}}\left(-z^{2}q^{2ns}\right)^{k}}{(q,q^{\nu+1};q)_{k}}.\label{eq:283}\end{equation}
Let \begin{equation}
s=\tau+\frac{\pi i\theta}{\log q},\quad\tau,\theta\in\mathbb{R},\label{eq:284}\end{equation}
 then for any complex number $z$ we have\begin{equation}
\frac{J_{\nu}^{(2)}(2q^{ns-\nu/2}z;q)(q;q)_{\infty}}{(q^{\nu+1};q)_{\infty}\left(q^{n\tau-\nu/2}ze^{n\pi i\theta}\right)^{\nu}}=\sum_{k=0}^{\infty}\frac{q^{k^{2}}(-q^{2n\tau}z^{2}e^{2n\pi i\theta})^{k}}{(q,q^{\nu+1};q)_{k}}.\label{eq:285}\end{equation}

\begin{thm*}
\label{thm:q-bessel}For any nonzero complex number $z$, let\begin{equation}
j_{n}=\left\lfloor \frac{q^{4}\log n}{-\log q}\right\rfloor ,\quad k_{n}=\left\lfloor \frac{-n\tau}{2}\right\rfloor ,\label{eq:286}\end{equation}
we have the following for the function $J_{\nu}^{(2)}(z;q)$:
\begin{enumerate}
\item When $\tau>0$ we have \begin{equation}
\frac{J_{\nu}^{(2)}(2q^{ns-\nu/2}z;q)(q;q)_{\infty}}{(q^{\nu+1};q)_{\infty}\left(q^{ns-\nu/2}z\right)^{\nu}}=1+r_{qb}(n|1),\label{eq:287}\end{equation}
and\begin{equation}
|r_{qb}(n|1)|\le\frac{|z|^{2}\exp\left(|z|^{2}q^{2n\tau+1}/(1-q)\right)}{(1-q)(q^{\nu+1};q)_{\infty}}q^{2n\tau+1}.\label{eq:288}\end{equation}

\item Assume that $\tau=0$. If for some fixed real numbers $\beta$ and
$\rho\ge1$ there are infinitely many positive integers $n$ such
that \begin{equation}
n\theta=m+\beta+b_{n},\quad|b_{n}|<\frac{1}{n^{\rho}},\quad m\in\mathbb{Z},\label{eq:289}\end{equation}
 then\begin{equation}
\frac{J_{\nu}^{(2)}(2q^{ns-\nu/2}z;q)(q;q)_{\infty}}{\left(q^{n\tau-\nu/2}ze^{n\pi i\theta}\right)^{\nu}}=\frac{J_{\nu}^{(2)}(2q^{-\nu/2}ze^{\pi i(m+\beta)};q)(q;q)_{\infty}}{\left(q^{-\nu/2}ze^{\pi i(m+\beta)}\right)^{\nu}}+r_{qb}(n|3),\label{eq:290}\end{equation}
and\begin{equation}
r_{qb}(n|3)|\le24\exp(q|z|^{2}/(1-q))\left\{ \frac{\log n}{n^{\rho}}+\frac{1}{j_{n}!}\left(\frac{q|z|^{2}}{1-q}\right)^{j_{n}}\right\} \label{eq:291}\end{equation}
for $n$ sufficiently large. 
\item Assume that $\tau<0$. If for some fixed real numbers $\lambda$ and
$\lambda_{1}$ there are infinitely many positive integers $n$ such
that \begin{equation}
-n\tau=m+\lambda,\quad m\in\mathbb{N},\quad n\theta=m_{1}+\lambda_{1},\quad m_{1}\in\mathbb{Z},\label{eq:292}\end{equation}
 then,\begin{align}
\frac{J_{\nu}^{(2)}(2q^{ns-\nu/2}z;q)(q;q)_{\infty}^{2}e^{-n\pi i\theta(2m+\nu)}}{(-1)^{m}z^{\nu+2m}q^{m(2n\tau+m)+n\nu\tau-\nu^{2}/2}} & =\theta_{4}\left(z^{-2}q^{2\lambda}e^{-2\pi i\lambda_{1}};q\right)+r_{qb}(n|4),\label{eq:293}\end{align}
and \begin{align}
|r_{qb}(n|4)| & \le12\theta_{3}\left(|z|^{-2}q^{2\lambda};q\right)\left\{ \frac{q^{k_{n}}}{(1-q)}+\frac{q^{k_{n}^{2}+2\lambda k_{n}}}{\left|z\right|^{2k_{n}}}\right\} \label{eq:294}\end{align}
 for $n$ sufficiently large.
\item Assume that $\tau<0$. If for some fixed real numbers $\beta$, $\lambda$
and $\rho\ge1$ there exist infinitely many positive integers $n$
such that\begin{equation}
-n\tau=m+\lambda,\quad m\in\mathbb{N},\quad n\theta=m_{1}+\beta+b_{n},\quad|b_{n}|<\frac{1}{n^{\rho}},\quad m_{1}\in\mathbb{Z},\label{eq:295}\end{equation}
 then,\begin{align}
\frac{J_{\nu}^{(2)}(2q^{ns-\nu/2}z;q)(q;q)_{\infty}^{2}e^{-n\pi i\theta(2m+\nu)}}{(-1)^{m}z^{\nu+2m}q^{m(2n\tau+m)+n\nu\tau-\nu^{2}/2}} & =\theta_{4}\left(z^{-2}q^{2\lambda}e^{-2\pi i\beta};q\right)+r_{rf}(n|5),\label{eq:296}\end{align}
 and\begin{align}
|r_{rf}(n|5)| & \le48\theta_{3}\left(|z|^{-2}q^{2\lambda};q\right)\left\{ \frac{q^{j_{n}^{2}+2\lambda j_{n}}}{\left|z\right|^{2j_{n}}}+\left|z\right|^{2j_{n}}q^{j_{n}^{2}-2\lambda j_{n}}+\frac{q^{k_{n}}}{(1-q)}+\frac{\log n}{n^{\rho}}\right\} \label{eq:297}\end{align}
 for $n$ sufficiently large.
\item Assume that $\tau<0$. If for some fixed real numbers $\beta$, $\lambda$
and $\rho\ge1$ there exist infinitely many positive integers $n$
such that\begin{equation}
-n\tau=m+\beta+a_{n},\quad m\in\mathbb{N},\quad|a_{n}|<\frac{1}{n^{\rho}},\quad n\theta=m_{1}+\lambda\quad m_{1}\in\mathbb{Z},\label{eq:298}\end{equation}
 then,\begin{align}
\frac{J_{\nu}^{(2)}(2q^{ns-\nu/2}z;q)(q;q)_{\infty}^{2}e^{-n\pi i\theta(2m+\nu)}}{(-1)^{m}z^{\nu+2m}q^{m(2n\tau+m)+n\nu\tau-\nu^{2}/2}} & =\theta_{4}(q^{2\beta}z^{-2}e^{-2\pi i\lambda};q)+r_{qb}(n|6),\label{eq:299}\end{align}
 and\begin{align}
|r_{qb}(n|6)|\le & 12\theta_{3}(|z|^{-2}q^{2\beta};q)\left\{ \frac{q^{j_{n}^{2}+2(\beta-1)j_{n}}}{|z|^{2j_{n}}}+|z|^{2j_{n}}q^{j_{n}^{2}-2(\beta+1)j_{n}}+\frac{q^{k_{n}}}{1-q}+\frac{\log n}{n^{\rho}}\right\} \label{eq:300}\end{align}
 for $n$ sufficiently large.
\item Assume that $\tau<0$. If for some fixed real numbers $\beta_{1}$,
$\beta_{2}$ and $\rho>0$ there exist infinitely many positive integers
$n$ such that\begin{equation}
-n\tau=m+\beta+a_{n},\quad m\in\mathbb{N},\quad|a_{n}|<\frac{1}{n^{\rho}},\label{eq:301}\end{equation}
 and\begin{equation}
n\theta=m_{1}+\beta+b_{n},\quad m_{1}\in\mathbb{Z},\quad|b_{n}|<\frac{1}{n^{\rho}}.\label{eq:302}\end{equation}
 Then, \begin{align}
\frac{J_{\nu}^{(2)}(2q^{ns-\nu/2}z;q)(q;q)_{\infty}^{2}e^{-n\pi i\theta(2m+\nu)}}{(-1)^{m}z^{\nu+2m}q^{m(2n\tau+m)+n\nu\tau-\nu^{2}/2}} & =\theta_{4}(q^{2\beta_{1}}z^{-2}e^{-2\pi i\beta_{2}};q)+r_{qb}(n|7),\label{eq:303}\end{align}
 and\begin{align}
|r_{qb}(n|7)| & \le156\theta_{3}(q^{2\beta_{1}}|z|^{-2};q)\left\{ \frac{q^{j_{n}^{2}+2(\beta_{1}-1)j_{n}}}{|z|^{2j_{n}}}+|z|^{2j_{n}}q^{j_{n}^{2}-2(\beta_{1}+1)j_{n}}+\frac{q^{k_{n}}}{1-q}+\frac{\log n}{n^{\rho}}\right\} \label{eq:304}\end{align}
 for $n$ sufficiently large.
\end{enumerate}
\end{thm*}
In the following corollaries we let \begin{equation}
z:=e^{\pi u},\quad u\in\mathbb{R}.\label{eq:305}\end{equation}

\begin{cor}
\label{cor:jackson-bessel-power}Assume that\begin{equation}
q=\exp(-n^{-a}\pi),\quad0<a<\frac{1}{2},\quad n\in\mathbb{N},\label{eq:306}\end{equation}
 we have the following results for $J_{\nu}^{(2)}(z;q)$:
\begin{enumerate}
\item Assume that $\tau>0$, we have\begin{align}
 & J_{\nu}^{(2)}(2\exp\pi(u-\tau n^{1-a}+\nu n^{-a}/2);\exp(-n^{-a}\pi))\label{eq:307}\\
 & =\frac{n^{a\nu}\exp(\pi\nu u-\pi\tau n^{1-a}+\nu^{2}n^{-a}\pi/2)}{(2\pi)^{\nu}}\left\{ \frac{1}{\Gamma(\nu+1)}+o(1)\right\} .\nonumber \end{align}
 for $n$ sufficiently large.
\item Assume that $\tau<0$ and for some fixed real number $\lambda$ there
are infinitely many positive integers $n$ such that \begin{equation}
-n\tau=m+\lambda,\quad m\in\mathbb{N}.\label{eq:308}\end{equation}
Then\begin{align}
 & J_{\nu}^{(2)}(2\exp\pi(u-\tau n^{1-a}+\nu n^{-a}/2);\exp(-n^{-a}\pi))\label{eq:309}\\
 & =\exp\left(\frac{\pi}{n^{a}}(n^{a}u+\nu/2-\tau n)^{2}+\frac{\pi n^{a}}{12}+\frac{\pi\nu^{2}n^{-a}}{4}-\frac{\pi n^{-a}}{12}\right)\nonumber \\
 & \times(-1)^{\lambda+n\tau}n^{-a}\left\{ \cos\pi(n^{a}u+\lambda)+\mathcal{O}(e^{-2\pi n^{a}})\right\} \nonumber \end{align}
for $n$ sufficiently large. 
\end{enumerate}
\end{cor}
Similarly, we have the following:

\begin{cor}
\label{cor:jackson-bessel-log}Assume that\begin{equation}
q=\exp(-\frac{\pi}{\gamma\log n}),\quad n\ge2,\quad\gamma>0,\label{eq:310}\end{equation}
we have the following results for $J_{\nu}^{(2)}(z;q)$:
\begin{enumerate}
\item Assume that $\tau>0$, we have\begin{align}
 & J_{\nu}^{(2)}\left(2\exp\pi\left(u-\frac{n\tau}{\gamma\log n}+\frac{\nu}{2\gamma\log n}\right);\exp\left(-\frac{\pi}{\gamma\log n}\right)\right)\label{eq:311}\\
 & =\left(\frac{\pi}{\gamma\log n}\right)^{\nu}\exp\left(\pi\nu u+\frac{\pi\nu^{2}}{2\gamma\log n}-\frac{\pi n\nu\tau}{\gamma\log n}\right)\left\{ \frac{1}{\Gamma(\nu+1)}+o(1)\right\} \nonumber \end{align}
 for $n$ sufficiently large.
\item Assume that $\tau<0$ and for some fixed real number $\lambda$ there
are infinitely many positive integers $n$ such that \begin{equation}
-n\tau=m+\lambda,\quad m\in\mathbb{N},\label{eq:312}\end{equation}
then,\begin{align}
 & J_{\nu}^{(2)}\left(2\exp\pi\left(u-\frac{n\tau}{\gamma\log n}+\frac{\nu}{2\gamma\log n}\right);\exp\left(-\frac{\pi}{\gamma\log n}\right)\right)\label{eq:313}\\
 & =\exp\left\{ \frac{\pi(u\gamma\log n-n\tau+\nu/2)^{2}}{\gamma\log n}-\frac{\pi}{12\gamma\log n}+\frac{\nu^{2}\pi}{4\gamma\log n}\right\} \nonumber \\
 & \times\frac{(-1)^{n\tau+\lambda}n^{\pi\gamma/12}}{\sqrt{\gamma\log n}}\left\{ \cos\pi(\gamma u\log n+\lambda)+\mathcal{O}(n^{-2\pi\gamma})\right\} ,\nonumber \end{align}
for $n$ sufficiently large.
\item Assume that $\tau<0$, if for some fixed real numbers $\beta$, $\rho\ge1$
and $\lambda$ there exist infinitely many positive integers $n$
such that\begin{equation}
-n\tau=m+\beta+a_{n},\quad|a_{n}|<\frac{1}{n^{\rho}},\quad m\in\mathbb{N},\label{eq:314}\end{equation}
then,\begin{align}
 & J_{\nu}^{(2)}\left(2\exp\pi\left(u-\frac{n\tau}{\gamma\log n}+\frac{\nu}{2\gamma\log n}\right);\exp\left(-\frac{\pi}{\gamma\log n}\right)\right)\label{eq:315}\\
 & =\exp\left\{ \frac{\pi(u\gamma\log n-n\tau+\nu/2)^{2}}{\gamma\log n}-\frac{\pi}{12\gamma\log n}+\frac{\nu^{2}\pi}{4\gamma\log n}\right\} \nonumber \\
 & \times\frac{(-1)^{m}n^{\rho/27}}{\sqrt{\gamma\log n}}\left\{ \cos\pi(\gamma u\log n+\lambda)+\mathcal{O}(n^{-8\rho/9})\right\} ,\nonumber \end{align}
for $n$ sufficiently large, where\begin{equation}
\gamma=\frac{4\rho}{9\pi}.\label{eq:316}\end{equation}

\end{enumerate}
\end{cor}

\subsection{Proofs for Theorem \ref{thm:q-bessel}}

Assume that $\tau<0$ and \begin{equation}
-n\tau=m+c_{n},\quad m\in\mathbb{N},\quad n\theta=m_{1}+d_{n},\quad m_{1}\in\mathbb{Z},\label{eq:317}\end{equation}
 then,\begin{align}
\frac{J_{\nu}^{(2)}(2q^{ns-\nu/2}z;q)(q;q)_{\infty}}{(q^{\nu+1};q)_{\infty}\left(q^{n\tau-\nu/2}ze^{n\pi i\theta}\right)^{\nu}} & =\sum_{k=0}^{\infty}\frac{q^{k^{2}}(-q^{2n\tau}z^{2}e^{2n\pi i\theta})^{k}}{(q,q^{\nu+1};q)_{k}}\label{eq:318}\\
 & =\sum_{k=0}^{m}\frac{q^{k^{2}}(-q^{2n\tau}z^{2}e^{2n\pi i\theta})^{k}}{(q,q^{\nu+1};q)_{k}}+\sum_{k=m+1}^{\infty}\frac{q^{k^{2}}(-q^{2n\tau}z^{2}e^{2n\pi i\theta})^{k}}{(q,q^{\nu+1};q)_{k}}\nonumber \\
 & =s_{1}+s_{2},\nonumber \end{align}
 and reverse summation order in $s_{1}$,\begin{align}
\frac{(q,q^{\nu+1};q)_{\infty}s_{1}}{(-z^{2}e^{2\pi in\theta})^{m}q^{m(2n\tau+m)}} & =\sum_{k=0}^{m}q^{k^{2}}(-q^{2c_{n}}z^{-2}e^{-2\pi id_{n}})^{k}e(k,n)\label{eq:319}\end{align}
 with\begin{equation}
e(k,n)=(q^{m-k+1},q^{\nu+1+m-k};q)_{\infty},\label{eq:320}\end{equation}
 then\begin{equation}
|e(k,n)|\le1\label{eq:321}\end{equation}
 for $0\le k\le m$. We also have \begin{equation}
|e(k,n)-1|\le\frac{6q^{k_{n}}}{1-q},\label{eq:322}\end{equation}
 for $0\le k\le k_{n}-1$ and $n$ sufficiently large. This could
be seen by expanding \begin{equation}
e(k,n)-1=\left\{ r_{1}\left(q;m-k\right)+1\right\} \left\{ r_{1}\left(q^{\nu+1};m-k\right)+1\right\} -1\label{eq:323}\end{equation}
and estimating each terms by Lemma \eqref{lem:1} to obtain

We shift the summation index from $k$ to $k+m$ in $s_{2}$,\begin{align}
\frac{(q,q^{\nu+1};q)_{\infty}s_{2}}{(-z^{2}e^{2\pi in\theta})^{m}q^{m(2n\tau+m)}} & =\sum_{k=1}^{\infty}q^{k^{2}}(-z^{2}q^{-2c_{n}}e^{2\pi id_{n}})^{k}f(k,n),\label{eq:324}\end{align}
and\begin{equation}
f(k,n)=(q^{m+k+1},q^{\nu+1+m+k};q)_{\infty}.\label{eq:325}\end{equation}
Thus,\begin{equation}
|f(k,n)|\le1\label{eq:326}\end{equation}
for $k\in\mathbb{N}$. Expand\begin{equation}
f(k,n)-1=\left\{ r_{1}\left(q;m+k\right)+1\right\} \left\{ r_{1}\left(q^{\nu+1};m+k\right)+1\right\} -1,\label{eq:327}\end{equation}
and estimate each term by Lemma \eqref{lem:1} to obtain\begin{equation}
|f(k,n)-1|\le\frac{6q^{k_{n}}}{1-q}\label{eq:328}\end{equation}
 for $k\in\mathbb{N}$ and $n$ sufficiently large. The rest of proof
is similar to the corresponding proof for $E_{q}(z)$.

\subsection{Proofs for Corollary \ref{cor:jackson-bessel-log}}

The $q$-Gamma function is defined as \cite{Andrews4,Gasper,Ismail2,Koekoek}\begin{equation}
\Gamma_{q}(x)=\frac{(q;q)_{\infty}}{(q^{x};q)_{\infty}}(1-q)^{1-x}\quad x\in\mathbb{C}.\label{eq:329}\end{equation}
It is a $q$-analogue of $\Gamma(x)$, \begin{equation}
\lim_{q\to1}\Gamma_{q}(x)=\Gamma(x).\label{eq:330}\end{equation}
 Thus\begin{equation}
\frac{(q^{\nu+1};q)_{\infty}(1-q)^{\nu}}{(q;q)_{\infty}}=\frac{1}{\Gamma(\nu+1)}+o(1)\label{eq:331}\end{equation}
 as $n$ sufficiently large. 

We also have\begin{align}
\theta_{4}(z^{-2}q^{2\lambda};q) & =\theta_{4}(ui+\lambda n^{-a}i|n^{-a}i)\label{eq:332}\\
 & =\sqrt{n^{a}}\exp\frac{\pi}{n^{a}}(n^{a}u+\lambda)^{2}\theta_{2}\left(n^{a}u+\lambda\mid n^{a}i\right)\nonumber \\
 & =2\sqrt{n^{a}}\exp\left(\frac{\pi}{n^{a}}(n^{a}u+\lambda)^{2}-\frac{\pi n^{a}}{4}\right)\nonumber \\
 & \times\cos\pi(n^{a}u+\lambda)\left\{ 1+\mathcal{O}(e^{-2\pi n^{a}})\right\} ,\nonumber \end{align}
 and\begin{align}
\theta_{3}\left(|z|^{-2}q^{2\lambda};q\right) & =\theta_{3}(ui+\lambda n^{-a}i|n^{-a}i)\label{eq:333}\\
 & =\sqrt{n^{a}}\exp\frac{\pi}{n^{a}}(n^{a}u+\lambda)^{2}\theta_{3}\left(n^{a}u+\lambda\mid n^{a}i\right)\nonumber \\
 & =\sqrt{n^{a}}\exp\frac{\pi}{n^{a}}(n^{a}u+\lambda)^{2}\left\{ 1+\mathcal{O}(e^{-\pi n^{a}})\right\} \nonumber \end{align}
 as $n\to\infty$. The rest of the proof is similar to corresponding
proof for $E_{q}(z)$ case.

\subsection{Proof for Corollary \ref{cor:jackson-bessel-log}}

Notice that \begin{equation}
\frac{1}{(q;q)_{\infty}}=\frac{n^{\pi\gamma/6}\left\{ 1+\mathcal{O}\left(n^{-4\pi\gamma}\right)\right\} }{\sqrt{2\gamma\log n}\exp\left(\frac{\pi}{24\gamma\log n}\right)},\label{eq:334}\end{equation}
and\begin{align}
\theta_{4}(z^{-2}q^{2\lambda};q) & =\theta_{4}\left(ui+\frac{i\lambda}{\gamma\log n}\mid\frac{i}{\gamma\log n}\right)\label{eq:335}\\
 & =\sqrt{\gamma\log n}\exp\frac{\pi(u\gamma\log n+\lambda)^{2}}{\gamma\log n}\theta_{2}\left(\gamma u\log n+\lambda\mid i\gamma\log n\right)\nonumber \\
 & =2n^{-\pi\gamma/4}\sqrt{\gamma\log n}\exp\left(\frac{\pi(u\gamma\log n+\lambda)^{2}}{\gamma\log n}\right)\nonumber \\
 & \times\cos\pi(\gamma u\log n+\lambda)\left\{ 1+\mathcal{O}(n^{-2\pi\gamma})\right\} ,\nonumber \end{align}
 and\begin{align}
\theta_{3}(z^{-2}q^{2\lambda};q) & =\theta_{3}\left(ui+\frac{i\lambda}{\gamma\log n}\mid\frac{i}{\gamma\log n}\right)\label{eq:336}\\
 & =\sqrt{\gamma\log n}\exp\frac{\pi(u\gamma\log n+\lambda)^{2}}{\gamma\log n}\theta_{3}\left(\gamma u\log n+\lambda\mid i\gamma\log n\right)\nonumber \\
 & =\sqrt{\gamma\log n}\exp\frac{\pi(u\gamma\log n+\lambda)^{2}}{\gamma\log n}\left\{ 1+\mathcal{O}(n^{-\pi\gamma})\right\} \nonumber \end{align}
as $n\to\infty$, and the rest of the proof is similar to the corresponding
proof for $E_{q}(z)$.

\section{Ismail-Masson Orthogonal Polynomials }

Ismail-Masson polynomials $\left\{ h_{n}(x|q)\right\} _{n=0}^{\infty}$
are defined as \cite{Ismail2} \begin{equation}
h_{n}(\sinh\xi|q)=\sum_{k=0}^{n}\frac{(q;q)_{n}q^{k(k-n)}(-1)^{k}e^{(n-2k)\xi}}{(q;q)_{k}(q;q)_{n-k}}.\label{eq:337}\end{equation}
Ismail-Masson polynomials come from an indeterminate moment problem.
Ismail-Masson orthogonal polynomials satisfy the following orthogonality

\begin{eqnarray}
\int_{-\infty}^{\infty}h_{m}(x|q)h_{n}(x|q)w_{im}(x)dx & = & q^{-n(n+1)/2}(q;q)_{n}\delta_{m,n},\label{eq:338}\end{eqnarray}
 where\begin{align}
w_{im}(x):= & q^{1/8}\sqrt{\frac{-2}{\pi\log q}}\exp\left(\frac{2\log^{2}(x+\sqrt{x^{2}+1})}{\log q}\right).\label{eq:339}\end{align}
The corresponding orthonormal Ismail-Masson functions are given by\begin{equation}
h_{n}(x\Vert q):=q^{n(n+1)/4}\sqrt{\frac{w_{im}(x)}{(q;q)_{n}}}h_{n}(x|q).\label{eq:340}\end{equation}

Let 

\begin{equation}
s=\frac{1+2\tau}{2}+i\frac{\theta\pi}{\log q},\quad\sinh\xi_{n}:=\frac{q^{-ns}z-q^{ns}z^{-1}}{2},\quad\tau,\theta\in\mathbb{R},\label{eq:341}\end{equation}
for any nonzero complex number $z$, then\emph{\begin{equation}
\frac{h_{n}(\sinh\xi_{n}|q)}{z^{n}q^{-n^{2}s}}=\sum_{k=0}^{n}\frac{q^{k^{2}}(q;q)_{n}e^{2nk\theta\pi i}}{(q;q)_{k}(q;q)_{n-k}}\left(-\frac{q^{2\tau n}}{z^{2}}\right)^{k}.\label{eq:342}\end{equation}
} Obviously,\emph{\begin{eqnarray}
\left|h_{n}(\sinh\xi_{n}|q)\right| & \le & \frac{|z|^{n}}{q^{n^{2}(\tau+1/2)}}\sum_{k=0}^{\infty}\frac{q^{k^{2}}}{(q;q)_{k}}\left(\frac{q^{2\tau n}}{|z|^{2}}\right)^{k}\le\frac{|z|^{n}A_{q}\left(-\frac{q^{2\tau n}}{|z|^{2}}\right)}{q^{n^{2}(\tau+1/2)}},\label{eq:343}\end{eqnarray}
 or\begin{equation}
\left|h_{n}(\sinh\xi_{n}|q)\right|\le\frac{(-\sqrt{q};q)_{\infty}|z|^{(4\tau+1)n}}{q^{2n^{2}(\tau+1/2)^{2}}}\exp\left(-\frac{2\log^{2}|z|}{\log q}\right).\label{eq:344}\end{equation}
 }

\subsection{Asymptotic Formulas For Ismail-Masson Polynomials\label{sec:ismail-masson asymptotics}}

\begin{thm}
\label{thm:ismail-masson}Given any nonzero complex number $z$, let
$s$ and $\xi_{n}$ be defined as in \eqref{eq:340} and \begin{equation}
j_{n}=\left\lfloor \frac{q^{4}\log n}{\log q^{-1}}\right\rfloor ,\quad k_{n}=\min\left\{ \left\lfloor \frac{(\tau+1)n}{2}\right\rfloor ,\left\lfloor \frac{-\tau n}{2}\right\rfloor \right\} ,\label{eq:345}\end{equation}
 we have the following results for Ismail-Masson polynomials:
\begin{enumerate}
\item Assume that $\tau>0$, we have \begin{equation}
\frac{h_{n}(\sinh\xi_{n}|q)}{z^{n}q^{-n^{2}s}}=1+r_{im}(n|1),\label{eq:346}\end{equation}
 and\begin{equation}
\left|r_{im}(n|1)\right|\le\frac{\exp\left(q^{2\tau n+1}/(\left|z\right|^{2}(1-q))\right)}{\left|z\right|^{2}(1-q)}q^{2\tau n+1}.\label{eq:347}\end{equation}

\item Assume that $\tau=0$. If for any fixed real number $\lambda$ there
are infinitely many positive integers $n$ such that \begin{equation}
n\theta=m+\lambda,\quad m\in\mathbb{Z},\label{eq:348}\end{equation}
then, \begin{equation}
\frac{h_{n}(\sinh\xi_{n}|q)}{z^{n}q^{-n^{2}s}}=A_{q}\left(\frac{e^{2\lambda\pi i}}{z^{2}}\right)+r_{im}(n|2),\label{eq:349}\end{equation}
and \begin{equation}
|r_{im}(n|2)|\le6\frac{A_{q}(-|z|^{-2})}{(q;q)_{\infty}}\left\{ q^{n/2}+\frac{q^{\left\lfloor n/2\right\rfloor ^{2}}}{|z|^{2\left\lfloor n/2\right\rfloor }}\right\} .\label{eq:350}\end{equation}
 for $n$ is sufficiently large.
\item Assume that $\tau=0$. If for fixed real numbers $\beta$ and $\rho\ge1$
there are infinitely positive integers $n$ such that\begin{equation}
n\theta=m+\beta+b_{n}\quad|b_{n}|<\frac{1}{n^{\rho}},\quad m\in\mathbb{Z},\label{eq:351}\end{equation}
then, \begin{equation}
\frac{h_{n}(\sinh\xi_{n}|q)}{z^{n}q^{-n^{2}s}}=A_{q}\left(\frac{e^{2\pi i\beta}}{z^{2}}\right)+e_{im}(n|3),\label{eq:352}\end{equation}
and\begin{equation}
|e_{im}(n|3)|\le24\exp\left(\frac{q}{\left|z\right|^{2}(1-q)}\right)\left\{ \frac{\log n}{n^{\rho}}+\frac{q^{n/2}}{1-q}+\frac{(|z|^{2}(q^{-1}-1))^{-j_{n}}}{j_{n}!}\right\} ,\label{eq:353}\end{equation}
for $n$ sufficiently large.
\item Assume that $-\frac{1}{2}<\tau<0$. If for some fixed real numbers
$\lambda$ and $\lambda_{1}$ there are infinite number of positive
integers $n$ such that\begin{equation}
-n\tau=m+\lambda,\quad m\in\mathbb{N},\quad n\theta=m_{1}+\lambda_{1},\quad m_{1}\in\mathbb{Z},\label{eq:354}\end{equation}
 then,\begin{align}
h_{n}(\sinh\xi_{n}|q) & =\frac{z^{n}q^{-n^{2}s+m\left(2\tau n+m\right)}}{(q;q)_{\infty}\left(-z^{2}e^{-2\pi i\lambda_{1}}\right)^{m}}\left\{ \theta_{4}\left(z^{2}q^{2\lambda}e^{-2\pi i\lambda_{1}};q\right)+r_{im}(n|4)\right\} ,\label{eq:355}\end{align}
 and\begin{align}
|r_{im}(n|4)| & \le28\theta_{3}\left(\left|z\right|^{2}q^{2\lambda};q\right)\left\{ \frac{q^{k_{n}}}{1-q}+|z|^{2k_{n}}q^{k_{n}^{2}+2\lambda k_{n}}+|z|^{-2k_{n}}q^{k_{n}^{2}-2\lambda k_{n}}\right\} .\label{eq:356}\end{align}
 for $n$ sufficiently large.
\item Assume that $-\frac{1}{2}<\tau<0$. If for some fixed real numbers
$\beta$, $\lambda$ and $\rho\ge1$ there are infinitely many positive
integers $n$ such that \begin{equation}
n\theta=m_{1}+\beta+b_{n},\quad|b_{n}|<\frac{1}{n^{\rho}},\quad m_{1}\in\mathbb{Z},\quad-n\tau=m+\lambda,\quad m\in\mathbb{N},\label{eq:357}\end{equation}
 then,\begin{align}
h_{n}(\sinh\xi_{n}|q) & =\frac{z^{n}q^{-n^{2}s+m\left(2\tau n+m\right)}}{(-z^{2}e^{-2n\theta\pi i})^{m}(q;q)_{\infty}}\left\{ \theta_{4}\left(z^{2}q^{2\lambda}e^{-2\beta\pi i};q\right)+r_{im}(n|5)\right\} ,\label{eq:358}\end{align}
 and\begin{align}
|r_{im}(n|5)| & \le24\theta_{3}(|z|^{2}q^{2\lambda};q)\left\{ |z|^{-2j_{n}}q^{j_{n}^{2}-2\lambda j_{n}}+|z|^{2j_{n}}q^{j_{n}^{2}+2\lambda j_{n}}+\frac{q^{k_{n}}}{1-q}+\frac{\log n}{n^{\rho}}\right\} \label{eq:359}\end{align}
for $n$ sufficiently large.
\item Assume that $-\frac{1}{2}<\tau<0$. If for fixed real numbers $\beta$,
$\lambda$ and $\rho\ge1$ there are infinitely many positive integers
$n$ such that \begin{equation}
-n\tau=m+\beta+a_{n},\quad|a_{n}|<\frac{1}{n^{\rho}},\quad m\in\mathbb{N},\quad n\theta=m_{1}+\lambda,\quad m_{1}\in\mathbb{Z},\label{eq:360}\end{equation}
 then,\begin{align}
h_{n}(\sinh\xi_{n}|q) & =\frac{z^{n}q^{-n^{2}s+m\left(2\tau n+m\right)}}{(-z^{2}e^{-2n\theta\pi i})^{m}(q;q)_{\infty}}\left\{ \theta_{4}\left(z^{2}q^{2\beta}e^{-2\lambda\pi i};q\right)+r_{im}(n|6)\right\} ,\label{eq:361}\end{align}
 and\begin{equation}
|r_{im}(n|6)|\le12\theta_{3}(|z|^{2}q^{2\beta};q)\left\{ \frac{\log n}{n^{\rho}}+\frac{q^{k_{n}}}{1-q}+|z|^{2j_{n}}q^{j_{n}^{2}+2(\beta-1)j_{n}}+|z|^{-2j_{n}}q^{j_{n}^{2}-2(\beta+1)j_{n}}\right\} \label{eq:362}\end{equation}
 for $n$ sufficiently large.
\item Assume that $-\frac{1}{2}<\tau<0$. If for fixed real numbers $\beta_{1},\,\beta_{2}$
and $\rho>0$ there are infinitely many positive integers $n$ such
that\begin{equation}
-\tau n=m+\beta_{1}+a_{n},\quad|a_{n}|<\frac{1}{n^{\rho}},\quad m\in\mathbb{N},\label{eq:363}\end{equation}
 and\begin{equation}
n\theta=m_{1}+\beta_{2}+b_{n},\quad|b_{n}|<\frac{1}{n^{\rho}},\quad m_{1}\in\mathbb{Z}.\label{eq:364}\end{equation}
 Then,\begin{align}
h_{n}(\sinh\xi_{n}|q) & =\frac{z^{n}q^{-n^{2}s+m\left(2\tau n+m\right)}}{(-z^{2}e^{-2n\theta\pi i})^{m}(q;q)_{\infty}}\left\{ \theta_{4}\left(z^{2}q^{2\beta_{1}}e^{-2\beta_{2}\pi i};q\right)+r_{im}(n|7)\right\} ,\label{eq:365}\end{align}
and\begin{align}
|r_{im}(n|7)| & \le54\theta_{3}(|z|^{2}q^{2\beta_{1}};q)\left\{ \frac{\log n}{n^{\rho}}+\frac{q^{k_{n}}}{1-q}+|z|^{2j_{n}}q^{j_{n}^{2}+2(\beta_{1}-1)j_{n}}+|z|^{-2j_{n}}q^{j_{n}^{2}-2(\beta_{1}+1)j_{n}}\right\} \label{eq:366}\end{align}
 for $n$ sufficiently large.
\end{enumerate}
\end{thm}
In the following corollaries we assume that\begin{equation}
z=e^{\pi u},\quad u\in\mathbb{R}.\label{eq:367}\end{equation}
 If we let \begin{equation}
q=\exp(-n^{-a}\pi),\quad0<a<\frac{1}{2},\quad n\in\mathbb{N},\label{eq:368}\end{equation}
then we have the following results for Ismail-Masson polynomials:

\begin{cor}
\label{cor:q-ismail-masson-power}
\begin{enumerate}
\item If$\tau>0$, then,\begin{align}
 & h_{n}\left(\sinh\pi\left(u+(\tau+1/2)n^{1-a}\right)|\exp(-\pi n^{-a})\right)\label{eq:369}\\
 & =\exp(n\pi u+(\tau+1/2)\pi n^{2-a})\left\{ 1+\mathcal{O}(e^{-4\pi n^{a}})\right\} ,\nonumber \end{align}
and\begin{align}
 & h_{n}\left(\sinh\pi(u+(\tau+1/2)n^{1-a})\Vert\exp(-\pi n^{-a})\right)\label{eq:370}\\
 & =\frac{1}{\sqrt{\pi}}\frac{\exp\left(-n^{-a}\pi(n^{a}u+\tau n)^{2}\right)}{\exp\frac{\pi}{12}\left(3n^{1-a}+n^{-a}-n^{a}\right)}\left\{ 1+\mathcal{O}(e^{-4\pi n^{a}})\right\} \nonumber \end{align}
 as $n\to\infty$.
\item Assume that  $-\frac{1}{2}<\tau<0$ and for some fixed real number
$\lambda$ there are infinitely many positive integers $n$ such that
\begin{equation}
-n\tau=m+\lambda,\quad m\in\mathbb{N},\label{eq:371}\end{equation}
then,\begin{align}
\mbox{} & h_{n}\left(\sinh\pi\left(u+(\tau+1/2)n^{1-a}\right)|\exp(-\pi n^{-a})\right)\label{eq:372}\\
 & =\frac{(-1)^{\tau n+\lambda}\sqrt{2}\left\{ \cos\pi(n^{a}u-\lambda)+\mathcal{O}(e^{-2\pi n^{a}})\right\} }{\exp\left\{ -\pi n^{-a}(n^{a}u+(\tau+1/2)n)^{2}-\frac{\pi n^{2-a}}{4}+\frac{\pi n^{a}}{12}+\frac{\pi n^{-a}}{24}\right\} },\nonumber \end{align}
 and\begin{align}
 & h_{n}\left(\sinh\pi(u+(\tau+1/2)n^{1-a})\Vert\exp(-\pi n^{-a})\right)\label{eq:373}\\
 & =\sqrt{\frac{2}{\pi}}\frac{(-1)^{\tau n+\lambda}\left\{ \cos\pi(n^{a}u-\lambda)+\mathcal{O}(e^{-2\pi n^{a}})\right\} }{\exp\left\{ \frac{n^{1-a}\pi}{4}+\frac{\pi n^{-a}}{8}\right\} }\nonumber \end{align}
as $n\to\infty$.
\end{enumerate}
\end{cor}
If we choose that \begin{equation}
q=\exp\left(-\frac{\pi}{\gamma\log n}\right),\quad n\ge2,\label{eq:374}\end{equation}
then we have the following:

\begin{cor}
\label{cor:q-ismail-masson-log}
\begin{enumerate}
\item If $\tau>0$, then we have\begin{align}
 & h_{n}\left(\sinh\pi\left(u+\frac{(\tau+1/2)n}{\gamma\log n}\right)\vert\exp\left(-\frac{\pi}{\gamma\log n}\right)\right)\label{eq:375}\\
 & =\exp\pi\left(nu+\frac{(\tau+1/2)n^{2}}{\gamma\log n}\right)\left\{ 1+\mathcal{O}\left(\exp(-2\pi\tau n/(\gamma\log n)\right)\right\} ,\nonumber \end{align}
and\begin{align}
 & h_{n}\left(\sinh\pi\left(u+\frac{(\tau+1/2)n}{\gamma\log n}\right)\Vert\exp\left(-\frac{\pi}{\gamma\log n}\right)\right)\label{eq:376}\\
 & ==\frac{n^{\gamma\pi/12}}{\sqrt{\pi}}\frac{\exp\left(-\frac{\pi}{\gamma\log n}(u\gamma\log n+n\tau)^{2}\right)}{\exp\left(\frac{n\pi}{4\gamma\log n}+\frac{\pi}{12\gamma\log n}\right)}\left\{ 1+\mathcal{O}(n^{-4\pi\gamma})\right\} \nonumber \end{align}
 for $n$ sufficiently large.
\item Assume that $-\frac{1}{2}<\tau<0$ and for some fixed real number
$\lambda$ there are infinitely many $m$ and $n>0$ such that \begin{equation}
-n\tau=m+\lambda,\quad m\in\mathbb{N}.\label{eq:377}\end{equation}
Then,\begin{align}
 & h_{n}\left(\sinh\pi\left(u+\frac{(\tau+1/2)n}{\gamma\log n}\right)\vert\exp\left(-\frac{\pi}{\gamma\log n}\right)\right)\label{eq:378}\\
 & =\frac{\sqrt{2}\exp\left(\frac{\pi}{\gamma\log n}(u\gamma\log n+(\tau+1/2)n)^{2}\right)}{(-1)^{\lambda+\tau n}n^{\pi\gamma/12}\exp(\frac{\pi}{24\gamma\log n}-\frac{n^{2}\pi}{4\gamma\log n})}\nonumber \\
 & \times\left\{ \cos\pi(u\gamma\log n-\lambda)+\mathcal{O}(n^{-2\pi\gamma})\right\} ,\nonumber \end{align}
and\begin{align}
 & h_{n}\left(\sinh\pi\left(u+\frac{(\tau+1/2)n}{\gamma\log n}\right)\Vert\exp\left(-\frac{\pi}{\gamma\log n}\right)\right)\label{eq:379}\\
 & =\sqrt{\frac{2}{\pi}}\frac{\left\{ \cos\pi(u\gamma\log n-\lambda)+\mathcal{O}(n^{-2\pi\gamma})\right\} }{(-1)^{\lambda+\tau n}\exp(\frac{\pi}{8\gamma\log n}+\frac{n\pi}{4\gamma\log n})}\nonumber \end{align}
for $n$ sufficiently large.
\item Assume that $-\frac{1}{2}<\tau<0$, if for some fixed real numbers
$\beta$, $\rho\ge1$ and $\lambda$ there exist infinitely many positive
integers $n$ such that\begin{equation}
-n\tau=m+\beta+a_{n},\quad|a_{n}|<\frac{1}{n^{\rho}},\quad m\in\mathbb{N},\label{eq:380}\end{equation}
and let \begin{equation}
\gamma=\frac{4\rho}{9\pi},\label{eq:381}\end{equation}
then,\begin{align}
 & h_{n}\left(\sinh\pi\left(u+\frac{(\tau+1/2)n}{\gamma\log n}\right)\vert\exp\left(-\frac{\pi}{\gamma\log n}\right)\right)\label{eq:382}\\
 & =\frac{\sqrt{2}\exp\left(\frac{\pi}{\gamma\log n}(u\gamma\log n+(\tau+1/2)n)^{2}\right)}{(-1)^{m}n^{\pi\gamma/12}\exp(\frac{\pi}{24\gamma\log n}-\frac{n^{2}\pi}{4\gamma\log n})}\nonumber \\
 & \times\left\{ \cos\pi(u\gamma\log n-\beta)+\mathcal{O}(n^{-8\rho/9}\log^{2}n)\right\} ,\nonumber \end{align}
and\begin{align}
 & h_{n}\left(\sinh\pi\left(u+\frac{(\tau+1/2)n}{\gamma\log n}\right)\Vert\exp\left(-\frac{\pi}{\gamma\log n}\right)\right)\label{eq:383}\\
 & =\sqrt{\frac{2}{\pi}}\frac{\left\{ \cos\pi(u\gamma\log n-\beta)+\mathcal{O}(n^{-8\rho/9}\log^{2}n)\right\} }{(-1)^{m}\exp(\frac{\pi}{8\gamma\log n}+\frac{n\pi}{4\gamma\log n})}\nonumber \end{align}
for $n$ sufficiently large. 
\end{enumerate}
\end{cor}

\subsection{Proof For Theorem \ref{thm:ismail-masson}}

In first three proofs we have used the following inequalities\begin{equation}
0<\frac{(q;q)_{n}}{(q;q)_{n-k}}\le1\label{eq:384}\end{equation}
 for $0\le k\le n$ and \begin{equation}
\left|\frac{(q;q)_{n}}{(q;q)_{n-k}}-1\right|\le\frac{6q^{1+n/2}}{1-q}\label{eq:385}\end{equation}
 for $0\le k\le\left\lfloor \frac{n}{2}\right\rfloor -1$ and $n$
sufficiently large. This can be seen by applying Lemma \ref{lem:1}
\begin{equation}
\frac{(q;q)_{n}}{(q;q)_{n-k}}-1=\left\{ r_{1}(q;n-k)+1\right\} \left\{ r_{2}(q;n)+1\right\} -1.\label{eq:386}\end{equation}
 Assume that\begin{equation}
-\tau n=m+c_{n},\quad m\in\mathbb{N},\quad n\theta=m_{1}+d_{n},\quad m_{1}\in\mathbb{Z},\label{eq:387}\end{equation}
then,\begin{align}
 & \frac{h_{n}(\sinh_{n}|q)}{z^{n}q^{-n^{2}s}}=\sum_{k=0}^{n}\frac{q^{k^{2}}(q;q)_{n}e^{2nk\theta\pi i}}{(q;q)_{k}(q;q)_{n-k}}\left(-\frac{q^{2\tau n}}{z^{2}}\right)^{k}\label{eq:388}\\
 & =\sum_{k=0}^{m}\frac{q^{k^{2}}(q;q)_{n}e^{2nk\theta\pi i}}{(q;q)_{k}(q;q)_{n-k}}\left(-\frac{q^{2\tau n}}{z^{2}}\right)^{k}+\sum_{k=m+1}^{n}\frac{q^{k^{2}}(q;q)_{n}e^{2nk\theta\pi i}}{(q;q)_{k}(q;q)_{n-k}}\left(-\frac{q^{2\tau n}}{z^{2}}\right)^{k}\nonumber \\
 & =s_{1}+s_{2}.\nonumber \end{align}
 We reverse the summation order in $s_{1}$ to obtain\begin{equation}
\frac{s_{1}(q;q)_{\infty}(-z^{2}e^{-2n\theta\pi i})^{m}}{q^{m(2\tau n+m)}}=\sum_{k=0}^{m}q^{k^{2}}\left(-z^{2}q^{2c_{n}}e^{-2\pi id_{n}}\right)^{k}e(k,n),\label{eq:389}\end{equation}
 and\begin{equation}
e(k,n)=\frac{(q;q)_{\infty}(q;q)_{n}}{(q;q)_{m-k}(q;q)_{n-m+k}}.\label{eq:390}\end{equation}
It is clear that\begin{equation}
|e(k,n)|\le1\label{eq:391}\end{equation}
 for $0\le k\le m$. Expand \begin{align}
e(k,n) & -1=\left\{ r_{2}\left(q;n\right)+1\right\} \left\{ r_{1}\left(q;m-k\right)+1\right\} \left\{ r_{1}\left(q;n-m+k\right)+1\right\} -1\label{eq:392}\end{align}
and estimate each term by Lemma \ref{lem:1} to get\begin{equation}
|e(k,n)-1|\le\frac{14q^{2+k_{n}}}{1-q}\label{eq:393}\end{equation}
 for $0\le k\le k_{n}-1$ and $n$ sufficiently large.

In sum $s_{2}$ we shift summation index from $k$ to $k+m$ to get
\begin{equation}
\frac{s_{2}(q;q)_{\infty}(-z^{2}e^{-2n\theta\pi i})^{m}}{q^{m(2\tau n+m)}}=\sum_{k=1}^{n-m}q^{k^{2}}\left(-z^{-2}q^{-2c_{n}}e^{2\pi id_{n}}\right)^{k}f(k,n),\label{eq:394}\end{equation}
 and\begin{equation}
f(k,n)=\frac{(q;q)_{\infty}(q;q)_{n}}{(q;q)_{m+k}(q;q)_{n-m-k}},\quad|f(k,n)|\le1,\label{eq:395}\end{equation}
 for $1\le k\le n-m$. Apply Lemma \ref{lem:1} to each term of \begin{align}
f(k,n)-1 & =\left\{ r_{2}\left(q;n\right)+1\right\} \left\{ r_{1}\left(q;m+k\right)+1\right\} \left\{ r_{1}\left(q;n-m-k\right)+1\right\} -1,\label{eq:396}\end{align}
to obtain \begin{equation}
|f(k,n)-1|\le\frac{14q^{2+k_{n}}}{1-q}\label{eq:397}\end{equation}
 for $1\le k\le k_{n}-1$ for $n$ sufficiently large.

\subsubsection{Proof for case 1. }

If we write\begin{equation}
\frac{h_{n}(\sinh\xi_{n}|q)}{z^{n}q^{-n^{2}s}}=1+r_{im}(n|1),\label{eq:398}\end{equation}
 then\begin{equation}
r_{im}(n|1)=\sum_{k=1}^{n}\frac{q^{k^{2}}(q;q)_{n}e^{2nk\theta\pi i}}{(q;q)_{k}(q;q)_{n-k}}\left(-\frac{q^{2\tau n}}{z^{2}}\right)^{k}.\label{eq:399}\end{equation}
 Thus\begin{align}
|r_{im}(n|1)| & \le\sum_{k=1}^{\infty}\frac{q^{k^{2}}}{(q;q)_{k}}\left(\frac{q^{2\tau n}}{\left|z\right|^{2}}\right)^{k}\le\sum_{k=1}^{\infty}\frac{1}{k!}\left(\frac{q^{2\tau n+1}}{\left|z\right|^{2}(1-q)}\right)^{k}\le\frac{q^{2\tau n+1}}{\left|z\right|^{2}(1-q)}\exp\left(\frac{q^{2\tau n+1}}{\left|z\right|^{2}(1-q)}\right).\label{eq:400}\end{align}

\subsubsection{Proof for case 2.}

We have\begin{align}
 & \frac{h_{n}(\sinh\xi_{n}|q)}{z^{n}q^{-n^{2}s}}=\sum_{k=0}^{n}\frac{q^{k^{2}}}{(q;q)_{k}}\left(-\frac{e^{2\lambda\pi i}}{z^{2}}\right)^{k}\frac{(q;q)_{n}}{(q;q)_{n-k}}\label{eq:401}\\
 & =\sum_{k=0}^{\infty}\frac{q^{k^{2}}}{(q;q)_{k}}\left(-\frac{e^{2\lambda\pi i}}{z^{2}}\right)^{k}-\sum_{k=\left\lfloor n/2\right\rfloor }^{\infty}\frac{q^{k^{2}}}{(q;q)_{k}}\left(-\frac{e^{2\lambda\pi i}}{z^{2}}\right)^{k}\nonumber \\
 & +\sum_{k=0}^{\left\lfloor n/2\right\rfloor -1}\frac{q^{k^{2}}}{(q;q)_{k}}\left(-\frac{e^{2\lambda\pi i}}{z^{2}}\right)^{k}\left\{ \frac{(q;q)_{n}}{(q;q)_{n-k}}-1\right\} +\sum_{k=\left\lfloor n/2\right\rfloor }^{n}\frac{q^{k^{2}}}{(q;q)_{k}}\left(-\frac{e^{2\lambda\pi i}}{z^{2}}\right)^{k}\frac{(q;q)_{n}}{(q;q)_{n-k}}\nonumber \\
 & =A_{q}\left(\frac{e^{2\lambda\pi i}}{z^{2}}\right)+s_{1}+s_{2}+s_{3}.\nonumber \end{align}
 Then,\begin{align}
|s_{1}+s_{3}| & \le2\sum_{k=\left\lfloor n/2\right\rfloor }^{\infty}\frac{q^{k^{2}}}{(q;q)_{k}}\left(\frac{1}{\left|z\right|^{2}}\right)^{k}\le\frac{2q^{\left\lfloor n/2\right\rfloor ^{2}}A_{q}(-|z|^{-2})}{(q;q)_{\infty}|z|^{2\left\lfloor n/2\right\rfloor }},\label{eq:402}\end{align}
 and \begin{equation}
|s_{2}|\le\frac{6q^{n/2}}{1-q}A_{q}(-|z|^{-2}).\label{eq:403}\end{equation}
 for $n$ sufficiently large. 

Let\begin{equation}
r_{im}(n|2)=s_{1}+s_{2}+s_{3},\label{eq:404}\end{equation}
then, \begin{equation}
\frac{h_{n}(\sinh\xi_{n}|q)}{z^{n}q^{-n^{2}s}}=A_{q}\left(\frac{e^{2\lambda\pi i}}{z^{2}}\right)+r_{im}(n|2),\label{eq:405}\end{equation}
 and\begin{equation}
|r_{im}(n|2)|\le6\frac{A_{q}(-|z|^{-2})}{(q;q)_{\infty}}\left\{ q^{n/2}+\frac{q^{\left\lfloor n/2\right\rfloor ^{2}}}{|z|^{2\left\lfloor n/2\right\rfloor }}\right\} .\label{eq:406}\end{equation}

\subsubsection{Proof for case 3.}

In this case we have\begin{align}
 & \frac{h_{n}(\sinh\xi_{n}|q)}{z^{n}q^{-n^{2}s}}=\sum_{k=0}^{n}\frac{q^{k^{2}}}{(q;q)_{k}}\left(-\frac{e^{2\pi i\beta}}{z^{2}}\right)^{k}\frac{(q;q)_{n}e^{2\pi ikb_{n}}}{(q;q)_{n-k}}\label{eq:407}\\
 & =\sum_{k=0}^{\infty}\frac{q^{k^{2}}}{(q;q)_{k}}\left(-\frac{e^{2\pi i\beta}}{z^{2}}\right)^{k}-\sum_{k=j_{n}}^{\infty}\frac{q^{k^{2}}}{(q;q)_{k}}\left(-\frac{e^{2\pi i\beta}}{z^{2}}\right)^{k}\nonumber \\
 & +\sum_{k=0}^{j_{n}-1}\frac{q^{k^{2}}}{(q;q)_{k}}\left(-\frac{e^{2\pi i\beta}}{z^{2}}\right)^{k}\left\{ \frac{(q;q)_{n}}{(q;q)_{n-k}}-1\right\} +\sum_{k=0}^{j_{n}-1}\frac{q^{k^{2}}}{(q;q)_{k}}\left(-\frac{e^{2\pi i\beta}}{z^{2}}\right)^{k}\frac{(q;q)_{n}}{(q;q)_{n-k}}\left\{ e^{2\pi ikb_{n}}-1\right\} \nonumber \\
 & +\sum_{k=j_{n}}^{n}\frac{q^{k^{2}}}{(q;q)_{k}}\left(-\frac{e^{2\pi i\beta}}{z^{2}}\right)^{k}\frac{(q;q)_{n}}{(q;q)_{n-k}}e^{2\pi ikb_{n}}=A_{q}\left(\frac{e^{2\pi i\beta}}{z^{2}}\right)+s_{1}+s_{2}+s_{3}+s_{4}.\nonumber \end{align}
Then,\begin{align}
|s_{1}+s_{4}| & \le2\sum_{k=j_{n}}^{\infty}\frac{q^{k^{2}}|z|^{-2k}}{(q;q)_{k}}\le\frac{2}{j_{n}!}\left(\frac{q}{\left|z\right|^{2}(1-q)}\right)^{j_{n}}\exp\left(\frac{q}{\left|z\right|^{2}(1-q)}\right).\label{eq:408}\end{align}
Given a positive number $\rho$, since that \begin{equation}
1<j_{n}<\frac{n}{2},\quad j_{n}\le\frac{n^{\rho}}{2\pi}\label{eq:409}\end{equation}
 for $n$ sufficiently large, clearly, \begin{equation}
|s_{2}|\le\frac{6q^{n/2}}{1-q}\exp\left(\frac{q}{\left|z\right|^{2}(1-q)}\right),\label{eq:410}\end{equation}
and\begin{equation}
|s_{3}|\le\frac{2\pi j_{n}}{n^{\rho}}e^{2\pi j_{n}/n^{\rho}}\sum_{k=0}^{\infty}\frac{q^{k^{2}}|z|^{-2k}}{(q;q)_{k}}\le\frac{24\log n}{n^{\rho}}\exp\left(\frac{q}{\left|z\right|^{2}(1-q)}\right).\label{eq:411}\end{equation}
 Let \begin{equation}
e_{im}(n|3)=s_{1}+s_{2}+s_{3}+s_{4},\label{eq:412}\end{equation}
 then\begin{equation}
\frac{h_{n}(\sinh\xi_{n}|q)}{z^{n}q^{-n^{2}s}}=A_{q}\left(\frac{e^{2\pi i\beta}}{z^{2}}\right)+e_{im}(n|3),\label{eq:413}\end{equation}
 and\begin{equation}
|e_{im}(n|3)|\le24\exp\left(\frac{q}{\left|z\right|^{2}(1-q)}\right)\left\{ \frac{\log n}{n^{\rho}}+\frac{q^{n/2}}{1-q}+\frac{(|z|^{2}(q^{-1}-1))^{-j_{n}}}{j_{n}!}\right\} \label{eq:414}\end{equation}
 for $n$ sufficiently large.

\subsubsection{Proof for case 4.}

Observe that\begin{align}
 & \frac{s_{1}\left(-z^{2}e^{-2\pi i\lambda_{1}}\right)^{m}(q;q)_{\infty}}{q^{m\left(2\tau n+m\right)}}=\sum_{k=0}^{\infty}q^{k^{2}}\left(-z^{2}q^{2\lambda}e^{-2\pi i\lambda_{1}}\right)^{k}\label{eq:415}\\
 & -\sum_{k=k_{n}}^{\infty}q^{k^{2}}\left(-z^{2}q^{2\lambda}e^{-2\pi i\lambda_{1}}\right)^{k}+\sum_{k=0}^{k_{n}-1}q^{k^{2}}\left(-z^{2}q^{2\lambda}e^{-2\pi i\lambda_{1}}\right)^{k}\left(e(k,n)-1\right)\nonumber \\
 & +\sum_{k=k_{n}}^{m}q^{k^{2}}\left(-z^{2}q^{2\lambda}e^{-2\pi i\lambda_{1}}\right)^{k}e(k,n)=\sum_{k=0}^{\infty}q^{k^{2}}\left(-z^{2}q^{2\lambda}e^{-2\pi i\lambda_{1}}\right)^{k}+s_{11}+s_{12}+s_{13}.\nonumber \end{align}
Then,\begin{align}
|s_{11}+s_{13}| & \le2\sum_{k=k_{n}}^{\infty}q^{k^{2}}\left(\left|z\right|^{2}q^{2\lambda}\right)^{k}\le2|z|^{2k_{n}}q^{k_{n}^{2}+2\lambda k_{n}}\sum_{k=0}^{\infty}q^{k^{2}}\left(\left|z\right|^{2}q^{2\lambda}\right)^{k}\label{eq:416}\\
 & \le2\theta_{3}\left(\left|z\right|^{2}q^{2\lambda}\mid q\right)|z|^{2k_{n}}q^{k_{n}^{2}+2\lambda k_{n}},\nonumber \end{align}
 and\begin{align}
|s_{12}| & \le\frac{14q^{k_{n}}}{1-q}\sum_{k=0}^{\infty}q^{k^{2}}(|z|^{2}q^{2\lambda})^{k}\le14\theta_{3}\left(\left|z\right|^{2}q^{2\lambda};q\right)\frac{q^{k_{n}}}{1-q}\label{eq:417}\end{align}
for $n$ sufficiently large.

Let \begin{equation}
r_{1}(n)=s_{11}+s_{12}+s_{13},\label{eq:418}\end{equation}
 then,\begin{equation}
\frac{s_{1}\left(-z^{2}e^{-2\pi i\lambda_{1}}\right)^{m}(q;q)_{\infty}}{q^{m\left(2\tau n+m\right)}}=\sum_{k=0}^{\infty}q^{k^{2}}\left(-z^{2}q^{2\lambda}e^{-2\pi i\lambda_{1}}\right)^{k}+r_{1}(n)\label{eq:419}\end{equation}
\begin{equation}
|r_{1}(n)|\le14\theta_{3}\left(\left|z\right|^{2}q^{2\lambda};q\right)\left\{ |z|^{2k_{n}}q^{k_{n}^{2}+2\lambda k_{n}}+\frac{q^{k_{n}}}{1-q}\right\} \label{eq:420}\end{equation}
 for $n$ sufficiently large.

Similarly,\begin{align}
 & \frac{s_{2}\left(-z^{2}e^{-2\pi i\lambda_{1}}\right)^{m}(q;q)_{\infty}}{q^{m\left(2\tau n+m\right)}}=\sum_{k=1}^{\infty}q^{k^{2}}\left(-z^{-2}q^{-2\lambda}e^{2\pi i\lambda_{1}}\right)^{k}\label{eq:421}\\
 & -\sum_{k=k_{n}}^{\infty}q^{k^{2}}\left(-z^{-2}q^{-2\lambda}e^{2\pi i\lambda_{1}}\right)^{k}+\sum_{k=1}^{k_{n}-1}q^{k^{2}}\left(-z^{-2}q^{-2\lambda}e^{2\pi i\lambda_{1}}\right)^{k}\left(f(k,n)-1\right)\nonumber \\
 & +\sum_{k=k_{n}}^{n-m}q^{k^{2}}\left(-z^{-2}q^{-2\lambda}e^{2\pi i\lambda_{1}}\right)^{k}f(k,n)=\sum_{k=-\infty}^{-1}q^{k^{2}}\left(-z^{2}q^{2\lambda}e^{-2\pi i\lambda_{1}}\right)^{k}+s_{21}+s_{22}+s_{23}.\nonumber \end{align}
Then, \begin{align}
|s_{21}+s_{23}| & \le2\sum_{k=k_{n}}^{\infty}q^{k^{2}}\left(\left|z\right|^{-2}q^{-2\lambda}\right)^{k}\le2\theta_{3}(|z|^{2}q^{2\lambda};q)\frac{q^{k_{n}^{2}-2\lambda k_{n}}}{|z|^{2k_{n}}},\label{eq:422}\end{align}
 and\begin{align}
|s_{22}| & \le\frac{14q^{k_{n}}}{1-q}\sum_{k=0}^{\infty}q^{k^{2}-2\lambda k}|z|^{-2k}\le14\theta_{3}(|z|^{2}q^{2\lambda};q)\frac{q^{k_{n}}}{1-q}.\label{eq:423}\end{align}
 for $n$ sufficiently large.

Let\begin{equation}
r_{2}(n)=s_{21}+s_{22}+s_{23},\label{eq:424}\end{equation}
then, \begin{equation}
\frac{s_{2}\left(-z^{2}e^{-2\pi i\lambda_{1}}\right)^{m}(q;q)_{\infty}}{q^{m\left(2\tau n+m\right)}}=\sum_{k=-\infty}^{-1}q^{k^{2}}\left(-z^{2}q^{2\lambda}e^{-2\pi i\lambda_{1}}\right)^{k}+r_{2}(n),\label{eq:425}\end{equation}
 and\begin{equation}
|r_{2}(n)|\le14\theta_{3}\left(\left|z\right|^{2}q^{2\lambda};q\right)\left\{ \frac{q^{k_{n}}}{1-q}+\frac{q^{k_{n}^{2}-2\lambda k_{n}}}{|z|^{2k_{n}}}\right\} \label{eq:426}\end{equation}
 for $n$ sufficiently large.

Thus,\begin{equation}
\frac{\left(-z^{2}e^{-2\pi i\lambda_{1}}\right)^{m}h_{n}(\sinh\xi_{n}|q)(q;q)_{\infty}}{z^{n}q^{-n^{2}s+m\left(2\tau n+m\right)}}=\theta_{4}\left(z^{2}q^{2\lambda}e^{-2\pi i\lambda_{1}};q\right)+r_{im}(n|4)\label{eq:427}\end{equation}
 with\begin{align}
|r_{im}(n|4)| & \le28\theta_{3}\left(\left|z\right|^{2}q^{2\lambda};q\right)\left\{ \frac{q^{k_{n}}}{1-q}+|z|^{2k_{n}}q^{k_{n}^{2}+2\lambda k_{n}}+|z|^{-2k_{n}}q^{k_{n}^{2}-2\lambda k_{n}}\right\} \label{eq:428}\end{align}
 for $n$ sufficiently large.

\subsubsection{Proof for case 5.}

In this case we have\begin{align}
 & \frac{s_{1}(q;q)_{\infty}(-z^{2}e^{-2n\theta\pi i})^{m}}{q^{m(2\tau n+m)}}=\sum_{k=0}^{m}q^{k^{2}}\left(-z^{2}q^{2\lambda}e^{-2\pi i\beta}\right)^{k}e^{-2k\pi ib_{n}}e(k,n)\label{eq:429}\\
 & =\sum_{k=0}^{\infty}q^{k^{2}}\left(-z^{2}q^{2\lambda}e^{-2\pi i\beta}\right)^{k}-\sum_{k=j_{n}}^{\infty}q^{k^{2}}\left(-z^{2}q^{2\lambda}e^{-2\pi i\beta}\right)^{k}\nonumber \\
 & +\sum_{k=0}^{j_{n}-1}q^{k^{2}}\left(-z^{2}q^{2\lambda}e^{-2\pi i\beta}\right)^{k}\left\{ e^{-2k\pi ib_{n}}-1\right\} +\sum_{k=0}^{j_{n}-1}q^{k^{2}}\left(-z^{2}q^{2\lambda}e^{-2\pi i\beta}\right)^{k}e^{-2k\pi ib_{n}}\left\{ e(k,n)-1\right\} \nonumber \\
 & +\sum_{k=j_{n}}^{m}q^{k^{2}}\left(-z^{2}q^{2\lambda}e^{-2\pi i\beta}\right)^{k}e^{-2k\pi ib_{n}}e(k,n)=\sum_{k=0}^{\infty}q^{k^{2}}\left(-z^{2}q^{2\lambda}e^{-2\pi i\beta}\right)^{k}+s_{11}+s_{12}+s_{13}+s_{14}.\nonumber \end{align}
 Then,\begin{align}
|s_{11}+s_{14}| & \le2\sum_{k=j_{n}}^{\infty}q^{k^{2}+2\lambda k}|z|^{2k}\le2\theta_{3}(|z|^{2}q^{2\lambda}\mid q)|z|^{2j_{n}}q^{j_{n}^{2}+2\lambda j_{n}}.\label{eq:430}\end{align}
 For sufficiently large $n$,\begin{align}
|s_{12}| & \le\frac{24j_{n}}{n^{\rho}}\sum_{k=0}^{\infty}q^{k^{2}+2\lambda k}|z|^{2k}\le24\theta_{3}(|z|^{2}q^{2\lambda}\mid q)\frac{\log n}{n^{\rho}},\label{eq:431}\end{align}
and\begin{align}
|s_{13}| & \le\frac{14q^{k_{n}}}{1-q}\sum_{k=0}^{\infty}q^{k^{2}+2\lambda k}|z|^{2k}\le14\theta_{3}(|z|^{2}q^{2\lambda}\mid q)\frac{q^{k_{n}}}{1-q}.\label{eq:432}\end{align}
 Let\begin{equation}
r_{1}(n)=s_{11}+s_{12}+s_{13}+s_{14},\label{eq:433}\end{equation}
 then,\begin{equation}
\frac{s_{1}(q;q)_{\infty}(-z^{2}e^{-2n\theta\pi i})^{m}}{q^{m(2\tau n+m)}}=\sum_{k=0}^{\infty}q^{k^{2}}\left(-z^{2}q^{2\lambda}e^{-2\pi i\beta}\right)^{k}+r_{1}(n),\label{eq:434}\end{equation}
 and\begin{equation}
|r_{1}(n)|\le24\theta_{3}(|z|^{2}q^{2\lambda}\mid q)\left\{ |z|^{2j_{n}}q^{j_{n}^{2}+2\lambda j_{n}}+\frac{q^{k_{n}}}{1-q}+\frac{\log n}{n^{\rho}}\right\} \label{eq:435}\end{equation}
 for $n$ sufficiently large.

Similarly,\begin{align}
 & \frac{s_{2}(q;q)_{\infty}(-z^{2}e^{-2n\theta\pi i})^{m}}{q^{m(2\tau n+m)}}=\sum_{k=1}^{n-m}q^{k^{2}}\left(-z^{-2}q^{-2\lambda}e^{2\beta\pi i}\right)^{k}e^{2k\pi ib_{n}}f(k,n)\label{eq:436}\\
 & =\sum_{k=1}^{\infty}q^{k^{2}}\left(-z^{-2}q^{-2\lambda}e^{2\beta\pi i}\right)^{k}-\sum_{k=j_{n}}^{\infty}q^{k^{2}}\left(-z^{-2}q^{-2\lambda}e^{2\beta\pi i}\right)^{k}\nonumber \\
 & +\sum_{k=1}^{j_{n}-1}q^{k^{2}}\left(-z^{-2}q^{-2\lambda}e^{2\beta\pi i}\right)^{k}\left\{ e^{2k\pi ib_{n}}-1\right\} +\sum_{k=1}^{j_{n}-1}q^{k^{2}}\left(-z^{-2}q^{-2\lambda}e^{2\beta\pi i}\right)^{k}e^{2k\pi ib_{n}}\left\{ f(k,n)-1\right\} \nonumber \\
 & +\sum_{k=j_{n}}^{n-m}q^{k^{2}}\left(-z^{-2}q^{-2\lambda}e^{2\beta\pi i}\right)^{k}e^{2k\pi ib_{n}}f(k,n)=\sum_{k=-\infty}^{-1}q^{k^{2}}\left(-z^{2}q^{2\lambda}e^{-2\beta\pi i}\right)^{k}+s_{21}+s_{22}+s_{23}+s_{24},\nonumber \end{align}
 then\begin{align}
|s_{21}+s_{24}| & \le2\sum_{k=j_{n}}^{\infty}q^{k^{2}-2\lambda k}|z|^{-2k}\le2\theta_{3}(|z|^{2}q^{2\lambda};q)\frac{q^{j_{n}^{2}-2\lambda j_{n}}}{|z|^{2j_{n}}}.\label{eq:437}\end{align}
 For sufficiently large $n$ \begin{align}
|s_{22}| & \le\frac{24j_{n}}{n^{\rho}}\sum_{k=0}^{\infty}q^{k^{2}-2\lambda k}|z|^{-2k}\le24\theta_{3}(|z|^{2}q^{2\lambda};q)\frac{\log n}{n^{\rho}},\label{eq:438}\end{align}
 and\begin{align}
|s_{23}| & \le14\frac{q^{k_{n}}}{1-q}\sum_{k=0}^{\infty}q^{k^{2}-2\lambda k}|z|^{-2k}\le14\theta_{3}(|z|^{2}q^{2\lambda};q)\frac{q^{k_{n}}}{1-q}.\label{eq:439}\end{align}
 for $n$ sufficiently large. 

Let\begin{equation}
r_{2}(n)=s_{21}+s_{22}+s_{23}+s_{24},\label{eq:440}\end{equation}
then,\begin{equation}
\frac{s_{2}(q;q)_{\infty}(-z^{2}e^{-2n\theta\pi i})^{m}}{q^{m(2\tau n+m)}}=\sum_{k=-\infty}^{-1}q^{k^{2}}\left(-z^{2}q^{2\lambda}e^{-2\beta\pi i}\right)^{k}+r_{2}(n),\label{eq:441}\end{equation}
and\begin{equation}
|r_{2}(n)|\le24\theta_{3}(|z|^{2}q^{2\lambda};q)\left\{ \frac{q^{j_{n}^{2}-2\lambda j_{n}}}{|z|^{2j_{n}}}+\frac{q^{k_{n}}}{1-q}+\frac{\log n}{n^{\rho}}\right\} \label{eq:442}\end{equation}
 for $n$ sufficiently large. 

Thus, \begin{equation}
\frac{(-z^{2}e^{-2n\theta\pi i})^{m}h_{n}(\sinh\xi_{n}|q)(q;q)_{\infty}}{z^{n}q^{-n^{2}s+m\left(2\tau n+m\right)}}=\theta_{4}\left(z^{2}q^{2\lambda}e^{-2\beta\pi i};q\right)+r_{im}(n|5),\label{eq:443}\end{equation}
 and\begin{align}
|r_{im}(n|5)| & \le48\theta_{3}(|z|^{2}q^{2\lambda};q)\left\{ \frac{q^{j_{n}^{2}-2\lambda j_{n}}}{|z|^{2j_{n}}}+|z|^{2j_{n}}q^{j_{n}^{2}+2\lambda j_{n}}+\frac{q^{k_{n}}}{1-q}+\frac{\log n}{n^{\rho}}\right\} \label{eq:444}\end{align}
 for sufficiently large $n$.

\subsubsection{Proof for case 6.}

Observe that\begin{align}
 & \frac{s_{1}(q;q)_{\infty}(-z^{2}e^{-2n\theta\pi i})^{m}}{q^{m(2\tau n+m)}}=\sum_{k=0}^{m}q^{k^{2}}\left(-z^{2}q^{2\beta}e^{-2\pi i\lambda}\right)^{k}q^{2ka_{n}}e(k,n)\label{eq:445}\\
 & =\sum_{k=0}^{\infty}q^{k^{2}}\left(-z^{2}q^{2\beta}e^{-2\pi i\lambda}\right)^{k}-\sum_{k=j_{n}}^{\infty}q^{k^{2}}\left(-z^{2}q^{2\beta}e^{-2\pi i\lambda}\right)^{k}\nonumber \\
 & +\sum_{k=0}^{j_{n}-1}q^{k^{2}}\left(-z^{2}q^{2\beta}e^{-2\pi i\lambda}\right)^{k}\left\{ q^{2ka_{n}}-1\right\} +\sum_{k=0}^{j_{n}-1}q^{k^{2}}\left(-z^{2}q^{2\beta}e^{-2\pi i\lambda}\right)^{k}q^{2ka_{n}}\left\{ e(k,n)-1\right\} \nonumber \\
 & +\sum_{k=j_{n}}^{m}q^{k^{2}}\left(-z^{2}q^{2\beta}e^{-2\pi i\lambda}\right)^{k}q^{2ka_{n}}e(k,n)=\sum_{k=0}^{\infty}q^{k^{2}}\left(-z^{2}q^{2\beta}e^{-2\pi i\lambda}\right)^{k}+s_{11}+s_{12}+s_{13}+s_{14}.\nonumber \end{align}
Since\begin{equation}
|q^{2ka_{n}}|\le q^{-2k}\label{eq:446}\end{equation}
for $0\le k\le n$, then,\begin{equation}
|s_{11}+s_{14}|\le2\sum_{k=j_{n}}^{\infty}q^{k^{2}}(|z|^{2}q^{2\beta-2})^{k}\le2\theta_{3}(|z|^{2}q^{2\beta};q)|z|^{2j_{n}}q^{j_{n}^{2}+2(\beta-1)j_{n}}.\label{eq:447}\end{equation}
For $n$ sufficiently large we have\begin{equation}
|q^{2ka_{n}}-1|\le\frac{6\log n}{n^{\rho}},\quad|q^{2ka_{n}}|\le3\label{eq:448}\end{equation}
for $0\le k\le j_{n}-1$, then,\begin{equation}
|s_{12}|\le6\theta_{3}(|z|^{2}q^{2\beta};q)\frac{\log n}{n^{\rho}},\label{eq:449}\end{equation}
and\begin{equation}
|s_{13}|\le3\theta_{3}(|z|^{2}q^{2\beta};q)\frac{q^{k_{n}}}{1-q}.\label{eq:450}\end{equation}
Thus,\begin{equation}
\frac{s_{1}(q;q)_{\infty}(-z^{2}e^{-2n\theta\pi i})^{m}}{q^{m(2\tau n+m)}}=\sum_{k=0}^{\infty}q^{k^{2}}\left(-z^{2}q^{2\beta}e^{-2\pi i\lambda}\right)^{k}+r_{1}(n),\label{eq:451}\end{equation}
and\begin{equation}
|r_{1}(n)|\le6\theta_{3}(|z|^{2}q^{2\beta};q)\left\{ \frac{\log n}{n^{\rho}}+\frac{q^{k_{n}}}{1-q}+|z|^{2j_{n}}q^{j_{n}^{2}+2(\beta-1)j_{n}}\right\} \label{eq:452}\end{equation}
for $n$ sufficiently large. 

Similarly,\begin{align}
 & \frac{s_{2}(q;q)_{\infty}(-z^{2}e^{-2n\theta\pi i})^{m}}{q^{m(2\tau n+m)}}=\sum_{k=1}^{n-m}q^{k^{2}}\left(-z^{-2}q^{-2\beta}e^{2\lambda\pi i}\right)^{k}q^{-2ka_{n}}f(k,n)\label{eq:453}\\
 & =\sum_{k=-\infty}^{-1}q^{k^{2}}\left(-z^{2}q^{2\beta}e^{-2\pi i\lambda}\right)^{k}-\sum_{k=j_{n}}^{\infty}q^{k^{2}}\left(-z^{-2}q^{-2\beta}e^{2\lambda\pi i}\right)^{k}\nonumber \\
 & +\sum_{k=1}^{j_{n}-1}q^{k^{2}}\left(-z^{-2}q^{-2\beta}e^{2\lambda\pi i}\right)^{k}\left\{ q^{-2ka_{n}}-1\right\} +\sum_{k=1}^{j_{n}-1}q^{k^{2}}\left(-z^{-2}q^{-2\beta}e^{2\lambda\pi i}\right)^{k}q^{-2ka_{n}}\left\{ f(k,n)-1\right\} \nonumber \\
 & +\sum_{k=j_{n}}^{n-m}q^{k^{2}}\left(-z^{-2}q^{-2\beta}e^{2\lambda\pi i}\right)^{k}q^{-2ka_{n}}f(k,n)=\sum_{k=-\infty}^{-1}q^{k^{2}}\left(-z^{2}q^{2\beta}e^{-2\pi i\lambda}\right)^{k}+s_{21}+s_{22}+s_{23}+s_{24}.\nonumber \end{align}
Then,\begin{equation}
|s_{21}+s_{24}|\le2\sum_{k=j_{n}}^{\infty}q^{k^{2}}(|z|^{-2}q^{-2\beta-2})^{k}\le2\theta_{3}(|z|^{2}q^{2\beta};q)|z|^{-2j_{n}}q^{j_{n}^{2}-2(\beta+1)j_{n}}.\label{eq:454}\end{equation}
For $n$ sufficiently large,\begin{equation}
|s_{22}|\le6\theta_{3}(|z|^{2}q^{2\beta};q)\frac{\log n}{n^{\rho}},\label{eq:455}\end{equation}
and\begin{equation}
|s_{23}|\le3\theta_{3}(|z|^{2}q^{2\beta};q)\frac{q^{k_{n}}}{1-q}.\label{eq:456}\end{equation}
 Therefore,

\begin{equation}
\frac{s_{2}(q;q)_{\infty}(-z^{2}e^{-2n\theta\pi i})^{m}}{q^{m(2\tau n+m)}}=\sum_{k=-\infty}^{-1}q^{k^{2}}\left(-z^{2}q^{2\beta}e^{-2\pi i\lambda}\right)^{k}+r_{2}(n),\label{eq:457}\end{equation}
 and\begin{equation}
|r_{2}(n)|\le6\theta_{3}(|z|^{2}q^{2\beta};q)\left\{ \frac{\log n}{n^{\rho}}+\frac{q^{k_{n}}}{1-q}+|z|^{-2j_{n}}q^{j_{n}^{2}-2(\beta+1)j_{n}}\right\} \label{eq:458}\end{equation}
for $n$ sufficiently large.

Hence,\begin{equation}
\frac{(-z^{2}e^{-2n\theta\pi i})^{m}h_{n}(\sinh\xi_{n}|q)(q;q)_{\infty}}{z^{n}q^{-n^{2}s+m\left(2\tau n+m\right)}}=\theta_{4}(z^{2}q^{2\beta}e^{-2\pi i\lambda};q)+r_{im}(n|6),\label{eq:459}\end{equation}
 and\begin{equation}
|r_{im}(n|6)|\le12\theta_{3}(|z|^{2}q^{2\beta};q)\left\{ \frac{\log n}{n^{\rho}}+\frac{q^{k_{n}}}{1-q}+|z|^{2j_{n}}q^{j_{n}^{2}+2(\beta-1)j_{n}}+|z|^{-2j_{n}}q^{j_{n}^{2}-2(\beta+1)j_{n}}\right\} \label{eq:460}\end{equation}
 for $n$ sufficiently large.

\subsubsection{Proof for case 7. }

In this case we have\begin{align}
 & \frac{s_{1}(q;q)_{\infty}(-z^{2}e^{-2n\theta\pi i})^{m}}{q^{m(2\tau n+m)}}=\sum_{k=0}^{m}q^{k^{2}}\left(-z^{2}q^{2\beta_{1}}e^{-2\pi i\beta_{2}}\right)^{k}e^{-2k\pi ib_{n}}q^{2ka_{n}}e(k,n)\label{eq:461}\\
 & =\sum_{k=0}^{\infty}q^{k^{2}}\left(-z^{2}q^{2\beta_{1}}e^{-2\pi i\beta_{2}}\right)^{k}-\sum_{k=j_{n}}^{\infty}q^{k^{2}}\left(-z^{2}q^{2\beta_{1}}e^{-2\pi i\beta_{2}}\right)^{k}\nonumber \\
 & +\sum_{k=0}^{j_{n}-1}q^{k^{2}}\left(-z^{2}q^{2\beta_{1}}e^{-2\pi i\beta_{2}}\right)^{k}\left\{ q^{2ka_{n}}-1\right\} +\sum_{k=0}^{j_{n}-1}q^{k^{2}}\left(-z^{2}q^{2\beta_{1}}e^{-2\pi i\beta_{2}}\right)^{k}\left\{ e^{-2k\pi ib_{n}}-1\right\} q^{2ka_{n}}\nonumber \\
 & +\sum_{k=0}^{j_{n}-1}q^{k^{2}}\left(-z^{2}q^{2\beta_{1}}e^{-2\pi i\beta_{2}}\right)^{k}e^{-2k\pi ib_{n}}q^{2ka_{n}}\left\{ e(k,n)-1\right\} +\sum_{k=j_{n}}^{m}q^{k^{2}}\left(-z^{2}q^{2\beta_{1}}e^{-2\pi i\beta_{2}}\right)^{k}e^{-2k\pi ib_{n}}q^{2ka_{n}}e(k,n)\nonumber \\
 & =\sum_{k=0}^{\infty}q^{k^{2}}\left(-z^{2}q^{2\beta_{1}}e^{-2\pi i\beta_{2}}\right)^{k}+s_{11}+s_{12}+s_{13}+s_{14}+s_{15}.\nonumber \end{align}
Then,\begin{equation}
\frac{s_{1}(q;q)_{\infty}(-z^{2}e^{-2n\theta\pi i})^{m}}{q^{m(2\tau n+m)}}=\sum_{k=0}^{\infty}q^{k^{2}}\left(-z^{2}q^{2\beta_{1}}e^{-2\pi i\beta_{2}}\right)^{k}+r_{1}(n),\label{eq:462}\end{equation}
and\begin{equation}
|r_{1}(n)|\le27\theta_{3}(|z|^{2}q^{2\beta_{1}};q)\left\{ \frac{\log n}{n^{\rho}}+\frac{q^{k_{n}}}{1-q}+|z|^{2j_{n}}q^{j_{n}^{2}+2(\beta_{1}-1)j_{n}}\right\} \label{eq:463}\end{equation}
for $n$ sufficiently large.

Similarly,\begin{align}
 & \frac{s_{2}(q;q)_{\infty}(-z^{2}e^{-2n\theta\pi i})^{m}}{q^{m(2\tau n+m)}}=\sum_{k=1}^{n-m}q^{k^{2}}\left(-z^{-2}q^{-2\beta_{1}}e^{2\beta_{2}\pi i}\right)^{k}e^{2k\pi ib_{n}}q^{-2ka_{n}}f(k,n)\label{eq:464}\\
 & =\sum_{k=-\infty}^{-1}q^{k^{2}}\left(-z^{2}q^{2\beta_{1}}e^{-2\beta_{2}\pi i}\right)^{k}-\sum_{k=j_{n}}^{\infty}q^{k^{2}}\left(-z^{-2}q^{-2\beta_{1}}e^{2\beta_{2}\pi i}\right)^{k}\nonumber \\
 & +\sum_{k=1}^{j_{n}-1}q^{k^{2}}\left(-z^{-2}q^{-2\beta_{1}}e^{2\beta_{2}\pi i}\right)^{k}\left\{ e^{2k\pi ib_{n}}-1\right\} +\sum_{k=1}^{j_{n}-1}q^{k^{2}}\left(-z^{-2}q^{-2\beta_{1}}e^{2\beta_{2}\pi i}\right)^{k}e^{2k\pi ib_{n}}\left\{ q^{-2ka_{n}}-1\right\} \nonumber \\
 & +\sum_{k=1}^{j_{n}-1}q^{k^{2}}\left(-z^{-2}q^{-2\beta_{1}}e^{2\beta_{2}\pi i}\right)^{k}e^{2k\pi ib_{n}}q^{-2ka_{n}}\left\{ f(k,n)-1\right\} +\sum_{k=j_{n}}^{n-m}q^{k^{2}}\left(-z^{-2}q^{-2\beta_{1}}e^{2\beta_{2}\pi i}\right)^{k}e^{2k\pi ib_{n}}q^{-2ka_{n}}f(k,n)\nonumber \\
 & =\sum_{k=-\infty}^{-1}q^{k^{2}}\left(-z^{2}q^{2\beta_{1}}e^{-2\beta_{2}\pi i}\right)^{k}+s_{21}+s_{22}+s_{23}+s_{24}+s_{25}.\nonumber \end{align}
Thus,\begin{equation}
\frac{s_{2}(q;q)_{\infty}(-z^{2}e^{-2n\theta\pi i})^{m}}{q^{m(2\tau n+m)}}=\sum_{k=-\infty}^{-1}q^{k^{2}}\left(-z^{2}q^{2\beta_{1}}e^{-2\beta_{2}\pi i}\right)^{k}+r_{2}(n),\label{eq:465}\end{equation}
and\begin{equation}
|r_{2}(n)|\le27\theta_{3}(|z|^{2}q^{2\beta_{1}};q)\left\{ \frac{\log n}{n^{\rho}}+\frac{q^{k_{n}}}{1-q}+|z|^{-2j_{n}}q^{j_{n}^{2}-2(\beta_{1}+1)j_{n}}\right\} \label{eq:466}\end{equation}
 for $n$ sufficiently large.

Therefore,

\begin{equation}
\frac{(-z^{2}e^{-2n\theta\pi i})^{m}h_{n}(\sinh\xi_{n}|q)(q;q)_{\infty}}{z^{n}q^{-n^{2}s+m\left(2\tau n+m\right)}}=\theta_{4}(z^{2}q^{2\beta_{1}}e^{-2\pi i\beta_{2}};q)+r_{im}(n|7),\label{eq:467}\end{equation}
and\begin{equation}
|r_{im}(n|7)|\le54\theta_{3}(|z|^{2}q^{2\beta_{1}};q)\left\{ \frac{\log n}{n^{\rho}}+\frac{q^{k_{n}}}{1-q}+|z|^{2j_{n}}q^{j_{n}^{2}+2(\beta_{1}-1)j_{n}}+|z|^{-2j_{n}}q^{j_{n}^{2}-2(\beta_{1}+1)j_{n}}\right\} \label{eq:468}\end{equation}
 For $n$ sufficiently large.

\subsection{Proof for Corollary \ref{cor:q-ismail-masson-power}}

In this case we have,

\begin{equation}
x_{n}(u)=\sinh\xi_{n}=\sinh\pi\left(u+(\tau+1/2)n^{1-a}\right),\label{eq:469}\end{equation}
and\begin{equation}
w_{im}(x_{n}(u))=\frac{\sqrt{2n^{a}}}{\pi}\exp\left\{ -\frac{n^{-a}\pi}{8}-2n^{a}\pi\left(u+(\tau+1/2)n^{1-a}\right)^{2}\right\} .\label{eq:470}\end{equation}
 From Lemma \ref{lem:1} and Lemma \ref{lem:2} to get \begin{equation}
\frac{1}{(q;q)_{n}}=\frac{r_{1}(q;n)+1}{(q;q)_{\infty}}=\frac{\exp(\pi n^{a}/6-\pi n^{-a}/24)}{\sqrt{2n^{a}}}\left\{ 1+\mathcal{O}(e^{-4\pi n^{a}})\right\} \label{eq:471}\end{equation}
as $n\to\infty$, thus,\begin{eqnarray}
q^{n(n+1)/4}\sqrt{\frac{w_{im}(x_{n}(u))}{(q;q)_{n}}} & = & \frac{\exp\left\{ -n^{-a}\pi(n^{a}u+(\tau+1/2)n)^{2}+\frac{\pi n^{a}}{12}-\frac{\pi n^{-a}}{12}-\frac{(n+1)n\pi}{4}\right\} }{\sqrt{\pi}}\label{eq:472}\\
 & \times & \left\{ 1+\mathcal{O}(e^{-4\pi n^{a}})\right\} \nonumber \end{eqnarray}
as $n\to\infty$.

\subsubsection{Proof of case 1.}

Observe that\begin{equation}
z^{n}q^{-n^{2}(\tau+1/2)}=\exp(n\pi u+(\tau+1/2)\pi n^{2-a}),\label{eq:473}\end{equation}
 then formulas \eqref{eq:344} and \eqref{eq:345} imply\begin{align}
 & h_{n}\left(\sinh\pi\left(u+(\tau+1/2)n^{1-a}\right)|\exp(-\pi n^{-a})\right)\label{eq:474}\\
 & =\exp(n\pi u+(\tau+1/2)\pi n^{2-a})\left\{ 1+\mathcal{O}(e^{-4\pi n^{a}})\right\} ,\nonumber \end{align}
and\begin{align}
 & h_{n}\left(\sinh\pi(u+(\tau+1/2)n^{1-a})\Vert\exp(-\pi n^{-a})\right)\label{eq:475}\\
 & =\frac{1}{\sqrt{\pi}}\frac{\exp\left(-n^{-a}\pi(n^{a}u+\tau n)^{2}\right)}{\exp\frac{\pi}{12}\left(3n^{1-a}+n^{-a}-n^{a}\right)}\left\{ 1+\mathcal{O}(e^{-4\pi n^{a}})\right\} \nonumber \end{align}
 as $n\to\infty$.

\subsubsection{Proof for case 2.}

From \eqref{eq:19}-- \eqref{eq:32} to get \begin{eqnarray}
\theta_{3}(|z|^{2}q^{2\lambda};q) & = & \theta_{3}(ui-\lambda n^{-a}i\mid n^{-a}i)\label{eq:476}\\
 & = & \sqrt{n^{a}}\exp\pi n^{-a}(n^{a}u-\lambda)^{2}\left\{ 1+\mathcal{O}(e^{-\pi n^{a}})\right\} ,\nonumber \end{eqnarray}
 and\begin{eqnarray}
\theta_{4}(z^{2}q^{2\lambda};q) & = & \theta_{4}(ui-\lambda n^{-a}i\mid n^{-a}i)\label{eq:477}\\
 & = & \sqrt{n^{a}}\exp\pi n^{-a}(n^{a}u-\lambda)^{2}\theta_{2}\left(n^{a}u-\lambda\vert n^{a}i\right)\nonumber \\
 & =2 & \sqrt{n^{a}}\exp\left\{ \pi n^{-a}(n^{a}u-\lambda)^{2}-\frac{\pi n^{a}}{4}\right\} \nonumber \\
 & \times & \cos\pi(n^{a}u-\lambda)\left\{ 1+\mathcal{O}(e^{-2\pi n^{a}})\right\} \nonumber \end{eqnarray}
 as $n\to\infty$.

From formulas \eqref{eq:353}--\eqref{eq:356}, we have\begin{align}
\mbox{} & h_{n}\left(\sinh\pi\left(u+(\tau+1/2)n^{1-a}\right)|\exp(-\pi n^{-a})\right)\label{eq:478}\\
 & =\frac{(-1)^{\tau n+\lambda}\sqrt{2}\left\{ \cos\pi(n^{a}u-\lambda)+\mathcal{O}(e^{-2\pi n^{a}})\right\} }{\exp\left\{ -\pi n^{-a}(n^{a}u+(\tau+1/2)n)^{2}-\frac{\pi n^{2-a}}{4}+\frac{\pi n^{a}}{12}+\frac{\pi n^{-a}}{24}\right\} }\nonumber \end{align}
 and\begin{align}
 & h_{n}\left(\sinh\pi(u+(\tau+1/2)n^{1-a})\Vert\exp(-\pi n^{-a})\right)\label{eq:479}\\
 & =\sqrt{\frac{2}{\pi}}\frac{(-1)^{\tau n+\lambda}\left\{ \cos\pi(n^{a}u-\lambda)+\mathcal{O}(e^{-2\pi n^{a}})\right\} }{\exp\left\{ \frac{n^{1-a}\pi}{4}+\frac{\pi n^{-a}}{8}\right\} }\nonumber \end{align}
as $n\to\infty$.

\subsection{Proof for Corollary \ref{cor:q-ismail-masson-log} }

In this case we have\begin{equation}
\sinh\xi_{n}=\sinh\pi\left(u+\frac{(\tau+1/2)n}{\gamma\log n}\right),\label{eq:480}\end{equation}
\begin{align}
w_{im}(\sinh\xi_{n})= & \frac{\sqrt{2\gamma\log n}}{\pi}\exp\left(-\frac{\pi}{8\gamma\log n}-\frac{2\pi}{\gamma\log n}\left(\gamma u\log n+(\tau+1/2)n\right)^{2}\right),\label{eq:481}\end{align}
and \begin{equation}
\frac{1}{(q;q)_{n}}=\frac{\exp\left(\frac{\gamma\pi\log n}{6}-\frac{\pi}{24\gamma\log n}\right)}{\sqrt{2\gamma\log n}}\left\{ 1+\mathcal{O}(n^{-4\pi\gamma})\right\} ,\label{eq:482}\end{equation}
\begin{align}
 & q^{n(n+1)/4}\sqrt{\frac{w_{im}(\sinh\xi_{n})}{(q;q)_{n}}}\label{eq:483}\\
 & =\frac{n^{\gamma\pi/12}\exp\left(-\frac{\pi}{\gamma\log n}\left(\gamma u\log n+(\tau+1/2)n\right)^{2}\right)}{\sqrt{\pi}\exp\left(\frac{\pi}{12\gamma\log n}+\frac{n(n+1)\pi}{4\gamma\log n}\right)}\left\{ 1+\mathcal{O}(n^{-4\pi\gamma})\right\} \nonumber \end{align}
 as $n\to\infty$.

\subsubsection{Proof for case 1. }

From formulas \eqref{eq:344} and \eqref{eq:345} to get\begin{align}
 & h_{n}\left(\sinh\pi\left(u+\frac{(\tau+1/2)n}{\gamma\log n}\right)\vert\exp\left(-\frac{\pi}{\gamma\log n}\right)\right)\label{eq:484}\\
 & =\exp\pi\left(nu+\frac{(\tau+1/2)n^{2}}{\gamma\log n}\right)\left\{ 1+\mathcal{O}\left(\exp(-2\pi\tau n/(\gamma\log n)\right)\right\} ,\nonumber \end{align}
and\begin{align}
 & h_{n}\left(\sinh\pi\left(u+\frac{(\tau+1/2)n}{\gamma\log n}\right)\Vert\exp\left(-\frac{\pi}{\gamma\log n}\right)\right)\label{eq:485}\\
 & ==\frac{n^{\gamma\pi/12}}{\sqrt{\pi}}\frac{\exp\left(-\frac{\pi(u\gamma\log n+n\tau)^{2}}{\gamma\log n}\right)}{\exp\left(\frac{n\pi}{4\gamma\log n}+\frac{\pi}{12\gamma\log n}\right)}\left\{ 1+\mathcal{O}(n^{-4\pi\gamma})\right\} \nonumber \end{align}
 as $n\to\infty$.

\subsubsection{Proof for case 2. }

In this case we have\begin{align}
\theta_{4}(z^{2}q^{2\lambda};q) & =\theta_{4}\left(ui-\frac{\lambda i}{\gamma\log n}\mid\frac{i}{\gamma\log n}\right)\label{eq:486}\\
 & =\frac{\sqrt{\gamma\log n}}{n^{\pi\gamma/4}}\exp\left(\frac{\pi(u\gamma\log n-\lambda)^{2}}{\gamma\log n}\right)\nonumber \\
 & \times\cos\pi(u\gamma\log n-\lambda)\left\{ 1+\mathcal{O}(n^{-2\pi\gamma})\right\} ,\nonumber \end{align}
 and\begin{align}
\theta_{3}(\left|z\right|^{2}q^{2\lambda};q) & =\theta_{3}\left(ui-\frac{\lambda i}{\gamma\log n}\mid\frac{i}{\gamma\log n}\right)\label{eq:487}\\
 & =\sqrt{\gamma\log n}\exp\left(\frac{\pi(u\gamma\log n-\lambda)^{2}}{\gamma\log n}\right)\left\{ 1+\mathcal{O}(n^{-\pi\gamma})\right\} .\nonumber \end{align}
Then,\begin{align}
 & h_{n}\left(\sinh\pi\left(u+\frac{(\tau+1/2)n}{\gamma\log n}\right)\vert\exp\left(-\frac{\pi}{\gamma\log n}\right)\right)\label{eq:488}\\
 & =\frac{\sqrt{2}\exp\left(\frac{\pi}{\gamma\log n}(u\gamma\log n+(\tau+1/2)n)^{2}\right)}{(-1)^{\lambda+\tau n}n^{\pi\gamma/12}\exp(\frac{\pi}{24\gamma\log n}-\frac{n^{2}\pi}{4\gamma\log n})}\nonumber \\
 & \times\left\{ \cos\pi(u\gamma\log n-\lambda)+\mathcal{O}(n^{-2\pi\gamma})\right\} ,\nonumber \end{align}
and\begin{align}
 & h_{n}\left(\sinh\pi\left(u+\frac{(\tau+1/2)n}{\gamma\log n}\right)\Vert\exp\left(-\frac{\pi}{\gamma\log n}\right)\right)\label{eq:489}\\
 & =\sqrt{\frac{2}{\pi}}\frac{\left\{ \cos\pi(u\gamma\log n-\lambda)+\mathcal{O}(n^{-2\pi\gamma})\right\} }{(-1)^{\lambda+\tau n}\exp(\frac{\pi}{8\gamma\log n}+\frac{n\pi}{4\gamma\log n})}\nonumber \end{align}
 for $n$ sufficiently large.

\subsubsection{Proof for case 3. }

Put $\theta,m_{1},\lambda=0$ in \eqref{eq:359}-\eqref{eq:361} to
get\begin{align}
 & h_{n}\left(\sinh\pi\left(u+\frac{(\tau+1/2)n}{\gamma\log n}\right)\vert\exp\left(-\frac{\pi}{\gamma\log n}\right)\right)\label{eq:490}\\
 & =\frac{\sqrt{2}\exp\left(\frac{\pi}{\gamma\log n}(u\gamma\log n+(\tau+1/2)n)^{2}\right)}{(-1)^{m}n^{\pi\gamma/12}\exp(\frac{\pi}{24\gamma\log n}-\frac{n^{2}\pi}{4\gamma\log n})}\nonumber \\
 & \times\left\{ \cos\pi(u\gamma\log n-\beta)+\mathcal{O}(n^{-8\rho/9}\log^{2}n)\right\} ,\nonumber \end{align}
and\begin{align}
 & h_{n}\left(\sinh\pi\left(u+\frac{(\tau+1/2)n}{\gamma\log n}\right)\Vert\exp\left(-\frac{\pi}{\gamma\log n}\right)\right)\label{eq:491}\\
 & =\sqrt{\frac{2}{\pi}}\frac{\left\{ \cos\pi(u\gamma\log n-\beta)+\mathcal{O}(n^{-8\rho/9}\log^{2}n)\right\} }{(-1)^{m}\exp(\frac{\pi}{8\gamma\log n}+\frac{n\pi}{4\gamma\log n})}\nonumber \end{align}
for $n$ sufficiently large, where $\gamma=\frac{8\rho}{9\pi}$.

\section{Stieltjes-Wigert Orthogonal Polynomials }

Stieltjes\emph{-}Wigert orthogonal polynomials $\left\{ S_{n}(x;q)\right\} _{n=0}^{\infty}$
are defined as \cite{Ismail2}\begin{equation}
S_{n}(x;q)=\sum_{k=0}^{n}\frac{q^{k^{2}}(-x)^{k}}{(q;q)_{k}(q;q)_{n-k}}.\label{eq:492}\end{equation}
Stieltjes-Wigert orthogonal polynomials come from an indeterminant
moment problem. They satisfy the orthogonality relation \begin{equation}
\int_{0}^{\infty}S_{m}(x;q)S_{n}(x;q)w_{sw}(x)dx=\frac{q^{-n}}{(q;q)_{n}}\delta_{m,n},\label{eq:493}\end{equation}
 where\begin{align}
w_{sw}(x) & :=\sqrt{\frac{-1}{2\pi\log q}}\exp\left(\frac{1}{2\log q}\left[\log\left(\frac{x}{\sqrt{q}}\right)\right]^{2}\right).\label{eq:494}\end{align}
 Clearly, the associated orthonormal Stieltjes-Wigert functions are
given by \begin{equation}
s_{n}(x\Vert q):=\sqrt{q^{n}(q;q)_{n}w_{sw}(x)}S_{n}(x;q).\label{eq:495}\end{equation}
 For any nonzero complex number $z$, let \begin{equation}
s=2\tau+2+i\frac{2\theta\pi}{\log q},\quad x_{n}(z,s)=zq^{-ns},\quad\tau,\theta\in\mathbb{R}.\label{eq:496}\end{equation}
 Reverse the summation order of \eqref{eq:492} to get\begin{equation}
\frac{S_{n}(x_{n}(z,s);q)}{(-z)^{n}q^{n^{2}(1-s)}}=\sum_{k=0}^{n}\frac{q^{k^{2}}e^{2nk\theta\pi i}}{(q;q)_{k}(q;q)_{n-k}}\left(-\frac{q^{2\tau n}}{z}\right)^{k}.\label{eq:497}\end{equation}
 Clearly,\begin{eqnarray}
|S_{n}(x_{n}(z,s);q)| & \le & \frac{|z|^{n}q^{-n^{2}(2\tau+1)}}{(q;q)_{\infty}}\sum_{k=0}^{\infty}\frac{q^{k^{2}}(q^{2\tau n}|z|^{-1})^{k}}{(q;q)_{k}}\le\frac{|z|^{n}A_{q}(-q^{2\tau n}|z|^{-1})}{(q;q)_{\infty}q^{n^{2}(2\tau+1)}},\label{eq:498}\end{eqnarray}
 or\begin{equation}
|S_{n}(x_{n}(z,s);q)|\le\frac{(-\sqrt{q};q)_{\infty}}{(q;q)_{\infty}}\frac{|z|^{(2\tau+1)n}\exp\left(-\frac{\log^{2}|z|}{2\log q}\right)}{q^{(2\tau^{2}+2\tau+1)n^{2}}}.\label{eq:499}\end{equation}

\subsection{Asymptotic Formulas For Stieltjes-Wigert Polynomials}

\begin{thm}
\label{thm:stieltjes-wigert}Given any nonzero complex number $z$,
let $s$ and $x_{n}(z,s)$ be defined as in \eqref{eq:493} and \eqref{eq:494}
and let\begin{equation}
j_{n}:=\left\lfloor \frac{q^{4}\log n}{\log q^{-1}}\right\rfloor ,\quad k_{n}=\min\left\{ \left\lfloor \frac{(1+\tau)n}{2}\right\rfloor ,\left\lfloor \frac{(-\tau)n}{2}\right\rfloor \right\} ,\label{eq:500}\end{equation}
 we have the following results for Stieltjes-Wigert polynomials: 
\begin{enumerate}
\item When $\tau>0$, we have \begin{equation}
\frac{S_{n}(x_{n}(z,s);q)(q;q)_{\infty}}{(-z)^{n}q^{n^{2}(1-s)}}=1+r_{sw}(n|1),\label{eq:501}\end{equation}
 and\begin{equation}
\left|r_{sw}(n|1)\right|\le\frac{q^{2\tau n+1}\exp(\left(q^{2\tau n+1}/((1-q)|z|)\right)}{(1-q)|z|}.\label{eq:502}\end{equation}

\item Assume that $\tau=0$, if for a fixed real number $\lambda$ there
are infinitely many positive integers $n$ such that \begin{equation}
n\theta=m+\lambda,\quad m\in\mathbb{Z},\label{eq:503}\end{equation}
 then, \begin{equation}
\frac{S_{n}(x_{n}(z,s);q)(q;q)_{\infty}}{(-z)^{n}q^{n^{2}(1-s)}}=A_{q}\left(\frac{e^{2\lambda\pi i}}{z}\right)+r_{sw}(n|2),\label{eq:504}\end{equation}
 and \begin{equation}
|r_{sw}(n|2)|\le\frac{2A_{q}(-|z|^{-1})}{(q;q)_{\infty}}\left\{ q^{n/2}+\frac{q^{\left\lfloor n/2\right\rfloor ^{2}}}{|z|^{\left\lfloor n/2\right\rfloor }}\right\} \label{eq:505}\end{equation}
 for $n$ is sufficiently large. 
\item Assume that $\tau=0$, if for any fixed real numbers $\beta$ and
$\rho\ge1$ there are infinitely many positive integers $n$ such
that \begin{equation}
n\theta=m+\beta+b_{n},\quad m\in\mathbb{Z},\quad|b_{n}|<\frac{1}{n^{\rho}},\label{eq:506}\end{equation}
then, \begin{equation}
\frac{S_{n}(x_{n}(z,s);q)(q;q)_{\infty}}{(-z)^{n}q^{n^{2}(1-s)}}=A_{q}\left(\frac{e^{2\pi i\beta}}{z}\right)+e_{sw}(n|3),\label{eq:507}\end{equation}
and\begin{equation}
|e_{sw}(n|3)|\le24\exp(q/((1-q)|z|))\left\{ \frac{\log n}{n^{\rho}}+\frac{q^{n/2}}{(1-q)}+\frac{((q^{-1}-1)|z|)^{-k}}{j_{n}!}\right\} \label{eq:508}\end{equation}
 for $n$ is sufficiently large. 
\item Assume that $-1<\tau<0$. If for some fixed real numbers $\lambda$
and $\lambda_{1}$ there are infinite number of positive integers
$n$ such that\begin{equation}
-n\tau=m+\lambda,\quad m\in\mathbb{N},\quad n\theta=m_{1}+\lambda_{1},\quad m_{1}\in\mathbb{Z},\label{eq:509}\end{equation}
 then, \begin{align}
S_{n}(x_{n}(z,s);q) & =\frac{(-z)^{n}q^{n^{2}(1-s)+m\left(2\tau n+m\right)}}{(q;q)_{\infty}^{2}(-ze^{-2n\theta\pi i})^{m}}\left\{ \theta_{4}\left(zq^{2\lambda}e^{-2\pi i\lambda_{1}};q\right)+r_{sw}(n|4)\right\} ,\label{eq:510}\end{align}
and\begin{align}
|r_{sw}(n|4)| & \le12\theta_{3}(|z|q^{2\lambda};q)\left\{ \frac{q^{k_{n}}}{1-q}+|z|^{k_{n}}q^{k_{n}^{2}+2\lambda k_{n}}+\frac{q^{k_{n}^{2}-2\lambda k_{n}}}{|z|^{k_{n}}}\right\} \label{eq:511}\end{align}
 for $n$ sufficiently large. 
\item Assume that $-1<\tau<0$ , if for any real numbers $\beta$, $\lambda$
and $\rho\ge1$ there are infinitely many positive integers $n$ such
that \begin{equation}
n\theta=m_{1}+\beta+b_{n},\quad|b_{n}|<\frac{1}{n^{\rho}},\quad m_{1}\in\mathbb{Z},\label{eq:512}\end{equation}
 and \begin{equation}
-n\tau=m+\lambda,\quad m\in\mathbb{N}.\label{eq:513}\end{equation}
 Then, \begin{align}
S_{n}(x_{n}(z,s);q) & =\frac{(-z)^{n}q^{n^{2}(1-s)+m\left(2\tau n+m\right)}}{(q;q)_{\infty}^{2}(-ze^{-2n\theta\pi i})^{m}}\left\{ \theta_{4}\left(zq^{2\lambda}e^{-2\beta\pi i};q\right)+r_{sw}(n|5)\right\} ,\label{eq:514}\end{align}
 and\begin{align}
|r_{sw}(n|5)| & \le48\theta_{3}(|z|q^{2\lambda};q)\left\{ \frac{q^{j_{n}^{2}-2\lambda j_{n}}}{|z|^{j_{n}}}+|z|^{j_{n}}q^{j_{n}^{2}+2\lambda j_{n}}+\frac{q^{k_{n}}}{1-q}+\frac{\log n}{n^{\rho}}\right\} ,\label{eq:515}\end{align}
 for $n$ sufficiently large. 
\item Assume that $-1<\tau<0$, if for any real numbers $\beta$, $\lambda$
and $\rho\ge1$ there are infinitely many positive integers $n$ such
that \begin{equation}
-n\tau=m+\beta+a_{n},\quad|a_{n}|<\frac{1}{n^{\rho}},\quad m\in\mathbb{N},\label{eq:516}\end{equation}
 and \begin{equation}
n\theta=m_{1}+\lambda,\quad m_{1}\in\mathbb{Z}.\label{eq:517}\end{equation}
 Then,\begin{align}
S_{n}(x_{n}(z,s);q) & =\frac{(-z)^{n}q^{n^{2}(1-s)+m\left(2\tau n+m\right)}}{(q;q)_{\infty}^{2}(-ze^{-2n\theta\pi i})^{m}}\left\{ \theta_{4}\left(zq^{2\beta}e^{-2\lambda\pi i};q\right)+r_{sw}(n|6)\right\} ,\label{eq:518}\end{align}
 and\begin{align}
|r_{sw}(n|6)| & \le36\theta_{3}(|z|q^{2\beta};q)\left\{ \frac{q^{j_{n}^{2}-2(\beta+1)j_{n}}}{|z|^{j_{n}}}+|z|^{j_{n}}q^{j_{n}^{2}+2(\beta-1)j_{n}}+\frac{q^{k_{n}}}{1-q}+\frac{\log n}{n^{\rho}}\right\} ,\label{eq:519}\end{align}
 for $n$ sufficiently large.
\item Assume that $-1<\tau<0$, If for fixed real numbers $\beta_{1}$,
$\beta_{2}$ and $\rho>0$ there are infinitely many positive integers
$n$ such that\begin{equation}
-\tau n=m+\beta_{1}+a_{n},\quad|a_{n}|<\frac{1}{n^{\rho}},\quad m\in\mathbb{N},\label{eq:520}\end{equation}
 and\begin{equation}
n\theta=m_{1}+\beta_{2}+b_{n},\quad|b_{n}|<\frac{1}{n^{\rho}},\quad m_{1}\in\mathbb{Z}.\label{eq:521}\end{equation}
 Then, \begin{align}
S_{n}(x_{n}(z,s);q) & =\frac{(-z)^{n}q^{n^{2}(1-s)+m\left(2\tau n+m\right)}}{(q;q)_{\infty}^{2}(-ze^{-2n\theta\pi i})^{m}}\left\{ \theta_{4}\left(zq^{2\beta_{1}}e^{-2\beta_{2}\pi i};q\right)+r_{sw}(n|7)\right\} ,\label{eq:522}\end{align}
 and\begin{align}
|r_{sw}(n|7)| & \le156\theta_{3}(|z|q^{2\beta_{1}};q)\left\{ \frac{q^{j_{n}^{2}-2(\beta_{1}+1)j_{n}}}{|z|^{j_{n}}}+|z|^{j_{n}}q^{j_{n}^{2}+2(\beta_{1}-1)j_{n}}+\frac{q^{k_{n}}}{1-q}+\frac{\log n}{n^{\rho}}\right\} \label{eq:523}\end{align}
for $n$ sufficiently large.
\end{enumerate}
\end{thm}
In the following corollaries we assume that\begin{equation}
z=e^{2\pi u},\quad u\in\mathbb{R}.\label{eq:524}\end{equation}
 If we let \begin{equation}
q=\exp(-n^{-a}\pi),\quad0<a<\frac{1}{2},\quad n\in\mathbb{N},\label{eq:525}\end{equation}
then we have the following results for Stieltjes-Wigert polynomials:

\begin{cor}
\label{cor:q-stieltjes-wigert-power}
\begin{enumerate}
\item If$\tau>0$, then\begin{align}
 & S_{n}(\exp2\pi(u+(\tau+1)n^{1-a});e^{-\pi n^{-a}})\label{eq:526}\\
 & =\frac{\exp\left(2n\pi u+(2\tau+1)\pi n^{2-a}+\frac{\pi n^{a}}{6}-\frac{\pi n^{-a}}{24}\right)\left\{ 1+\mathcal{O}\left(e^{-4\pi n^{a}}\right)\right\} }{(-1)^{n}\sqrt{2n^{a}}},\nonumber \end{align}
 and\begin{align}
 & s_{n}(\exp2\pi(u+(\tau+1)n^{1-a})\Vert e^{-\pi n^{-a}})\label{eq:527}\\
 & =\frac{(-1)^{n}\exp\left(-\frac{\pi}{n^{a}}(n^{a}u+\tau n)^{2}-\frac{\pi u}{2}\right)\left\{ 1+\mathcal{O}\left(e^{-4\pi n^{a}}\right)\right\} }{\sqrt{2\pi}\exp\left(\frac{\pi(\tau+2)n^{1-a}}{2}+\frac{\pi n^{-a}}{12}-\frac{\pi n^{a}}{12}\right)}\nonumber \end{align}
 as $n\to\infty$.
\item Assume that $-1<\tau<0$ and for some fixed real number $\lambda$,
there are infinitely many positive integers $n$ such that \begin{equation}
-n\tau=m+\lambda,\quad m\in\mathbb{N}.\label{eq:528}\end{equation}
Then,\begin{align}
 & S_{n}\left(\exp2\pi(u+(\tau+1)n^{1-a});e^{-\pi n^{-a}}\right)\label{eq:529}\\
 & =\frac{\exp\left(\pi n^{-a}\left(n^{a}u+(\tau+1)n\right)^{2}\right)}{(-1)^{(1+\tau)n+\lambda}n^{a/2}\exp\left(\frac{\pi n^{-a}}{12}-\frac{\pi n^{a}}{12}\right)}\nonumber \\
 & \times\left\{ \cos\pi(n^{a}u-\lambda)+\mathcal{O}\left(e^{-2\pi n^{a}}\right)\right\} ,\nonumber \end{align}
and\begin{align}
 & s_{n}\left(\exp2\pi(u+(\tau+1)n^{1-a})\Vert e^{-\pi n^{-a}}\right)\label{eq:530}\\
 & =\frac{(-1)^{(1+\tau)n+\lambda}\left\{ \cos\pi(n^{a}u-\lambda)+\mathcal{O}\left(e^{-2\pi n^{a}}\right)\right\} }{\sqrt{\pi}\exp\left(\frac{\pi u}{2}+\frac{(\tau+2)\pi n^{1-a}}{2}+\frac{\pi n^{-a}}{8}\right)}\nonumber \end{align}
 for $n$ sufficiently large. 
\end{enumerate}
\end{cor}
If we choose that \begin{equation}
q=\exp\left(-\frac{\pi}{\gamma\log n}\right),\quad n\ge2,\label{eq:531}\end{equation}
then we have the following:

\begin{cor}
\label{cor:q-stieltjes-wigert-log}
\begin{enumerate}
\item If $\tau>0$, then we have\begin{align}
 & S_{n}\left(\exp2\pi\left(u+\frac{(\tau+1)n}{\gamma\log n}\right);\exp\left(-\frac{\pi}{\gamma\log n}\right)\right)\label{eq:541}\\
 & =\frac{\exp\left(2\pi nu+\frac{(2\tau+1)n^{2}}{\gamma\log n}+\frac{\gamma\pi\log n}{6}-\frac{\pi}{24\gamma\log n}\right)\left\{ 1+\mathcal{O}\left(n^{-4\pi\gamma}\right)\right\} }{(-1)^{n}\sqrt{2\gamma\log n}},\nonumber \end{align}
 and \begin{align}
 & s_{n}\left(\exp2\pi\left(u+\frac{(\tau+1)n}{\gamma\log n}\right)\Vert\exp\left(-\frac{\pi}{\gamma\log n}\right)\right)\label{eq:542}\\
 & =\frac{\exp\left(-\frac{\pi}{\gamma\log n}(\gamma u\log n+\tau n)^{2}-\frac{\pi u}{2}\right)\left\{ 1+\mathcal{O}\left(n^{-4\pi\gamma}\right)\right\} }{(-1)^{n}\sqrt{2\pi}\exp\left(\frac{(\tau+2)n\pi}{2\gamma\log n}+\frac{\pi}{12\gamma\log n}-\frac{\gamma\pi\log n}{12}\right)}\nonumber \end{align}
for $n$ sufficiently large.
\item Assume that $-1<\tau<0$ and for some fixed real number $\lambda$
there are infinitely many positive integers $n$ such that \begin{equation}
-n\tau=m+\lambda,\quad m\in\mathbb{N}.\label{eq:543}\end{equation}
Then,\begin{align}
 & S_{n}\left(\exp2\pi\left(u+\frac{(\tau+1)n}{\gamma\log n}\right);\exp\left(-\frac{\pi}{\gamma\log n}\right)\right)\label{eq:544}\\
 & =\frac{(-1)^{(1+\tau)n+\lambda}n^{\gamma\pi/12}}{\sqrt{\gamma\log n}}\exp\left(\frac{\pi(u\gamma\log n+(\tau+1)n)^{2}}{\gamma\log n}-\frac{\pi}{12\gamma\log n}\right)\nonumber \\
 & \times\left\{ \cos\pi(\gamma u\log n-\lambda)+\mathcal{O}\left(n^{-2\pi\gamma}\right)\right\} ,\nonumber \end{align}
 and\begin{align}
 & s_{n}\left(\exp2\pi\left(u+\frac{(\tau+1)n}{\gamma\log n}\right)\Vert\exp\left(-\frac{\pi}{\gamma\log n}\right)\right)\label{eq:545}\\
 & =\frac{(-1)^{(1+\tau)n+\lambda}\left\{ \cos\pi(\gamma u\log n-\lambda)+\mathcal{O}\left(n^{-2\pi\gamma}\right)\right\} }{\sqrt{\pi}\exp\left(\frac{\pi u}{2}+\frac{\pi(\tau+2)n}{2\gamma\log n}+\frac{\pi}{8\gamma\log n}\right)}\nonumber \end{align}
for $n$ sufficiently large.
\item Assume that $-1<\tau<0$, if for some fixed real numbers $\beta$,
$\rho\ge1$ and $\lambda$ there exist infinitely many positive integers
$n$ such that\begin{equation}
-n\tau=m+\beta+a_{n},\quad|a_{n}|<\frac{1}{n^{\rho}},\quad m\in\mathbb{N},\label{eq:546}\end{equation}
and let \begin{equation}
\gamma=\frac{4\rho}{9\pi},\label{eq:547}\end{equation}
then for each of such $n$, we have\begin{align}
 & S_{n}\left(\exp2\pi\left(u+\frac{(\tau+1)n}{\gamma\log n}\right);\exp\left(-\frac{\pi}{\gamma\log n}\right)\right)\label{eq:548}\\
 & =\frac{(-1)^{n-m}n^{\rho/27}}{\sqrt{\gamma\log n}}\exp\left(\frac{\pi(u\gamma\log n+(\tau+1)n)^{2}}{\gamma\log n}-\frac{\pi}{12\gamma\log n}\right)\nonumber \\
 & \times\left\{ \cos\pi(\gamma u\log n-\beta)+\mathcal{O}\left(n^{-8\rho/9}\log n\right)\right\} ,\nonumber \end{align}
 and\begin{align}
 & s_{n}\left(\exp2\pi\left(u+\frac{(\tau+1)n}{\gamma\log n}\right)\Vert\exp\left(-\frac{\pi}{\gamma\log n}\right)\right)\label{eq:549}\\
 & =\frac{(-1)^{n-m}\left\{ \cos\pi(\gamma u\log n-\beta)+\mathcal{O}\left(n^{-8\rho/9}\log n\right)\right\} }{\sqrt{\pi}\exp\left(\frac{\pi u}{2}+\frac{\pi(\tau+2)n}{2\gamma\log n}+\frac{\pi}{8\gamma\log n}\right)}\nonumber \end{align}
 for $n$ sufficiently large. 
\end{enumerate}
\end{cor}

\subsection{Proofs for Theorem \ref{thm:stieltjes-wigert}}

In the first three proofs we will use the inequalities\begin{equation}
0<\frac{(q;q)_{\infty}}{(q;q)_{n-k}}\le1\label{eq:550}\end{equation}
 for $0\le k\le n$ and \begin{equation}
\left|\frac{(q;q)_{\infty}}{(q;q)_{n-k}}-1\right|\le\frac{2q^{1+n/2}}{1-q}\label{eq:551}\end{equation}
 for $0\le k\le\left\lfloor \frac{n}{2}\right\rfloor -1$ for $n$
sufficiently large, which is a directly from Lemma \ref{lem:1}.

In the last four proofs we have\begin{equation}
-\tau n=m+c_{n},\quad m\in\mathbb{N},\quad n\theta=m_{1}+d_{n},\quad m_{1}\in\mathbb{Z}.\label{eq:552}\end{equation}
 Then,\begin{align}
 & \frac{S_{n}(x_{n}(z,s);q)}{(-z)^{n}q^{n^{2}(1-s)}}=\sum_{k=0}^{n}\frac{q^{k^{2}}e^{2nk\theta\pi i}}{(q;q)_{k}(q;q)_{n-k}}\left(-\frac{q^{2\tau n}}{z}\right)^{k}\label{eq:553}\\
 & =\sum_{k=0}^{m}\frac{q^{k^{2}}e^{2nk\theta\pi i}}{(q;q)_{k}(q;q)_{n-k}}\left(-\frac{q^{2\tau n}}{z}\right)^{k}+\sum_{k=m+1}^{n}\frac{q^{k^{2}}e^{2nk\theta\pi i}}{(q;q)_{k}(q;q)_{n-k}}\left(-\frac{q^{2\tau n}}{z}\right)^{k}\nonumber \\
 & =s_{1}+s_{2}.\nonumber \end{align}
 We reverse the summation order in $s_{1}$ to obtain\begin{equation}
\frac{s_{1}(q;q)_{\infty}^{2}(-ze^{-2n\theta\pi i})^{m}}{q^{m(2\tau n+m)}}=\sum_{k=0}^{m}q^{k^{2}}\left(-zq^{2c_{n}}e^{-2\pi id_{n}}\right)^{k}e(k,n),\label{eq:554}\end{equation}
and\begin{equation}
e(k,n)=\frac{(q;q)_{\infty}^{2}}{(q;q)_{m-k}(q;q)_{n-m+k}}.\label{eq:555}\end{equation}
 It is clear that\begin{equation}
|e(k,n)|\le1\label{eq:556}\end{equation}
 for $0\le k\le m$. Using Lemma \ref{lem:1} to expand and estimate
each term of \begin{eqnarray}
e(k,n) & -1= & \left\{ r_{1}\left(q;m-k\right)+1\right\} \left\{ r_{1}\left(q;n-m+k\right)+1\right\} -1,\label{eq:557}\end{eqnarray}
to get\begin{equation}
|e(k,n)-1|\le\frac{6q^{k_{n}+2}}{1-q}\label{eq:558}\end{equation}
 for $0\le k\le k_{n}-1$ and $n$ sufficiently large.

In sum $s_{2}$ we shift summation index from $k$ to $k+m$ to get
\begin{equation}
\frac{s_{2}(q;q)_{\infty}^{2}(-ze^{-2n\theta\pi i})^{m}}{q^{m(2\tau n+m)}}=\sum_{k=1}^{n-m}q^{k^{2}}\left(-z^{-1}q^{-2c_{n}}e^{2\pi id_{n}}\right)^{k}f(k,n),\label{eq:559}\end{equation}
 and\begin{equation}
f(k,n)=\frac{(q;q)_{\infty}^{2}}{(q;q)_{m+k}(q;q)_{n-m-k}}\label{eq:560}\end{equation}
 for $1\le k\le n-m$. Hence\begin{equation}
|f(k,n)|\le1\label{eq:561}\end{equation}
 for $1\le k\le n-m$. Expand and apply Lemma \ref{lem:1} to each
terms of \begin{eqnarray}
f(k,n) & -1= & \left\{ r_{1}\left(q;m+k\right)+1\right\} \left\{ r_{1}\left(q;n-m-k\right)+1\right\} -1\label{eq:562}\end{eqnarray}
 to obtain \begin{equation}
|f(k,n)-1|\le\frac{6q^{k_{n}+2}}{1-q}\label{eq:563}\end{equation}
 for $1\le k\le k_{n}-1$ and $n$ sufficiently large.

The rest are similar to the corresponding proofs for the Ismail-Masson
polynomials.

\subsection{Proof for Corollary \ref{cor:q-stieltjes-wigert-power} }

In this case we have\begin{equation}
x_{n}(z,s)=\exp2\pi(u+(\tau+1)n^{1-a})\label{eq:564}\end{equation}
and\begin{align}
w_{sw}(x_{n}(z,s)) & =\frac{n^{a/2}}{\pi\sqrt{2}}\exp\left(-2\pi n^{-a}(n^{a}u+(\tau+1)n+\frac{1}{4})^{2}\right).\label{eq:565}\end{align}
 Clearly,\begin{align}
(q;q)_{\infty} & =\sqrt{2n^{a}}\exp\left(\frac{\pi n^{-a}}{24}-\frac{\pi n^{a}}{6}\right)\left\{ 1+\mathcal{O}\left(e^{-4\pi n^{a}}\right)\right\} ,\label{eq:566}\\
\frac{1}{(q;q)_{\infty}} & =\frac{\exp\left(\frac{\pi n^{a}}{6}-\frac{\pi n^{-a}}{24}\right)\left\{ 1+\mathcal{O}\left(e^{-4\pi n^{a}}\right)\right\} }{\sqrt{2n^{a}}},\label{eq:567}\\
q^{n}(q;q)_{n} & =\sqrt{2n^{a}}\exp\left(-\pi n^{1-a}+\frac{\pi n^{-a}}{24}-\frac{\pi n^{a}}{6}\right)\left\{ 1+\mathcal{O}\left(e^{-4\pi n^{a}}\right)\right\} ,\label{eq:568}\end{align}
then,

\begin{align}
 & \sqrt{q^{n}(q;q)_{n}w_{sw}(x_{n}(z,s))}\label{eq:569}\\
 & =\sqrt{\frac{n^{a}}{\pi}}\exp\left(\frac{\pi n^{-a}}{48}-\frac{\pi n^{a}}{12}-\frac{\pi n^{1-a}}{2}\right)\nonumber \\
 & \times\exp\left(-\pi n^{-a}(n^{a}u+(\tau+1)n+\frac{1}{4})^{2}\right)\left\{ 1+\mathcal{O}\left(e^{-4\pi n^{a}}\right)\right\} \nonumber \end{align}
as $n\to\infty$. 

Notice that\begin{align}
\theta_{4}\left(zq^{2\lambda};q\right) & =\theta_{4}\left(ui-i\lambda n^{-a}\mid n^{-a}i\right)\label{eq:570}\\
 & =n^{a/2}\exp\left(\frac{\pi(n^{a}u-\lambda)^{2}}{n^{a}}\right)\theta_{2}\left(n^{a}u-\lambda\mid n^{a}i\right)\nonumber \\
 & =2n^{a/2}\exp\left(\frac{\pi(n^{a}u-\lambda)^{2}}{n^{a}}-\frac{\pi n^{a}}{4}\right)\nonumber \\
 & \times\cos\pi(n^{a}u-\lambda)\left\{ 1+\mathcal{O}\left(e^{-2\pi n^{a}}\right)\right\} ,\nonumber \end{align}
 and\begin{eqnarray}
\theta_{3}\left(zq^{2\lambda};q\right) & = & \theta_{3}\left(ui-i\lambda n^{-a}\mid n^{-a}i\right)\label{eq:571}\\
 & = & n^{a/2}\exp\left(\frac{\pi(n^{a}u-\lambda)^{2}}{n^{a}}\right)\theta_{3}\left(n^{a}u-\lambda\mid n^{a}i\right)\nonumber \\
 & = & 2n^{a/2}\exp\left(\frac{\pi(n^{a}u-\lambda)^{2}}{n^{a}}\right)\left\{ 1+\mathcal{O}\left(e^{-\pi n^{a}}\right)\right\} \nonumber \end{eqnarray}
 as $n\to\infty$. The proof for case 1 follows from \eqref{eq:501}
and \eqref{eq:502}, while the proof for case 2 comes from \eqref{eq:510}
and \eqref{eq:511}.

\subsection{Proof for Corollary \ref{cor:q-stieltjes-wigert-log} }

In this case we have\begin{equation}
x_{n}(z,s)=\exp2\pi\left(u+\frac{(\tau+1)n}{\gamma\log n}\right),\label{eq:572}\end{equation}
 and\begin{equation}
w_{sw}(x_{n}(z,s))=\frac{1}{\pi}\sqrt{\frac{\gamma\log n}{2}}\exp\left(-\frac{2\pi}{\gamma\log n}(\gamma u\log n+(\tau+1)n+\frac{1}{4})^{2}\right).\label{eq:573}\end{equation}
 It is clear that\begin{eqnarray}
(q;q)_{\infty} & = & \sqrt{2\gamma\log n}\exp\left(\frac{\pi}{24\gamma\log n}-\frac{\gamma\pi\log n}{6}\right)\left\{ 1+\mathcal{O}\left(n^{-4\pi\gamma}\right)\right\} ,\label{eq:574}\\
\frac{1}{(q;q)_{\infty}} & = & \frac{1}{\sqrt{2\gamma\log n}}\exp\left(\frac{\gamma\pi\log n}{6}-\frac{\pi}{24\gamma\log n}\right)\left\{ 1+\mathcal{O}\left(n^{-4\pi\gamma}\right)\right\} \label{eq:575}\end{eqnarray}
 as $n\to\infty$.\begin{align}
 & \sqrt{q^{n}(q;q)_{n}w_{sw}(x_{n}(z,s))}\label{eq:576}\\
 & =\sqrt{\frac{\gamma\log n}{\pi}}\exp\left(\frac{\pi}{48\gamma\log n}-\frac{\pi n}{2\gamma\log n}-\frac{\gamma\pi\log n}{12}\right)\nonumber \\
 & \times\exp\left(-\frac{\pi}{\gamma\log n}(\gamma u\log n+(\tau+1)n+\frac{1}{4})^{2}\right)\left\{ 1+\mathcal{O}\left(n^{-4\pi\gamma}\right)\right\} \nonumber \end{align}
 as $n\to\infty$.\begin{align}
\theta_{4}(zq^{2\lambda};q) & =\theta_{4}\left(ui-\frac{\lambda i}{\gamma\log n}\mid\frac{i}{\gamma\log n}\right)\label{eq:577}\\
 & =\sqrt{\gamma\log n}\exp\frac{\pi(\gamma u\log n-\lambda)^{2}}{\gamma\log n}\theta_{2}\left(\gamma u\log n-\lambda\mid i\gamma\log n\right)\nonumber \\
 & =\sqrt[24]{\frac{(\gamma\log n)^{2}}{n^{\pi\gamma}}}\exp\frac{\pi(\gamma u\log n-\lambda)^{2}}{\gamma\log n}\cos\pi(\gamma u\log n-\lambda)\left\{ 1+\mathcal{O}\left(n^{-2\pi\gamma}\right)\right\} ,\nonumber \end{align}
 and\begin{align}
\theta_{3}(zq^{2\lambda};q) & =\theta_{3}\left(ui-\frac{\lambda i}{\gamma\log n}\mid\frac{i}{\gamma\log n}\right)\label{eq:578}\\
 & =\sqrt{\gamma\log n}\exp\frac{\pi(\gamma u\log n-\lambda)^{2}}{\gamma\log n}\theta_{3}\left(\gamma u\log n-\lambda\mid i\gamma\log n\right)\nonumber \\
 & =\sqrt{\gamma\log n}\exp\frac{\pi(\gamma u\log n-\lambda)^{2}}{\gamma\log n}\left\{ 1+\mathcal{O}\left(n^{-\pi\gamma}\right)\right\} ,\nonumber \end{align}
as $n\to\infty$. The proofs for three cases follow from the formulas
\eqref{eq:501}, \eqref{eq:510} and \eqref{eq:518} respectively.

\section{$q$-Laguerre Orthogonal Polynomials }

The $q$-Laguerre orthogonal polynomials $\left\{ L_{n}^{(\alpha)}(x;q)\right\} _{n=0}^{\infty}$
are defined as \cite{Ismail2} \begin{equation}
L_{n}^{(\alpha)}(x;q)=\sum_{k=0}^{n}\frac{q^{k^{2}+\alpha k}(-x)^{k}(q^{\alpha+1};q)_{n}}{(q;q)_{k}(q,q^{\alpha+1};q)_{n-k}},\label{eq:579}\end{equation}
for $\alpha>-1$. For any nonzero complex number $z$, let \begin{equation}
s=2\tau+2+i\frac{2\theta\pi}{\log q},\quad x_{n}(z,s,\alpha)=zq^{-ns-\alpha},\quad\tau,\theta\in\mathbb{R},\label{eq:580}\end{equation}
reverse the summation order to get\begin{equation}
\frac{L_{n}^{(\alpha)}(x_{n}(z,s,\alpha);q)}{(-z)^{n}q^{n^{2}(1-s)}}=\sum_{k=0}^{n}\frac{q^{k^{2}}e^{2nk\theta\pi i}}{(q;q)_{k}(q;q)_{n-k}}\frac{(q^{\alpha+1};q)_{n}}{(q^{\alpha+1};q)_{n-k}}\left(-\frac{q^{2\tau n}}{z}\right)^{k},\label{eq:581}\end{equation}
 then\begin{eqnarray}
\left|L_{n}^{(\alpha)}(x_{n}(z,s,\alpha);q)\right| & \le & \frac{|z|^{n}\sum_{k=0}^{\infty}\frac{q^{k^{2}}\left(|z|^{-1}q^{2\tau n}\right)^{k}}{(q;q)_{k}}}{(q;q)_{\infty}q^{n^{2}(2\tau+1)}}\le\frac{|z|^{n}A_{q}\left(-|z|^{-1}q^{2\tau n}\right)}{(q;q)_{\infty}q^{n^{2}(2\tau+1)}},\label{eq:582}\end{eqnarray}
 or\begin{equation}
\left|L_{n}^{(\alpha)}(x_{n}(z,s,\alpha);q)\right|\le\frac{(-\sqrt{q};q)_{\infty}}{(q;q)_{\infty}}\frac{|z|^{(2\tau+1)n}\exp\left(-\frac{\log^{2}|z|}{2\log q}\right)}{q^{(2\tau^{2}+2\tau+1)n^{2}}}.\label{eq:583}\end{equation}

\subsection{Asymptotic Formulas $q$-Laguerre Polynomials}

\begin{thm}
\label{thm:q-laguerre}Assume that $\alpha>-1$, given any complex
numbers $z\neq0$, let $s$ and $x_{n}(z,s,\alpha)$ be defined as
in \eqref{eq:580} and let \begin{equation}
j_{n}=\left\lfloor \frac{q^{4}\log n}{\log q^{-1}}\right\rfloor ,\quad k_{n}=\min\left\{ \left\lfloor \frac{(1+\tau)n}{2}\right\rfloor ,\left\lfloor \frac{-\tau n}{2}\right\rfloor \right\} ,\label{eq:584}\end{equation}

we have the following results for $q$-Laguerre polynomials:
\begin{enumerate}
\item Assume that $\tau>0$, we have \begin{equation}
\frac{L_{n}^{(\alpha)}(x_{n}(z,s,\alpha);q)(q;q)_{\infty}}{(-z)^{n}q^{n^{2}(1-s)}}=1+r_{ql}(n|1)\label{eq:585}\end{equation}
 with\begin{equation}
\left|r_{ql}(n|1)\right|\le\frac{q^{2\tau n+1}\exp\left(q^{2\tau n+1}/((1-q)|z|)\right)}{(1-q)|z|}.\label{eq:586}\end{equation}

\item Assume that $\tau=0$, if for some fixed real number $\lambda$ there
are infinitely positive integers $n$ such that \begin{equation}
n\theta=m+\lambda,\quad m\in\mathbb{Z},\label{eq:587}\end{equation}
then, \begin{equation}
\frac{L_{n}^{(\alpha)}(x_{n}(z,s,\alpha);q)(q;q)_{\infty}}{(-z)^{n}q^{n^{2}(1-s)}}=A_{q}\left(\frac{e^{2\lambda\pi i}}{z}\right)+r_{ql}(n|2),\label{eq:588}\end{equation}
 and\begin{align}
|r_{ql}(n|2)| & \le\frac{14A_{q}(-|z|^{-1})}{(q;q)_{\infty}}\left\{ q^{n/2}+\frac{q^{\left\lfloor n/2\right\rfloor ^{2}}}{|z|^{\left\lfloor n/2\right\rfloor }}\right\} \label{eq:589}\end{align}
 for $n$ sufficiently large.
\item Assume that $\tau=0$, if for some fixed real numbers $\beta$ and
$\rho\ge1$ there are infinitely many positive integers $n$ such
that\begin{equation}
n\theta=m+\beta+b_{n}\quad|b_{n}|<\frac{1}{n^{\rho}},\quad m\in\mathbb{Z},\label{eq:560}\end{equation}
then, \begin{equation}
\frac{L_{n}^{(\alpha)}(x_{n}(z,s,\alpha);q)(q;q)_{\infty}}{(-z)^{n}q^{n^{2}(1-s)}}=A_{q}\left(\frac{e^{2\pi i\beta}}{z}\right)+e_{ql}(n|3),\label{eq:561}\end{equation}
and\begin{align}
|e_{ql}(n|3)| & \le24\exp\left(\frac{q}{(1-q)z}\right)\left\{ \frac{\log n}{n^{\rho}}+\frac{q^{n/2}}{1-q}+\frac{\left((q^{-1}-1)|z|\right)^{-j_{n}}}{j_{n}!}\right\} \label{eq:562}\end{align}
for $n$ sufficiently large.
\item Assume that $-1<\tau<0$, if for some fixed real numbers $\lambda$
and $\lambda_{1}$ there are infinite number of positive integers
$n$ such that\begin{equation}
-n\tau=m+\lambda,\quad m\in\mathbb{N},\quad n\theta=m_{1}+\lambda_{1},\quad m_{1}\in\mathbb{Z},\label{eq:563}\end{equation}
 then, \begin{align}
L_{n}^{(\alpha)}(x_{n}(z,s,\alpha);q) & =\frac{(-z)^{n}q^{n^{2}(1-s)+m\left(2\tau n+m\right)}}{(q;q)_{\infty}^{2}\left(-ze^{-2n\theta\pi i}\right)^{m}}\left\{ \theta_{4}\left(zq^{2\lambda}e^{-2\pi i\lambda_{1}};q\right)+r_{ql}(n|4)\right\} ,\label{eq:564}\end{align}
and\begin{align}
|r_{ql}(n|4)| & \le60\theta_{3}\left(|z|q^{2\lambda};q\right)\left\{ \frac{q^{k_{n}}}{1-q}+|z|^{-k_{n}}q^{k_{n}^{2}-2\lambda k_{n}}+|z|^{k_{n}}q^{k_{n}^{2}+2\lambda k_{n}}\right\} \label{eq:565}\end{align}
 for $n$ sufficiently large.
\item Assume that $-1<\tau<0$ , if for some fixed real numbers $\beta$,
$\lambda$ and $\rho\ge1$ there are infinitely many positive integers
$n$ such that \begin{equation}
n\theta=m_{1}+\beta+b_{n},\quad|b_{n}|<\frac{1}{n^{\rho}},\quad m_{1}\in\mathbb{Z},\quad-n\tau=m+\lambda,\quad m\in\mathbb{N},\label{eq:566}\end{equation}
then, \begin{align}
L_{n}^{(\alpha)}(x_{n}(z,s,\alpha);q) & =\frac{(-z)^{n}q^{n^{2}(1-s)+m\left(2\tau n+m\right)}}{(q;q)_{\infty}^{2}(-ze^{-2n\theta\pi i})^{m}}\left\{ \theta_{4}\left(zq^{2\lambda}e^{-2\beta\pi i};q\right)+r_{ql}(n|5)\right\} ,\label{eq:567}\end{align}
 and\begin{align}
|r_{ql}(n|5)| & \le60\theta_{3}(|z|q^{2\lambda};q)\left\{ |z|^{j_{n}}q^{j_{n}^{2}+2\lambda j_{n}}+|z|^{-j_{n}}q^{j_{n}^{2}-2\lambda j_{n}}+\frac{q^{k_{n}}}{1-q}+\frac{\log n}{n^{\rho}}\right\} \label{eq:568}\end{align}
 for $n$ sufficiently large.
\item Assume that $-1<\tau<0$, if for any fixed real numbers $\beta$,$\lambda$
and $\rho\ge1$ there are infinitely many positive integers $n$ such
that \begin{equation}
-n\tau=m+\beta+a_{n},\quad|a_{n}|<\frac{1}{n^{\rho}},\quad m\in\mathbb{N},\quad n\theta=m_{1}+\lambda,\quad m_{1}\in\mathbb{Z},\label{eq:569}\end{equation}
then,\begin{align}
L_{n}^{(\alpha)}(x_{n}(z,s,\alpha);q) & =\frac{(-z)^{n}q^{n^{2}(1-s)+m\left(2\tau n+m\right)}}{(q;q)_{\infty}^{2}(-ze^{-2n\theta\pi i})^{m}}\left\{ \theta_{4}\left(zq^{2\beta}e^{-2\lambda\pi i};q\right)+r_{ql}(n|6)\right\} ,\label{eq:570}\end{align}
and\begin{align}
|r_{ql}(n|6)| & \le180\theta_{3}(|z|q^{2\beta};q)\left\{ |z|^{-j_{n}}q^{j_{n}^{2}-2(\beta+1)j_{n}}+|z|^{j_{n}}q^{j_{n}^{2}+2(\beta-1)j_{n}}+\frac{q^{k_{n}}}{1-q}+\frac{\log n}{n^{\rho}}\right\} \label{eq:571}\end{align}
 for $n$ sufficiently large.
\item Assume that $-1<\tau<0$, if for some fixed real numbers $\beta_{1}$,
$\beta_{2}$ and $\rho>0$ there are infinitely many positive integers
$n$ such that\begin{equation}
-\tau n=m+\beta_{1}+a_{n},\quad|a_{n}|<\frac{1}{n^{\rho}},\quad m\in\mathbb{N},\label{eq:572}\end{equation}
 and\begin{equation}
n\theta=m_{1}+\beta_{2}+b_{n},\quad|b_{n}|<\frac{1}{n^{\rho}},\quad m_{1}\in\mathbb{Z},\label{eq:573}\end{equation}
 then, \begin{align}
L_{n}^{(\alpha)}(x_{n}(z,s,\alpha);q) & =\frac{(-z)^{n}q^{n^{2}(1-s)+m\left(2\tau n+m\right)}}{(q;q)_{\infty}^{2}(-ze^{-2n\theta\pi i})^{m}}\left\{ \theta_{4}\left(zq^{2\beta_{1}}e^{-2\beta_{2}\pi i};q\right)+r_{ql}(n|7)\right\} ,\label{eq:574}\end{align}
 and\begin{align}
|r_{ql}(n|7)| & \le180\left\{ |z|^{-j_{n}}q^{j_{n}^{2}-2(\beta_{1}+1)j_{n}}+|z|^{j_{n}}q^{j_{n}^{2}+2(\beta_{1}-1)j_{n}}+\frac{q^{k_{n}}}{1-q}+\frac{\log n}{n^{\rho}}\right\} \label{eq:575}\end{align}
 for $n$ sufficiently large. 
\end{enumerate}
\end{thm}
In the following corollaries we assume that\begin{equation}
z=e^{2\pi u},\quad u\in\mathbb{R}.\label{eq:576}\end{equation}
 If we let \begin{equation}
q=\exp(-n^{-a}\pi),\quad0<a<\frac{1}{2},\quad n\in\mathbb{N},\label{eq:577}\end{equation}
then we have the following results for $q$-Laguerre orthogonal polynomials:

\begin{cor}
\label{cor:q-laguerre-power}
\begin{enumerate}
\item If$\tau>0$, then\begin{align}
 & L_{n}^{(\alpha)}\left(\exp2\pi(u+(\tau+1)n^{1-a}+\alpha n^{-a}/2);\exp(-n^{-a}\pi)\right)\label{eq:578}\\
 & =\frac{\exp2n\pi\left(u+(\tau+\frac{1}{2})n^{1-a}\right)\left\{ 1+\mathcal{O}(e^{-4\pi n^{a}})\right\} }{(-1)^{n}\sqrt{2n^{a}}\exp\left(-\frac{\pi}{24}(n^{-a}-4n^{a})\right)}\nonumber \end{align}
 as $n\to\infty$.
\item Assume that $-1<\tau<0$ and for some fixed real number $\lambda$
there are infinitely many positive integers $n$ such that \begin{equation}
-n\tau=m+\lambda,\quad m\in\mathbb{N}.\label{eq:579}\end{equation}
Then, \begin{align}
 & L_{n}^{(\alpha)}(\exp2\pi(u+(\tau+1)n^{1-a}+\alpha n^{-a}/2);\exp(-n^{-a}\pi))\label{eq:580}\\
 & =\frac{\exp\left(\pi n^{-a}(n^{a}u+(\tau+1)n)^{2}\right)\left\{ \cos\pi(n^{a}u-\lambda)+\mathcal{O}\left(e^{-2\pi n^{a}}\right)\right\} }{(-1)^{(1+\tau)n+\lambda}\sqrt{n^{a}}\exp\left(\frac{\pi n^{-a}}{24}+\frac{\pi n^{a}}{12}\right)}\nonumber \end{align}
 for $n$ sufficiently large. 
\end{enumerate}
\end{cor}
If we choose \begin{equation}
q=\exp\left(-\frac{\pi}{\gamma\log n}\right),\quad n\ge2,\label{eq:581}\end{equation}
then we have the following:

\begin{cor}
\label{cor:q-laguerre-log}
\begin{enumerate}
\item If $\tau>0$, then we have\begin{align}
 & L_{n}^{(\alpha)}\left(\exp\frac{2\pi(\gamma u\log n+(\tau+1)n+\alpha/2)}{\gamma\log n};\exp\left(-\frac{\pi}{\gamma\log n}\right)\right)\label{eq:582}\\
 & =\frac{n^{\pi\gamma/6}\exp\left(2n\pi u+\frac{(2\tau+1)n^{2}\pi}{\gamma\log n}-\frac{\pi}{24\gamma\log n}\right)}{(-1)^{n}\sqrt{2\gamma\log n}}\left\{ 1+\mathcal{O}\left(n^{-4\pi\gamma}\right)\right\} \nonumber \end{align}
 for $n$ sufficiently large.
\item Assume that $-1<\tau<0$ and for some fixed real number $\lambda$
there are infinitely many positive integers $n$ such that \begin{equation}
-n\tau=m+\lambda,\quad m\in\mathbb{N},\label{eq:583}\end{equation}
then,\begin{align}
 & L_{n}^{(\alpha)}\left(\exp\left(\frac{2\pi(\gamma u\log n+(\tau+1)n+\alpha/2)}{\gamma\log n}\right);\exp\left(-\frac{\pi}{\gamma\log n}\right)\right)\label{eq:584}\\
 & =\frac{n^{\pi\gamma/12}\exp\left(\frac{\pi(u\gamma\log n+(\tau+1)n)^{2}}{\gamma\log n}-\frac{\pi}{12\gamma\log n}\right)\left\{ \cos\pi(\gamma u\log n-\lambda)+\mathcal{O}\left(n^{-2\pi\gamma}\right)\right\} }{(-1)^{(1+\tau)n+\lambda}\sqrt{\gamma\log n}}\nonumber \end{align}
 for $n$ sufficiently large.
\item Assume that $-1<\tau<0$, if for some fixed real numbers $\beta$,
$\rho\ge1$ and $\lambda$ there exist infinitely many positive integers
$n$ such that\begin{equation}
-n\tau=m+\beta+a_{n},\quad|a_{n}|<\frac{1}{n^{\rho}},\quad m\in\mathbb{N},\label{eq:585}\end{equation}
and let \begin{equation}
\gamma=\frac{4\rho}{9\pi},\label{eq:586}\end{equation}
then for each of such $n$, we have\begin{align}
 & L_{n}^{(\alpha)}\left(\exp\left(\frac{2\pi(\gamma u\log n+(\tau+1)n+\alpha/2)}{\gamma\log n}\right);\exp\left(-\frac{\pi}{\gamma\log n}\right)\right)\label{eq:587}\\
 & =\frac{\exp\left(\frac{\pi(u\gamma\log n+(\tau+1)n)^{2}}{\gamma\log n}-\frac{\pi}{12\gamma\log n}\right)}{(-1)^{n-m}n^{-\rho/27}\sqrt{\gamma\log n}}\left\{ \cos\pi(\gamma u\log n-\beta)+\mathcal{O}\left(\frac{\log n}{n^{8\rho/9}}\right)\right\} \nonumber \end{align}
 for $n$ sufficiently large. 
\end{enumerate}
\end{cor}

\subsection{Proofs for $q$-Laguerre Polynomial Case}

In the first three proofs we will use the inequalities\begin{equation}
0<\frac{(q;q)_{\infty}}{(q;q)_{n-k}}\frac{(q^{\alpha+1};q)_{n}}{(q^{\alpha+1};q)_{n-k}}\le1\label{eq:588}\end{equation}
 for $0\le k\le n$ and \begin{equation}
\left|\frac{(q;q)_{\infty}}{(q;q)_{n-k}}\frac{(q^{\alpha+1};q)_{n}}{(q^{\alpha+1};q)_{n-k}}-1\right|\le\frac{14q^{n/2}}{1-q}\label{eq:589}\end{equation}
 for $0\le k\le\left\lfloor \frac{n}{2}\right\rfloor $-1 and $n$
sufficiently large. This can be seen by expanding \begin{align}
\frac{(q;q)_{\infty}}{(q;q)_{n-k}}\frac{(q^{\alpha+1};q)_{n}}{(q^{\alpha+1};q)_{n-k}} & -1=\left\{ r_{1}(q;n-k)+1\right\} \left\{ r_{2}(q^{\alpha+1};n)+1\right\} \left\{ r_{1}(q^{\alpha+1};n-k)+1\right\} -1\label{eq:590}\end{align}
 and estimating each term by Lemma \ref{lem:1}.

Assume that\begin{equation}
-\tau n=m+c_{n},\quad m\in\mathbb{N},\quad n\theta=m_{1}+d_{n},\quad m_{1}\in\mathbb{Z},\label{eq:591}\end{equation}
then,\begin{align}
 & \frac{L_{n}^{(\alpha)}(x_{n}(z,s,\alpha);q)}{(-z)^{n}q^{n^{2}(1-s)}}=\sum_{k=0}^{n}\frac{q^{k^{2}}e^{2nk\theta\pi i}}{(q;q)_{k}(q;q)_{n-k}}\frac{(q^{\alpha+1};q)_{n}}{(q^{\alpha+1};q)_{n-k}}\left(-\frac{q^{2\tau n}}{z}\right)^{k}\label{eq:592}\\
 & =\sum_{k=0}^{m}\frac{q^{k^{2}}e^{2nk\theta\pi i}}{(q;q)_{k}(q;q)_{n-k}}\frac{(q^{\alpha+1};q)_{n}}{(q^{\alpha+1};q)_{n-k}}\left(-\frac{q^{2\tau n}}{z}\right)^{k}\nonumber \\
 & +\sum_{k=m+1}^{n}\frac{q^{k^{2}}e^{2nk\theta\pi i}}{(q;q)_{k}(q;q)_{n-k}}\frac{(q^{\alpha+1};q)_{n}}{(q^{\alpha+1};q)_{n-k}}\left(-\frac{q^{2\tau n}}{z}\right)^{k}=s_{1}+s_{2}.\nonumber \end{align}
 We reverse the summation order in $s_{1}$ to obtain\begin{equation}
\frac{s_{1}(q;q)_{\infty}^{2}(-ze^{-2n\theta\pi i})^{m}}{q^{m(2\tau n+m)}}=\sum_{k=0}^{m}q^{k^{2}}\left(-zq^{2c_{n}}e^{-2\pi id_{n}}\right)^{k}e(k,n),\label{eq:593}\end{equation}
 and\begin{equation}
e(k,n)=\frac{(q;q)_{\infty}^{2}(q^{\alpha+1};q)_{n}}{(q;q)_{m-k}(q;q)_{n-m+k}(q^{\alpha+1};q)_{n-m+k}}.\label{eq:594}\end{equation}
 Clearly,\begin{equation}
|e(k,n)|\le1\label{eq:595}\end{equation}
 for $0\le k\le m$. Expand and estimate each term of\begin{align}
e(k,n) & -1=\left\{ r_{2}\left(q^{\alpha+1};n\right)+1\right\} \left\{ r_{1}\left(q;m-k\right)+1\right\} \left\{ r_{1}\left(q;n-m+k\right)+1\right\} \left\{ r_{1}\left(q^{\alpha+1};n-m+k\right)+1\right\} -1\label{eq:596}\end{align}
 by Lemma \ref{lem:1} to get\begin{equation}
|e(k,n)-1|\le\frac{30q^{k_{n}+1}}{1-q}\label{eq:597}\end{equation}
 for $0\le k\le k_{n}-1$.

In sum $s_{2}$ we shift summation from $k$ to $k+m$ to obtain \begin{equation}
\frac{s_{2}(q;q)_{\infty}^{2}(-ze^{-2n\theta\pi i})^{m}}{q^{m(2\tau n+m)}}=\sum_{k=1}^{n-m}q^{k^{2}}\left(-z^{-1}q^{-2c_{n}}e^{2\pi id_{n}}\right)^{k}f(k,n),\label{eq:598}\end{equation}
 and\begin{equation}
f(k,n)=\frac{(q;q)_{\infty}^{2}(q^{\alpha+1};q)_{n}}{(q;q)_{m+k}(q;q)_{n-m-k}(q^{\alpha+1};q)_{n-m-k}}\label{eq:599}\end{equation}
 for $1\le k\le n-m$. Hence\begin{equation}
|f(k,n)|\le1\label{eq:600}\end{equation}
 for $1\le k\le n-m$. Apply Lemma \ref{lem:1} to each term of \begin{align}
f(k,n) & -1=\left\{ r_{2}\left(q^{\alpha+1};n\right)+1\right\} \left\{ r_{1}\left(q;m+k\right)+1\right\} \left\{ r_{1}\left(q;n-m-k\right)+1\right\} \left\{ r_{1}\left(q^{\alpha+1};n-m-k\right)+1\right\} -1,\label{eq:601}\end{align}
we obtain \begin{equation}
|f(k,n)-1|\le\frac{30q^{k_{n}+1}}{1-q}\label{eq:602}\end{equation}
 for $1\le k\le k_{n}-1$.

The rest of the proofs are very similar to the corresponding proofs
for Ismail-Masson polynomials and we will not repeat them here.

\end{document}